\documentclass{birkmult}
\usepackage[dvips]{epsfig,color}
\usepackage{pstricks,pst-text,pst-grad,pst-char}
\usepackage{amsmath}
\usepackage{amsfonts}
\usepackage[psamsfonts]{amssymb}
\usepackage{amsthm}
\usepackage[mathscr]{eucal}
\usepackage{exscale}
\usepackage{graphics,graphicx}

\usepackage[russian,english]{babel}

\newtheorem{theorem}{THEOREM}[section]
\newtheorem{lemma}{LEMMA}[section]
\newtheorem{definition}{DEFINITION}[section]
\newtheorem{remark}{REMARK}[section]
\newtheorem{corollary}{COROLLARY}[section]
\def\theequation{\arabic{section}.\arabic{equation}}

\DeclareMathOperator*{\Vee}{\bigvee}
\DeclareMathOperator*{\Vvee}{\Large\text{ \(\bigvee\)}}
\DeclareMathOperator*{\Pprod}{\Huge\text{\(\prod\)}}
\newcommand{\Curvearrowleft}%
    {\LARGE\text{\(\curvearrowleft\)}}%
\DeclareSymbolFont{rsfs}{U}{rsfs}{m}{n}
\DeclareSymbolFontAlphabet{\mathscrn}{rsfs}
 \DeclareFixedFont{\Azh}{OT1}{ptm}{b}{n}{2.2ex}
\newcommand{\Bbo}%
   {\pscharpath[fillstyle=none,linewidth=0.4pt,linecolor=black]%
      {\text{\Azh {1}}}}
 \newcommand{\Bbe}%
   {\pscharpath[fillstyle=none,linewidth=0.2pt,linecolor=black]%
     {\text{\large\Azh {\(\boldsymbol\varepsilon\)}}}}%
\title{The Schur Algorithm in Terms of System \linebreak Realizations}

\author{Bernd Fritzsche}
  \address{ Mathematisches Institut\\
   Universit\"at Leipzig\\
   D-04009, Leipzig,\\ Germany}
\email{bernd.fritzsche@mathematik.uni-leipzig.de}
\author{Victor Katsnelson}
\address{Department of Mathematics\\
 the Weizmann Institute\\
  Rehovot 76100\\
   Israel}
\email{victor.katsnelson@weizmann.ac.il,
victorkatsnelson@gmail.com}%
\author{Bernd Kirstein}
  \address{ Mathematisches Institut\\
   Universit\"at Leipzig\\
   D-04009, Leipzig,\\ Germany}
\email{bernd.kirstein@mathematik.uni-leipzig.de}
\date{}
\dedicatory{Dedicated to Moshe Liv\v{s}ic: A great man,
              thinker, philosopher and mathematician.}

\begin{document}
\begin{abstract}
The main goal of this paper is to demonstrate the usefulness of
certain ideas from System Theory in the study of problems from
complex analysis. With this paper, we also aim to encourage analysts, who
might not be familiar with System Theory, colligations
or operator models to take a closer look at these topics.
For this reason, we present a short introduction to the necessary
background. The method of system realizations of analytic
functions often provides new insights into and interpretations of
results relating to the objects under consideration. In this paper
we will use a well-studied topic from classical analysis as an
example. More precisely, we will look at the classical
Schur algorithm from the perspective of System Theory.
We will confine our considerations to rational inner
functions. This will allow us to avoid questions
involving limits and will enable us to concentrate on
the algebraic aspects of the problem at hand.
Given a non-negative integer \(n\), we describe all
system realizations of a given rational inner
function of degree \(n\) in terms of an appropriately
constructed equivalence relation in the set of all unitary
\( (n + 1) \times (n + 1) \)-matrices.
The concept of Redheffer coupling of colligations gives us
the possibility to choose a particular representative
from each equivalence class. The Schur algorithm
for a rational inner function is, consequently, described
in terms of the state space representation.
\end{abstract}
\maketitle
 NOTATION:
\begin{enumerate}
\item[] \(\mathbb{T}\) is the unit circle in the complex plane:
\(\mathbb{T}=\{t\in\mathbb{C}:\ |t|=1\}\)
\item[]  \(\mathbb{D}\) \ \ is the unit disc in the complex plane:
\(\mathbb{D}=\{z\in\mathbb{C}:\ |z|<1\}\)
\item[]  \(\mathbb{D}^-\) \ is the exterior of the unit circle:
\(\mathbb{D}^-=\{z: 1<|z|\leq\infty\}\).
\item[]  \(\mathfrak{M}_{p\times{}q}\) is the set of all \(p\times{}q\)
(\(p\) rows, \(q\) columns) matrices with complex entries\,.
\item[]  \(I_n\) - the identity \(n\times{}n\) matrix.
\end{enumerate}

\vspace{2.0ex}%
 \centerline{\Large\textbf{Table of Contents.}}
\begin{enumerate}
\setcounter{enumi}{-1}
\item
Introduction.
\item
Rational Inner Functions.
\item
The Schur Algorithm.
\item
The System Representation of a Rational Inner Function.
\item
Coupled Systems and The Schur Transfopmation: Input-Output
Mappings.
\item
The Redheffer Coupling of Unitary Colligations.
\item
The Inverse Schur Transformation and Redheffer Couplings of
Colligations.
\item
One Step of the Schur Algorithm, Expressed in the Language of
Colligations.
\item
Hessenberg Matrices. The Householder Algorithm.
\item
The Schur Algorithm in Terms of System Representations.
\item
An Expression for the Colligation Matrix in Terms of the Schur
Parameters.
\item
On Work Related to System Theoretic Interpretations
of the Schur Algorithm
\item
Appendix: System Realizations of Inner Rational Functions.
\end{enumerate}

\setcounter{section}{-1}
\section{Introduction\label{INT}}

Up until the 1960s System Theory suggested that a system be
considered only in terms of its input and output. A system was
treated as a `black box' with input and output terminals.
Associated with each system was an `input-output' mapping,
considered to be of primary importance to the theory at the time.
This approach, however, did not take the internal state of the
system into account. It is to be assumed that an input signal
will, in some way, influence the internal state of a system.
Nevertheless, there was little discussion of the relationship
between input and the inner state of a system until the
introduction of State Space System Theory. This theory not only
incorporated input and output spaces, serving, respectively, as
`domains' for input and output `signals', but also a
`state space'. This state space was introduced to describe the
interior state of the system.

State Space System Theory (both the linear and general non-linear
variations of the theory) was developed in the early 1960s. Two
names closely associated with the early development of this theory
are those of R. Kalman and M.S. Livshitz. Kalman's first
publications pertaining to State Space System Theory include
\cite{Kal1, Kal2, Kal3}. The monograph \cite{KFA} summarizes these
papers, among others. R. Kalman's approach to State Space System
Theory was from the perspective of Control Theory. This approach
suggested that the questions of a system's controllability and
observability be given the most attention. Control Theory does
not, however, put much emphasis on energy relations and, as a
result, Kalman's work does not address the subject of energy
balance relations (Kalman's approach to System Theory was
abstract. He develops the theory over arbitrary fields, not
specifically over the field of complex numbers). In Kalman's
theory, one first starts from the input-output behavior (i.e.
transfer function) and then constructs the state operator. In
Livshitz's theory, the reverse approach is used: The
characteristic function (which is the analogue of the transfer
function) is produced from the main operator (which is the
analogue of the state operator). Kalman's theory is mainly
finite-dimensional and affine, whereas Livshitz's theory is mainly
infinite-dimensional and metric. It took some decades before the connections
between these two theories were discovered in the 1970s. Among
others, Dewilde \cite{Dew1, Dew2} and Helton \cite{He1,He2,He3} produced
much of the work leading to this discovery. The connections
between the two approaches were made explicit in the monograph
\cite{BGK}.

M. S. Livshitz, a pioneer in the theory of non-self-adjoint
operators, chose to approach State Space Theory from the
perspective of Operator Theory. For a particular class of
non-self-adjoint operators, Livshitz was able to associate each
operator of this class with an analytic function in the upper
half-plane or unit disc. These analytic functions were called
`characteristic functions'. Livshitz was, furthermore, able to
determine a correspondence between the invariant subspaces of a
linear operator and the factors of its characteristic function
(See \cite{Liv3} and references within \cite{Liv3}). Using the
framework provided for by these results for characteristic
functions, Livshitz constructed triangular models of
non-self-adjoint operators (Triangular models were later partially
supplanted by functional models. See \cite{SzNFo}). Following
this, Livshitz focused on questions in both mathematics and
physics. Oscillation and wave propagation problems in linear
isolated systems are related to self-adjoint operators. In the
mid-1950s M.S. Livshitz began to look for a physical example to
which his theory of non-self-adjoint operators could be applied.
This lead him to consider a number of concrete linear systems.
These systems were not isolated systems, but were such that they
allowed for the exchange of energy with the `external world'. The
model of the dynamical behavior of a system of this type makes use
of an operator and this `principal' operator is, in general,
non-self-adjoint. The energy exchange of the system is reflected
in the non-self-adjointness of the operator. Livshitz worked on
problems involving the scattering of elementary particles (See
\cite{Liv4}, \cite{Liv5}, \cite{BrLi}), problems in electrical
networks (See \cite{LiFl}) and questions dealing with wave
propagation in wave-guides (See \cite{Liv6}). It was at this
juncture that the notion of an `operator colligation' (also common
are the terms `operator node' and `operator cluster') was
introduced to provide further clarity. An operator colligation
consists of the aforementioned `principal' operator, but also
`channel' spaces and `channel' operators, of which the latter two
objects describe the non-self-adjointness of the `principal'
operator. The introduction of this concept allowed a
characteristic function to be associated with an operator
colligation, as opposed to its respective `principal' operator
(See \cite{BrLi}, \cite{Br}, \cite{LiYa} and references therein).
At much the same time, the concept of an `open system' was then
being established (What Livshitz then referred to as an `open
system' was, in essence, what is now known as a stationary linear
dynamical system). Livshitz first introduced the notion of an
`open system' in his influential paper, \cite{Liv9} (See
Definition 1 on p. 1002 of the original Russian paper [Liv7] or
p. ??? of the English translation in the present volume). To each
system there is an associated colligation and in \cite{Liv7} it is
shown that a system's transfer operator coincides with the
characteristic function of the system's colligation. \cite{Liv8}
introduces the operation of coupling open systems as well as the
concept of closing coupling channels. \cite{Liv8} furthermore
introduces the `kymological resolution' of an open system, i.e.
the resolution of this system into a chain of simpler coupled open
systems. These simpler systems correspond to the invariant
subspaces of the `inner-state' operator of the original open
system. To emphasize that the notion of an open system is closely
related to oscillations and to wave-propagation processes,
Livshitz uses the terminology `kymological', `kymmer' and
`kymmery', derived from the Greek word
`\(\kappa\upsilon\mu\alpha\)', meaning `wave'.
 We quote from page 15 of the English translation of \cite{Liv9}
  and mention
 that: "the appropriate representation of an open system, transforming a
  known input into a known output, depends on which are known and unknown
  variables, so that the concept of an open system is `physico-logical' rather
  than purely physical in nature."

The relevant theory of open systems and operator colligations, as
developed by Livshitz and other mathematicians, is presented in
the monographs \cite{Liv9}, \cite{LiYa} and \cite{Br}. Chapter 2
of the monograph \cite{Liv9} deals with the details of the
kymological resolution of open systems (a concept of which much
use is made in the following). A detailed presentation of
Scattering Theory for linear stationary dynamical systems (with an
emphasis on applications to the Wave Equation in \(\mathbb{R}^n)\)
can be found in \cite{LaPhi}.

General State Space System Theory, as developed by R. Kalman and
M. S. Livshitz provides us with the proper setting and the
necessary language for the further study of physical systems and
various aspects of Control Theory. Despite the fact that
State Space System Theory does not immediately lead to a solution
of the initial physical or control problem, it does lead to some
interesting related questions (mostly analytic). It should,
furthermore, be noted that general State Space Theory's importance
extends beyond its significance within Control Theory and when
applied to physical systems. M. S. Livshitz was very likely the
first to understand that this theory had wide-reaching
applications within mathematics, e.g. in Complex Analysis.

Analytic functions can be represented or specified in many ways,
e.g. as Taylor-series, by decomposition into continuous fractions,
or via representations as Cauchy or Fourier integrals. In the
early half of the 1970s an additional method for representing an
analytic function was introduced, namely the method of `system
realization'. This theory has its origins in Synthesis Theory for
linear electrical networks, the theory of linear control systems
and the theory of operator colligations (and associated
characteristic functions). M. S. Livshitz established the Theory
of System Realizations and L. A. Sakhnovich, a former Ph.D.
student of Livshitz's, later made further important progress in
the theory (See \cite{Sakh1} and also \cite{Sakh2} for a more
detailed presentation of these results). L. A. Sakhnovich studied
the spectral factorization of a given rational matrix-function
\(R\), where both \(R\) and the inverse function \(R^{-1}\) are
transfer functions corresponding to linear systems (operator
colligations). Unfortunately, the paper \cite{Sakh1} did not
garner the attention it deserved at the time. L. A Sakhnovich's
factorization theorem is a predecessor to a fundamental result due
to Bart/Gohberg/Kaashoek/van Dooren \cite{BGKV}, which was
remembered as Theorem 2 in the Editorial Introduction to
\cite{CWHF}, where one can also find a detailed account of the
history of the state space factorization theorem.

Our goal is not to provide a comprehensive survey of the
history of System Theory, so that we have focused on the
period leading up to the mid-1970s (with particular emphasis
on the contributions of M. S. Livshitz and his co-workers).
His work on open systems was unkown in the western
world until his monograph \cite{Liv9} was translated in 1973.
His fundamental papers \cite{Liv7} and \cite{Liv8} remained
untranslated up until this memorial volume.

The subsequent
development of the Theory of System Realizations is generally
associated with the name I. Gohberg, who produced and inspired
much in the way of new work and results for this theory and its
applications. As a result, the theory experienced a period of
accelerated growth, beginning in the late 1970s. Published in
1979, the monograph \cite{BGK}\footnote{ N.b. There is now an
extended version of this monograph. See \cite{BGKR}. } dealt with
general factorizations of a rational matrix-functions as well as with
the Wiener-Hopf factorization of rational matrix-function,
where, in both cases, this function is a transfer function for a linear
system (operator colligation).

I. Gohberg and his co-workers have shown that State Space Theory has
a much wider range and goes far beyond System Theory and the theory
of operator colligations. We list a few topics to which State Space
Theory can be applied:
\begin{enumerate}
    \item
        Methods of factorization of matrix- and operator-valued functions;
        solutions of Wiener-Hopf and singular integral equations.

    \item
        Interpolation in the complex plane and generalizations.

    \item
        Limit formulas of Akhiezer/Kac/Widom type.

    \item
        Projection methods, Bezoutiants, resultants.

    \item
        Inverse problems.
\end{enumerate}

The monograph \cite{BGR} offers a
detailed discussion of interpolation problems and many other questions.
Matrix-function factorization is a tool applied in discussions of
many other problems as well, e.g. in the theory of inverse
problems for differential equations and also in prediction theory
for stationary stochastic processes. If a matrix-function is
rational, then this factorization can be attained using system
realizations. These system realizations, in turn, play a certain
role in the solution of the original problem (See, for example,
\cite{AG}). The Theory of Isoprincipal Deformations of Rational
Matrix-Functions (which is, in particular, a useful tool for
investigating rational solutions of Schlesinger systems) is
formulated in terms of the Theory of System Realizations (See
\cite{KaVo1} and \cite{KaVo2}. For our purposes, the theory
developed in \cite{Ka} is most relevant). The current state of
System Theory, as a branch of pure mathematics, is presented in
\cite{Nik}.

In the present paper we show how the Schur algorithm for
contractive holomorphic functions in the unit disc can be
described in terms of system realizations. In the following, we
consider only rational inner functions, which allows us to avoid
questions involving limits and enables us to concentrate on the
algebraic aspects of the problem at hand. At first glance the
formulas here presented might seem rather complicated and, to some
degree, less than intuitive. This is, however, from the
perspective of System Theory, not the case. The aforementioned
formulas serve as the function-theoretical counterpart to
Livshitz's kymological resolution as applied to the system
(represented by the original inner function) corresponding to the
cascade coupling, i.e. the Redheffer coupling, of open systems.
The elementary open systems, which make up this cascade (or chain)
correspond to the steps of the Schur algorithm.

This paper is organized as follows. In Section \ref{RIF},
we state some facts relating to rational inner functions.
In Section \ref{SA}, we discuss some aspects of the
classical Schur algorithm. Section \ref{SRIF} is
devoted to a short introduction to operator colligations
and their characteristic functions, where particular
attention is paid to finite-dimensional unitary colligations.
The characteristic functions of finite-dimensional unitary
colligations are shown to be rational inner matrix-functions
(See Theorem \ref{PCFUC}). Theorem \ref{ReThU} shows that
an arbitrary rational inner matrix-function can, on the other
hand, be realized as a characteristic function of a
finite-dimensional minimal unitary colligation.
The scalar rational inner functions of degree \(n\) are just
the finite Blaschke products of \(n\) elementary Blaschke factors.
The essential facts on the realization of scalar inner rational
functions of degree \(n\) as characteristic functions of
minimal unitary colligations are summarized in Theorem \ref{SyR}.
These minimal unitary colligations can be equivalently described
by equivalence classes of minimal unitary
\( (n + 1) \times (n + 1) \)-matrices. A proof for Theorem \ref{SyR}
can be found in the Appendix at the end of the paper.

The main objective of this paper can be described as follows. The
application of the Schur algorithm to a given rational inner function
\(s(z)\) of degree \(n\) produces a sequence \,
\(s_k(z) \, , \  k = 0, \, 1, \, \ldots \, , n \) \,
of rational inner functions with \(s_0(z) = s(z)\) and
\(\deg{s_i(z)} = n - k\). In particular, the function \(s_n(z)\)
is constant with unimodular value.
In Section \ref{SRIF}, it will be shown that each of the functions
\(s_k(z)\) admits a system representation
\begin{equation*}
    s_k(z)  =   A_k +   z B_k
                        \left(
                            I   -   z   D_k
                        \right)^{-1}
                        C_k
\end{equation*}
in terms of the blocks of some minimal unitary matrix
\(U_k   \in \mathfrak{M}_{(n+1)\times(n+1)}\),
\begin{equation*}
    U_k         =   \begin{pmatrix}
                            A_k&B_k\\
                            C_k&D_k
                        \end{pmatrix}   .
\end{equation*}
We assume that \(U_0\) is given. The goal is to recursively produce
the sequence matrices \(U_k\). In other words, the steps of the
Schur algorithm have to be described in terms of the state space
representation. Since the unitary matrices \(U_k\) are defined only
up to an equivalence relation, we have to find corresponding
operations for the arithmetic of these equivalence classes.

In Section \ref{SASRio} we discuss the means by which the linear-fractional
transformation associated with the Schur algorithm can be described
in terms of the input-output mapping of linear systems. The
Redheffer coupling of linear systems will be introduced as a useful
tool in these considerations.

In Section \ref{RCOC}, the Redheffer coupling of linear systems will
be translated into the language of unitary systems.

In Section \ref{SASRoc}, we apply the concept of Redheffer couplings
of colligations to the linear-fractional transformation associated
with the inverse of the Schur algorithm. In so doing, we will
describe the `degrees of freedom' of unitary equivalence. A closer look shows
us that amongst all the unitary matrices which realize the
desired system realization, there are some distinguished by the fact
that they are, in a sense, associated with the concept of Redheffer coupling.

In Section \ref{SSA}, the basic step of the Schur algorithm will
be described in the language of colligations. This requires that
we solve a particular equation for unitary matrices, suggested
by the results of Section \ref{SASRoc}. The solution to
this matrix equation is given in Theorem \ref{fnt}.
Together with Lemma \ref{RMCF}, Theorem \ref{StScAl}
describes the basic step of the Schur algorithm in terms
of system representations.

The investigations of Section \ref{HaHe} show that a certain
normalization procedure has to be performed at every step of the
Schur algorithm if the Schur algorithm is to be dealt with
in the language of system realizations. We consider the degrees
of freedom for this normalization procedure. It turns out that we can use
these degrees of freedom to make the normalization procedure a
one-time-procedure, so that it might be dealt with during
preprocessing for further step-by-step recurrence. A
one-time-normalization of this kind is related to the reduction
of the `initial' colligation matrix to the lower Hessenberg matrix.

In Section \ref{SchA}, we will be well-positioned to present the
Schur algorithm in terms of unitary colligations representing
the appropriate functions.

In Section \ref{ESchP} we express the colligation matrix in terms of
the Schur parameters.

In the final section (Section \ref{RelWork}) we discuss some connections
between the present work and other work relating to the Schur
algorithm as expressed in terms of system realizations. In particular,
we discuss the results presented in Alpay/Azizov/Dijksma/Langer \cite{AADL}
and Killip/Nenciu \cite{KiNe}.

\section{Rational Inner Functions\label{RIF}}
\setcounter{equation}{0}

We say that a function \(\,s:\mathbb{D}\to\mathbb{C}\), where \(s\) is
holomorphic in \(\mathbb{D}\), is \linebreak
\textsf{contractive} if \[|s(z)|\leq 1 \text{ for every } z\in\mathbb{D}.\]
A contractive function \(s\) is called an \textsf{inner} function if
\[|s(t)|= 1 \text{ for every } t\in\mathbb{T}.\]

In the following we consider \textsf{rational inner functions},
so that \(s(t)\) is defined for every \(t\in\mathbb{T}\).

A rational function is representable as a quotient of
irreducible polynomials and we call the order of the highest-degree polynomial
\textsf{the degree of the rational function}.

If a rational function \(s\) is an inner function, then the degree of its
numerator and the degree of its denominator are equal.

An inner rational function \(s\) is representable as a
\textsf{finite Blaschke product}, i.e. in the form
\begin{equation}
\label{BP}
s(z)=c\prod\limits_{1\leq{}k\leq{}n}\frac{z_k-z}{1-z\overline{z}_k} .
\end{equation}
\(z_1,\,\dots\,,\,z_n\) are points in \(\mathbb{D}\), or, in other words,
complex numbers satisfying the condition
\begin{equation}
\label{zbp} |z_1|<1,\,\dots\,,\,|z_n|<1,
\end{equation}
\(c\) is a unimodular complex number, i.e.
\begin{equation}
\label{abp} |c|=1.
\end{equation}
Conversely, given complex numbers \(z_1,\,\dots\,,\,z_n\) and
\(c\) satisfying the conditions \eqref{zbp} and \eqref{abp},
respectively, the function \(s\) in \eqref{BP} is an inner rational function
of degree \(n\).

The number \(c\) and the \textit{set} \(\{z_1,\,\dots\,,\,z_n\}\)
are uniquely defined by the inner function \(s\) (the
\textit{sequence} of numbers \((z_1,\,\dots\,,\,z_n)\)) up to permutation.

The notions of contractive and inner functions can also be defined for matrix-functions:

We say that a matrix-function \(S:
\mathbb{D}\to\mathfrak{M}_{p\times{}p}\)\,, where \(S\) is
holomorphic in \(\mathbb{D}\), is \textsf{contractive} if
\[I_p-S^{\ast}(z)S(z)\geq{}0 \text{ for every } z\in\mathbb{D}\,.\]
A contractive matrix-function \(S: \mathbb{D}\to\mathfrak{M}_{p\times{}p}\)\,,
is called an \textsf{inner} function if%
\footnote{%
For a contractive holomorphic function \(S\) in
\(\mathbb{D}\), the boundary values \(S(t)\stackrel{\text{\tiny
def}}{=} \lim\limits_{r\to{}1-0}S(rt)\) exist for almost every
 \(t\in\mathbb{T}\) (with respect to the Lebesgue measure).}%
\[I_p-S^{\ast}(t)S(t)={}0 \text{ for almost every } t\in\mathbb{T}\,.\]

\section{The Schur Algorithm\label{SA}}
In this section, we present a short introduction to the classical
Shur algorithm, which orginated in Issai Schur's renowned paper,
\cite{Sch}. In so doing, we will mainly emphasize those aspects of
the Schur algorithm, which are essential for this paper.
For comprehensive treatments of the Schur algorithm and its
matricial generalizations, we refer the reader to \cite{BFK1},
\cite{BFK2}, \cite{Con2}, \cite{DFK}, \cite{S:Meth} and the
references therein.

\setcounter{equation}{0}%
 Let \(s(z)\) be a contractive holomorphic function in \(\mathbb{D}\) and
 \begin{equation}
 \label{zsp}
        s_0=s(0)    .
 \end{equation}
 Then \(|s_0|\leq 1\), where \(|s_0|=1\) only if \(s(z)\equiv{}s_0\).
  If \(|s_0|<1\), then the function
 \begin{equation}       \label{ie}
        \omega(z)=\frac{1}{z}\frac{s(z)-s_0}{1-s(z)\overline{s_0}}
 \end{equation}
    is well-defined. Moreover, it is contractive holomorphic in \(\mathbb{D}\).
    The function \(s(z)\) can be expressed in terms of these \(\omega(z)\)
    and \(s_0\):
 \begin{equation}        \label{fsa}
        s(z)=\frac{s_0+z\,\omega(z)}{1+z\,\overline{s_0}\,\omega(z)}\,.
 \end{equation}
 If the function \(s(z)\) is an inner function, then \(\omega(z)\)
 is also an inner function. If \(s(z)\) is an inner rational function of degree
 \(n\), then \(\omega(z)\) is an inner rational function of degree
 \(n-1\).

 Conversely, if \(\omega(z)\) is an \textit{arbitrary} contractive
 holomorphic function in \(\mathbb{D}\) and \(s_0\) is an
 \textit{arbitrary} complex number satisfying the condition \(|s_0|<1\),
 then the expression on the right-hand side of \eqref{fsa}
 defines the function \(s(z)\), which is holomorphic and
 contractive in \(\mathbb{D}\). Furthermore, if \(\omega(z)\) is
 an inner function, then \(s(z)\) is an inner function as well.
 \begin{definition}\ \\ %
     \label{DST}%
        \hspace*{2.0ex}\textrm{I.}
            We call the transformation \(s(z) \longmapsto \omega(z)\),
            defined by \eqref{ie}, where \linebreak
            \(s_0=s(0)\), the (direct)
            {\sf Schur transformation}.  \\
        \hspace*{2.0ex}\textrm{II.}
            We call the transformation \ \(\omega(z) \longmapsto{} s(z)\),
            defined by \eqref{fsa}, where \(s_0\) is a given complex number,
            the \textsf{inverse Schur transformation.}
 \end{definition}%

 The correspondence \(s(z)\Longleftrightarrow (s(0),\,\omega(z)\))
 describes the elementary step of the Schur algorithm.

 The Schur algorithm is applied to a holomorphic
 function \(s(z)\), which is contractive in \(\mathbb{D}\).
 This algorithm inductively produces the sequence (finite or infinite)
  of contractive holomorphic functions \(s_k(z)\) in \(\mathbb{D}\)
 and contractive numbers \(s_k=s_k(0)\), \(k=0,\,1,\,2\,,\,\dots\,\,\).
 The algorithm terminates only if \(s(z)\) is a rational inner
 function. Starting from \(s(z)\), we define
 \[s_0(z)=s(z), \qquad  s_0=s_0(0)\,. \]
  If the functions \(s_i(z),\quad i=0,\,1,\,\dots\,,\,k\) are
  already constructed and \linebreak
  \(|s_k(0)|<1\), then we construct the function \(s_{k+1}(z)\) as follows:
 \begin{equation}
     \label{kas}
      s_{k+1}(z)=\frac{1}{z}\,\frac{s_{k}-s_{k}(z)}{1-s_{k}(z)\overline{s_{k}}}\,,\quad
      s_{k+1}=s_{k+1}(0)\,.
  \end{equation}

If \(s(z)\) is not a rational inner function, then the algorithm
does not terminate: On the \(k\)-th step we obtain the function \(s_{k}(z)\),
for which \linebreak \(|s_{k}(0)|<1\),
so that we can construct the function \(s_{k+1}(z)\) and still have \linebreak \(|s_{k+1}(0)|<1\).

If \(s(z)\) is a rational inner function of degree \(n\), then we
can define the functions \(s_i(z)\) for \(i=0,\,1,\,\dots\,,\,n\) such that
 \begin{equation*}
        \deg{}s_i(z)=n-i, \qquad    i=0,\,1,\,\dots\,,\,n\,.
 \end{equation*}
The numbers \(s_i=s_i(0)\) satisfy the conditions
\[|s_i|<1,\,    \qquad  i=0,\,1\,\dots\,,\,n-1\,.\]
However, in this case
\[|s_n|=1, \qquad s_n(z)\equiv s_n.\]
So, for \(k=n\) the numerator and the denominator of the
expression on the right-hand side of \eqref{kas} vanish
identically. The function \(s_{n+1}(z)\) is thus not
defined and the Schur algorithm terminates.\\

The numbers \(s_k=s_k(0)\) are called the \textsf{Schur parameters}
of the function \(s(z)\).\\

 If \(s(z)\) is not an inner rational function,
then the sequence of its Schur parameters is infinite and these
parameters \(s_k\) satisfy the inequality \linebreak \(|s_k|<1\) for all
\(k:\,0\leq k<\infty\). If \(s(z)\) is an inner
rational function \linebreak with \( \deg s(z)=n\), then its Schur parameters
\(s_k\) are defined only for \linebreak \(k=0,\,1,\,\ldots\,,n\) and
\begin{equation}
    \label{Sprif}%
     |s_k|<1,\quad k=0,\,1,\,\dots\,,\,n-1\,, \ \qquad
    |s_n|=1.
\end{equation}
Conversely, given complex numbers \(s_0,\,s_1,\,\dots\,,\,s_n\)
satisfying the conditions \eqref{Sprif}, one can construct the
inner rational function of degree \(n\), having Schur parameters
\(s_0,\,s_1,\,\dots\,,\,s_n\). This function \(s(z)\) can be
constructed inductively: First, we set
\[s_n(z)\equiv{}s_n\,.\]
If the functions \(s_i(z)\) for \(i=n,\,n-1,\,\dots\,,\,k\)  are
already constructed, then we set
\[s_{k-1}(z)=\frac{s_{k-1}+z\,s_k(z)}{1+z\,\overline{s_{k-1}}\,s_k(z)}\,.\]
In the final step we construct the function \(s_0(z)\) and set
 \begin{equation*}
        s(z)=s_0(z)     .
 \end{equation*}

Thus, \textit{there exists a one-to-one correspondence between
rational inner functions of degree \(n\) and sequences of complex
numbers \(\{s_k\}_{0\leq{}k\leq{n}}\) satisfying the conditions}
\eqref{Sprif}.
\section{The System Representation of a Rational Inner Function.\label{SRIF}}
\setcounter{equation}{0}%

Contractive holomorphic functions appear in several roles. In
particular, such functions appear in Operator Theory as the
\emph{characteristic functions of operator colligations}. The
notion of an operator colligation is closely related to that of
a linear stationary dynamical system. There is a correspondence
between the theory of operator colligations and the theory of
linear stationary dynamical systems. The concepts and results of one
theory can be translated into the language of the other.
There are interesting connections to be made between these theories.
Definitions and constructions, which are well-motivated and natural
in the framework of one theory may look  artificial in the
framework of the other. In particular, the notion of the
characteristic function of a colligation and of the
coupling of colligations are more transparent in the language of
System Theory.

In this section, the term `operator' means `continuous linear operator'.

\begin{definition}\label{DOC}
Let \(\mathcal{H},\ \mathcal{E}^{in},\ \mathcal{E}^{out}\) be Hilbert
spaces and \(U\) be an operator:
\begin{equation}
\label{UnCol}
U:\,\mathcal{E}^{in}\oplus{}\mathcal{H}\to{}\mathcal{E}^{out}\oplus{}\mathcal{H}\,,
\end{equation}
Let%
\begin{equation}%
\label{CODB}%
 U=
\begin{bmatrix}
A&B\\[1.0ex]
C&D
\end{bmatrix}
\end{equation}%
 be the block decomposition of the operator \(U\),
corresponding to \eqref{UnCol}:
\begin{equation}%
\label{FWTW}%
A:\,\mathcal{E}^{in}\to{}\mathcal{E}^{out},\
B:\,\mathcal{H}\to{}\mathcal{E}^{out}, \
C:\,\mathcal{E}^{in}\to{}\mathcal{H},
D:\,\mathcal{H}\to{}\mathcal{H}\,.
\end{equation}

The quadruple \((\mathcal{E}^{in},\ \mathcal{E}^{out},\ %
\mathcal{H},\,U)\) is called an \textsf{operator colligation}.

\(\mathcal{E}^{in}\) and  \(\mathcal{E}^{out}\) are, respectively,
the \textsf{input} and \textsf{output spaces} of the colligation.
We call \(\mathcal{H}\) the \textsf{state space} of the colligation and
\(A\) the \textsf{exterior operator}.
We call \(B\) and \(C\) \textsf{channel
operators}, while \(D\) is referred to as the \textsf{principal operator} of
the colligation.
Finally, we call \(U\) the \textsf{colligation operator}.

If the input and the output spaces \(\mathcal{E}^{in}\) and
\(\mathcal{E}^{out}\) coincide:
\(\mathcal{E}^{in}=\mathcal{E}^{out}=\mathcal{E}\), we call the space
\(\mathcal{E}\) the \textsf{exterior space} of the colligation and denote the
colligation by the triple \((\mathcal{E},\,\mathcal{H},\,U)\)
\end{definition}
 \begin{definition}
\label{DeOpCo}%
Let \((\mathcal{E}^{in},\,\mathcal{E}^{out},\,\mathcal{H},\,U)\) be an
operator colligation.

The operator-function
\begin{equation}
\label{CharFunc}%
    S(z)=A+zB(I_{\mathcal{H}}-zD)^{-1}C\,
\end{equation}
is called the \textsf{characteristic function} of the colligation
\((\mathcal{E}^{in},\,\mathcal{E}^{out},\,\mathcal{H},\,U)\).

The function \(S(z)\) is defined for the \(z\in\mathbb{C}\)
where the operator \((I_{\mathcal{H}}-zD)^{-1}\) exists.
The values of \(S\) are operators acting from \(\mathcal{E}^{in}\) into
\(\mathcal{E}^{out}\).
\end{definition}

\begin{remark}
    \label{DDH}
    The function \(S(z)\) is defined and holomorphic in some neighborhood of the point \(z=0\).
    Furthermore, \(S(0)=A\).
\end{remark}

The notion of a colligation's characteristic function draws on
the framework of the theory of linear stationary
dynamical systems (LSDS). (The theory of open systems, in the
terminology of M.S.Liv\^sic).  The theory of LSDS,  which we are
dealing with is not a `black box theory', where only the input
signals, output signals and the mapping `input\,\(\to\)\,output'
are considered. The theory of LSDS also takes 'interior states'
of the system into account. The input and output signals are
described (in the discrete time case, where the
index \(k\) serves as time) by sequences
\(\{\varphi_k\}_{0\leq{}k<\infty}\) and
\(\{\psi_k\}_{0\leq{}k<\infty}\) of vectors belonging to some
Hilbert spaces \(\mathcal{E}^{in}\) and \(\mathcal{E}^{out}\) (the
\textit{input} and the \textit{output} spaces of the system). The
`interior states' are described by vectors \(h\) of a Hilbert
space \(\mathcal{H}\), called \textit{the state space} of the system.

The dynamics of a \textit{linear stationary} system is
described by the linear equations
\begin{equation}
        \label{LiStDS}
    \begin{bmatrix}
            \psi_k\\[0.7ex]
            h_{k+1}
    \end{bmatrix}=
    \begin{bmatrix}
            A&B\\[0.7ex]
            C&D
    \end{bmatrix}
    \begin{bmatrix}
            \varphi_k\\[0.7ex]
            h_{k}
    \end{bmatrix}\,,
    \quad k=0,\,1,\,2\,,\ \ldots\,,
\end{equation}
where the operators \(A,\,B,\,C,\,D\) do not depend on
\(k\) (`time') and are defined in \eqref{FWTW}.

It is natural to consider the four operators \(A,\,B,\,C,\,D\) as
blocks of the `unified' operator, say \(U\) as in \eqref{CODB},
from the space \(\mathcal{E}^{in}\oplus\mathcal{H}\) into the space
\(\mathcal{E}^{out}\oplus\mathcal{H}\). The operator colligation
\((\mathcal{E}^{in},\,\mathcal{E}^{out},\,\mathcal{H},\,U)\) then
corresponds to the LSDS \eqref{LiStDS}, \eqref{FWTW}. Given the
sequence \(\{\varphi_k\}_{0\leq{}k\leq{}m}\) and the initial value
\(h_0\), the system \eqref{LiStDS} uniquely determines the
sequences \(\{\psi_k\}_{0\leq{}k\leq{}m}\) and
\(\{h_k\}_{0\leq{}k\leq{}m+1}\). In the case \(h_0=0\),
\begin{equation}%
\label{RSTD}%
 \psi_0=A\varphi_0,\quad
\psi_m=A\varphi_m+\sum\limits_{1\leq{}k\leq{}m-1}BD^kC\varphi_{m-k-1}\,,
\ m\geq{}1\,.
\end{equation}
The relation \eqref{RSTD} can be considered as the description of
the evolution of the LSDS \eqref{LiStDS} in the \textit{time domain}.
The description of the evolution is, however, especially
transparent in the \textit{frequency domain}. Since the considered
sequences are unilateral, the Fourier transforms of these
sequences are the (formal) power series
\begin{equation}
\label{FoTr}%
 \varphi(z)=\sum\limits_{0\leq{}k<\infty}\varphi_kz^k\,,\ \
 \psi(z)=\sum\limits_{0\leq{}k<\infty}\psi_kz^k\,, \ \
h(z)=\sum\limits_{0\leq{}k<\infty}h_kz^k\,.
\end{equation}
The complex variable \(z\) can be interpreted as the frequency. Under
the extra assumption that \(h_0=0\) we can rewrite \eqref{LiStDS}
in terms of the Fourier representations:
\begin{equation}
\label{LSDSFr}
\begin{bmatrix}
\psi(z)\\[1.0ex]
z^{-1}h(z)
\end{bmatrix}=
\begin{bmatrix}
A&B\\[1.0ex]
C&D
\end{bmatrix}
\begin{bmatrix}
\varphi(z)\\[1.0ex]
h(z)
\end{bmatrix}
\,.
\end{equation}
 From \eqref{LSDSFr} we obtain
\begin{subequations}
\label{LinS}
\begin{align}
    \psi(z)&=A\varphi(z)+Bh(z),\label{LinS1}\\
    h(z)&=z\,(I-zD)^{-1}C\varphi(z)\,.\label{LinS2}
\end{align}
\end{subequations}
Eliminating \(h(z)\), we get
\begin{equation}
    \label{TrF}
     \psi(z)=S(z)\varphi(z)\,,
\end{equation}
where \(S(z)\) is expressed in terms of the matrix \eqref{CODB} as
in \eqref{CharFunc}:
\[S(z)=A+zB(I-zD)^{-1}C\,.\]
 The operator function \(S(z)\) describes the input-output mapping
corresponding to LSDS \eqref{LiStDS}.

\begin{definition}
    \label{TrFu} In the framework of System Theory, the function
    \(S(z)\) in \eqref{TrF} is called the \textsf{transfer matrix}
    of the LSDS \eqref{LiStDS}.
\end{definition}

In the theory of operator colligations the operator function \(S(z)\)
is called the \textit{characteristic function}, while in the theory of LSDS
it is called the \textit{transfer function}. This notion, however,
makes more sense in the theory of LSDS. Along with the
input-output mapping described by the transfer function \(S(z)\),
the input-state mapping:
\begin{equation*}%
 \varphi(z)\to{}h(z) \, ,
 \quad \text{where} \quad
 h(z)=z\,(I-zD)^{-1}C\varphi(z)\,,
\end{equation*}%
is also naturally related to the system \eqref{LiStDS}.

 If the dimensions \( \ \dim \mathcal{E}^{in}\) and \( \ \dim \mathcal{E}^{out}\) of
the input and output spaces are
finite, then, choosing  bases in \(\mathcal{E}^{in}\) and
\(\mathcal{E}^{out}\), we can consider \(S(z)\) as a matrix-valued
function. If, moreover, the dimension \(\dim \mathcal{H}\) of the
state space is finite, then \(S(z)\) is a rational matrix-function.
\begin{definition}
    The colligation
    \((\mathcal{E}^{in},\,\mathcal{E}^{out},\,\mathcal{H},\,U)\) is said to
    be \textsf{finite-dimensional} if \(\ \dim{\mathcal{E}^{in}}<\infty\),
    \(\ \dim{\mathcal{E}^{out}}<\infty \ \) and \(\ \dim{\mathcal{H}}<\infty\).
\end{definition}

The dimension \(\ \dim \mathcal{H}\) of the state space of the
finite-dimensional colligation \((\mathcal{E}^{in},\,\mathcal{E}^{out},\,
\mathcal{H},\,U)\) is related to the degree of its characteristic
function. Here we use the notion of the \emph{McMillan
degree} of a rational matrix-valued function as it is defined in
\cite{McM}. The notion of the degree of a rational matrix-function is
discussed in \cite{DuHa} and \cite{Kal4}. See also \cite{BGK}. In
the case when \(\dim\mathcal{E}=1\), i.e. in the case when the
considered rational function is scalar (or \(\mathbb{C}\)-valued),
the McMillan degree of this function coincides with its `standard' degree.

To precisely formulate how the dimension of the state space \(\mathcal{H}\)
and the degree of the characteristic function \(S(z)\) are related, we need
to introduce the notion of a \textsl{minimal} colligation
\((\mathcal{E}^{in},\,\mathcal{E}^{out},\,\mathcal{H},\,U)\).

\begin{definition}
Let \((\mathcal{E}^{in},\,\mathcal{E}^{out},\,\mathcal{H},\,U)\) be a colligation.
We define the following subspaces of the state space \(\mathcal{H}\):
\begin{equation}
    \label{COS}
    \mathcal{H}^c=\textup{clos}\big(\Vee\limits_{0\leq{}k<\infty}{(D^kC)\mathcal{E}^{in}}\big),\quad
    \mathcal{H}^o=
    \textup{clos}\big(\Vee\limits_{0\leq{}k<\infty}{(D^{\ast{}k}B^\ast)\mathcal{E}^{out}}\big)\,,
\end{equation}
where \(\Vee\limits_kf_k\) denotes the linear hull of the vectors
\(f_k\) and \(\text{clos}(M)\) denotes the closure of the set \(M\).

 The subspaces \(\mathcal{H}^c\) and \(\mathcal{H}^o\)
are, respectively, called the \textsf{controllability} and
\textsf{observability subspaces} of the colligation \((\mathcal{E}^{in},\,\mathcal{E}^{out},\,\mathcal{H},\,U)\).
\end{definition}
\begin{remark}
If the state space \(\mathcal{H}\) is finite-dimensional, say \linebreak
\(\dim{}\mathcal{H}=n<\infty\),
then it is enough to restrict our considerations in \eqref{COS} to the linear hull of the vectors
\((D^kC)\mathcal{E}^{in}\) and \((D^{\ast{}k}B^\ast)\mathcal{E}^{out}\) with \(k<{}n\).
In this case there is no need to make use of the closure in \eqref{COS}.
\end{remark}
\begin{definition}
\label{DeCoOb}%
 We say that a colligation \((\mathcal{E}^{in},\,\mathcal{E}^{out},\,\mathcal{H},\,U)\)
is \textsf{controllable} if \(\mathcal{H}^c=\mathcal{H}\) and \textsf{observable} if
\(\mathcal{H}^o=\mathcal{H}\).

We say that a colligation is \textsf{simple} if the sum of the
controllability and the observability subspaces is dense in the
state space, i.e. if
\[\textup{clos}\big(\mathcal{H}^c+\mathcal{H}^o\big)=\mathcal{H}\,.
\]

We say that a colligation
\((\mathcal{E}^{in},\,\mathcal{E}^{out},\,\mathcal{H},\,U)\) is
\textsf{minimal} if it is both controllable and observable, i.e. if
\begin{equation*}
    \mathcal{H}^c=\mathcal{H}   \quad and \quad \mathcal{H}^o=\mathcal{H}\,.
\end{equation*}
\end{definition}

\begin{theorem}
\label{dEd}%
 Let \((\mathcal{E}^{in},\,\mathcal{E}^{out},\,\mathcal{H},\,U)\) be a finite-dimensional
 colligation and let \(S(z)\) be the characteristic function of this colligation.

\(S(z)\) is then a rational matrix-function, which is holomorphic at \(z=0\) and
such that
\begin{equation}
    \label{IDeDi} \deg S\leq \dim \mathcal{H}
\end{equation}
 Equality holds in \eqref{IDeDi} if and only if the colligation \((\mathcal{E}^{in},\,\mathcal{E}^{out},\,\mathcal{H},\,U)\)
is minimal.
\end{theorem}

\begin{theorem}
\label{RealTh}
Let \(\mathcal{E}_1\) and \(\mathcal{E}_2\) be finite-dimensional spaces
and let \(S(z)\) be a rational function, whose values are operators acting
from \(\mathcal{E}_1\) to \(\mathcal{E}_2\) and which is holomorphic at the point
\(z=0\).

 There then exists a finite-dimensional minimal operator colligation\\
\((\mathcal{E}_1^{in},\,\mathcal{E}_1^{out},\,\mathcal{H},\,U)\),
\eqref{UnCol}\,-\,\eqref{CODB}\,-\,\eqref{FWTW}, with
 \(\mathcal{E}^{in}=\mathcal{E}_1\) and \(\mathcal{E}^{out}=\mathcal{E}_2\),
whose characteristic function \(S_U(z)=A+zB(I-zD)^{-1}C\) coincides with
 the original function \(S(z)\). In other words, \(S\) can be expressed
 in the form \eqref{CharFunc}.

\end{theorem}

\begin{definition}
\label{DeStSpR}
The representation of a given function \(S(z)\) as a characteristic function
of an operator colligation is called the \textsf{state space representation of \(S(z)\)}
or the \textsf{state space realization of \(S(z)\)}. If the representative operator
colligation is minimal, then we say that the state space realization of \(S(z)\) is
\textsf{minimal}.
\end{definition}

Let us discuss the uniqueness of the state space representation.

\begin{definition}
Let \((\mathcal{E}_1^{in},\,\mathcal{E}_1^{out},\,\mathcal{H}_1,\,U_1)\) and
\((\mathcal{E}_2^{in},\,\mathcal{E}_2^{out},\,\mathcal{H}_2,\,U_2)\) be two operator
colligations:
\begin{equation}
\label{DCO}%
 U_1=
\begin{bmatrix}
A_1&B_1\\[1.0ex]
C_1&D_1
\end{bmatrix},\qquad
U_2=
\begin{bmatrix}
A_2&B_2\\[1.0ex]
C_2&D_2
\end{bmatrix}\,,
\end{equation}
where
\begin{multline}
\label{fwtw}
A_i:\,\mathcal{E}_i^{in}\to\mathcal{E}_i^{out}\,,\ \ B_i:\,\mathcal{H}_i\to\mathcal{E}_i^{out},\ \  %
C_i:\,\mathcal{E}_i^{in}\to\mathcal{H}_i,\ \ D_i:\,\mathcal{H}_i\to\mathcal{H}_i,\\ %
     i=1,\,2\,.
\end{multline}
We consider the colligations \((\mathcal{E}_1^{in},\,\mathcal{E}_1^{out},\,
\mathcal{H}_1,\,U_1)\) and \((\mathcal{E}_2^{in},\,\mathcal{E}_2^{out},\,
\mathcal{H}_2,\,U_2)\) to be \textsf{equivalent}
if \emph{invertible} operators \(E^{in},\,E^{out}\) and \(V\):
\begin{equation}%
\label{FWTWb}
E^{in}:\mathcal{E}_2^{in}\to\mathcal{E}_1^{in},\quad \ E^{out}:\mathcal{E}_2^{out}\to\mathcal{E}_1^{in}, \quad \
V:\mathcal{H}_2=\mathcal{H}_1,
\end{equation}
exist, such that the intertwining relation
\begin{equation}
\label{LinEq}
\begin{bmatrix}
E^{out}&0\\[1.0ex]
0&V
\end{bmatrix}
\begin{bmatrix}
A_2&B_2\\[1.0ex]
C_2&D_2
\end{bmatrix}=
\begin{bmatrix}
A_1&B_1\\[1.0ex]
C_1&D_1
\end{bmatrix}
\begin{bmatrix}
E^{in}&0\\[1.0ex]
0&V
\end{bmatrix}
\end{equation}
holds.
\end{definition}

Clearly, given two
equivalent operator colligations, one of these colligations
is controllable, observable, simple or minimal if and only if the other
colligation possesses the same respective property.\\

The following result is evident:
\begin{theorem}%
  \label{ChSO}%
Let  \((\mathcal{E}_1^{in},\,\mathcal{E}_1^{out},\,
\mathcal{H}_1,\,U_1)\) and
 \((\mathcal{E}_2^{in},\,\mathcal{E}_2^{out},\,
\mathcal{H}_2,\,U_2)\) be operator colligations.
Assume that these colligations
 are equivalent, i.e. that the intertwining relation \eqref{LinEq} holds
with some invertible operators \(E^{in},\,E^{out}\) and \(V\).

 Then the characteristic functions
\(S_1(z)\) and \(S_2(z)\) of these colligations,
\begin{equation}%
\label{CFs}
S_i(z)=A_i+zB_i(I-zD_i)^{-1}C_i,\qquad \ i=1,\,2,
\end{equation}
satisfy the intertwining relation:
\begin{equation}%
\label{Intw}%
E^{out}S_2(z)=S_1(z)E^{in}.
\end{equation}
for all \(z\) where \(S_1\) and \(S_2\) are defined.
\end{theorem}
Under the extra assumptions that the colligations
are minimal and finite-dimensional we can show that for
Theorem \ref{ChSO} the converse assertion also holds.
\begin{theorem}
 \label{ChSOIn}
 Let
\((\mathcal{E}_1^{in},\,\mathcal{E}_1^{out},\,
\mathcal{H}_1,\,U_1)\) and
 \((\mathcal{E}_2^{in},\,\mathcal{E}_2^{out},\,
\mathcal{H}_2,\,U_2)\) be \emph{finite-dimensional}
operator colligations.  Let \(S_1(z)\) and \(S_2(z)\),
\eqref{CFs}, be the characteristic functions of these colligations.
We make the following assumptions:
\begin{enumerate}
\item[\textup{1.}]
The functions
\(S_1(z)\) and \(S_2(z)\) satisfy
the intertwining relation \eqref{Intw} for  all \(z\)
small enough,
where \(E^{in}:\mathcal{E}_2^{in}\to{}\mathcal{E}_1^{in}\) and
\(E^{out}:\mathcal{E}_2^{out}\to{}\mathcal{E}_1^{out}\) are some invertible operators.
\item[\textup{2.}]
The colligations \((\mathcal{E}_1^{in},\,\mathcal{E}_1^{out},\,\mathcal{H}_1,\,U_1)\) and
\((\mathcal{E}_2^{in},\,\mathcal{E}_2^{out},\, \mathcal{H}_2,\,U_2)\) are minimal.
\end{enumerate}

These colligations are then equivalent, i.e. there exists an invertible
operator \(V:\mathcal{H}_2\to\mathcal{H}_1\) such that the intertwining
relation \eqref{LinEq} holds.
\end{theorem}

Up to this point, we have not taken advantage of any scalar products that may
be defined in the input, output and state spaces. From this point forward,
we will focus more on these scalar products and the benefits they bring
when we have them at our disposal. In what follows, we consider rational
\emph{inner} functions. Operator colligations representing such functions
are \emph{unitary, finite-dimensional} operator colligations.\\

For convenience, we recall the definition of a unitary operator:

\textit{Let \(\mathcal{L}_1\) and
\(\mathcal{L}_2\) be Hilbert spaces and \(T:\mathcal{L}_1\to\mathcal{L}_2\)
be an operator. We say that \(T\) is \textsf{unitary} if it satisfies the
following two conditions:\\
\hspace*{2.0ex} a) \(T\) preserves the scalar product, i.e.
\begin{equation*}%
\langle\,Tx,\,Ty\rangle_{\mathcal{L}_2}=\langle\,x,\,y\rangle_{\mathcal{L}_1}
\quad \forall x\in\mathcal{L}_1,\,y\in\mathcal{L}_1\,.
\end{equation*}
\hspace*{2.0ex} b) \(T\) maps \(\mathcal{L}_1\) \textsf{onto} \(\mathcal{L}_2\),
i.e. \(T\) is invertible.}\\

The unitarity property of a linear operator \(T\) can also be characterized as
follows:
\begin{equation*}%
T^{\ast}T=I_{\mathcal{L}_1}, \qquad TT^{\ast}=I_{\mathcal{L}_2}\,.
\end{equation*}
\begin{definition}
\label{DefUOpC}
Let \((\mathcal{E}^{in},\,\mathcal{E}^{out},\,\mathcal{H},\,U)\),
\eqref{CODB}\,-\,\eqref{FWTW}, be an
operator colligation.
We call
\((\mathcal{E}^{in},\,\mathcal{E}^{out},\,\mathcal{H},\,U)\)
a \textsf{unitary colligation} if
the colligation operator \(U\) is a unitary operator, i.e. if
\begin{equation}
\label{UOp} U^{\ast}U=I_{\mathcal{E}^{in}\oplus{}\mathcal{H}},\qquad \
UU^{\ast}=I_{\mathcal{E}^{out}\oplus{}\mathcal{H}}\,.
\end{equation}
\end{definition}

\begin{definition}
\label{DeInMF}
Let \(\mathcal{E}_1\) and \(\mathcal{E}_2\) be finite-dimensional
Hilbert spaces and let \(S(z)\) be a rational function whose values are
operators acting from \(\mathcal{E}_1\) to \(\mathcal{E}_2\).

The matrix-function \(S\) is called an \textsf{inner} function if its values
\(S(z)\) are contractive operators for \(z\in\mathbb{D}\) and
unitary operators for \(t\in\mathbb{T}\), i.e. if the conditions
\begin{subequations}
\label{IForm}
\begin{alignat}{3}%
\label{IForm1}
& & I_{\mathcal{E}_1}-S^\ast(z)S(z)\geq{}0\,,\ \ %
&I_{\mathcal{E}_2}-S(z)S^\ast(z)\geq{}0,
\ \ &\text{for} \ \ z\in\mathbb{D}, \\[1.0ex]
& & I_{\mathcal{E}_1}-S^\ast(t)S(t)\,=0\,,\ \ &%
I_{\mathcal{E}_2}-S(t)S^\ast(t)\,=0 \ , \ \ &\text{for} \ \ %
t\in\mathbb{T}\,.
\label{IForm2}
 \end{alignat}
 \end{subequations}
 hold. (In particular, \(S\) has no singularities in \(\mathbb{D}\cup\mathbb{T}\).)
\end{definition}

\begin{remark}
Since unitary operators are invertible,
 \(\mathcal{E}_1\,\text{-}\,\mathcal{E}_2\) inner functions
exist only if \(\dim\mathcal{E}_1=\dim\mathcal{E}_1\).
\end{remark}

\begin{theorem}
\label{PCFUC}%
Let \((\mathcal{E}^{in},\,\mathcal{E}^{out},\,\mathcal{H},\,U)\),
\eqref{CODB}\,-\,\eqref{FWTW}, be a
finite-dimensional unitary colligation and \(S(z)\),
\eqref{CharFunc}, be its characteristic function.

 Then the function \(S(z)\) is a rational inner function.
\end{theorem}
\noindent
 \textsf{PROOF.} The proof of this lemma is based on
identity \eqref{LSDSFr}, where \(h(z)\) is expressed in terms of
\(\varphi(z)\) as in \eqref{LinS2}. Let \(z\) and \(\zeta\) be such that
the operators \(I-zD\) and \(I-\zeta{}D\) are invertible (These operators
are invertible if \(z\in\mathbb{D},\zeta\in\mathbb{D}\). Also, since the
spectrum of the operator \(D\) is a finite set, the operators \(I-zD\) and \(I-\zeta{}D\)
are invertible for all but finitely many \(z\in\mathbb{T},\zeta\in\mathbb{T}\).)
 Because the operator \(U\) is unitary, \eqref{LSDSFr} yields
\[\langle\psi(z),\psi(\zeta)\rangle_{\mathcal{E}^{out}}+
 (z\overline{\zeta})^{-1}\langle\,h(z),h(\zeta)\,\rangle_{\mathcal{H}}=
\langle\varphi(z),\varphi(\zeta)\rangle_{\mathcal{E}}+
\langle\,h(z),h(\zeta)\,\rangle_{\mathcal{H}}\,,\] or
\begin{multline}
\label{id1}%
 \big\langle\varphi(z),\varphi(\zeta)\big\rangle_{\mathcal{E}^{in}}-
\big\langle{}S(z)\varphi(z),S(\zeta)\varphi(\zeta)\big\rangle_{\mathcal{E}^{out}}=\\[1.0ex]
=(1-z\overline{\zeta})\big\langle(I-zA)^{-1}C\varphi(z),
(I-\zeta{}A)^{-1}C\varphi(\zeta) \big\rangle_{\mathcal{H}}\,.
\end{multline}
In particular, taking
\(\varphi(z)\equiv\varphi^\prime\) and \(\varphi(\zeta)\equiv\varphi^{\prime\prime}\),
where \(\varphi^\prime\) and \(\varphi^{\prime\prime}\) are arbitrary
vectors in \(\mathcal{E}^{in}\), we obtain the equality
\begin{equation}%
\label{Id1}
\frac{I_{\mathcal{E}^{in}}-S^{\ast}(\zeta)S(z)}{1-\overline{\zeta}z}=
C^{\ast}(I-\overline{\zeta}{}D^{\ast})^{-1}(I-zD)^{-1}C\,.
\end{equation}
In the same way we obtain the equality
\begin{equation}%
\label{Id2}
\frac{I_{\mathcal{E}^{out}}-S(z)S^{\ast}(\zeta)}{1-z\overline{\zeta}}=
B(I-zD)^{-1}(I-\overline{\zeta}{}D^{\ast})^{-1}B^{\ast}\,.
\end{equation}
Using the identity
\(\textstyle{\frac{\zeta(I-\zeta{}D)^{-1}-z(I-zD)^{-1}}{\zeta-z}=(I-\zeta{}D)^{-1}(I-zD)^{-1}},\)
we obtain
\begin{equation}%
\label{Id3}%
 \frac{S(\zeta)-S(z)}{\zeta-z}=
B(I-\zeta{}D)^{-1}(I-zD)^{-1}C\,,
\end{equation}
and
\begin{equation}%
\label{Id4}%
 \frac{S^{\ast}(\zeta)-S^{\ast}(z)}{\overline{\zeta}-\overline{z}}=
C^{\ast}(I-\overline{\zeta}{}D^{\ast})^{-1}(I-\overline{z}D^{\ast})^{-1}B^{\ast}\,,
\end{equation}

 To get \eqref{IForm} we let \(\zeta=z\) in
\eqref{Id1}\,-\,\eqref{Id2}:
\begin{subequations}
\label{SpRe}%
\begin{gather}
\label{SpRe1}%
I_{\mathcal{E}^{in}}-S^{\ast}(z)S(z)=(1-|z|^2)\,%
C^{\ast}(I-\overline{z}{}A^{\ast})^{-1}(I-zA)^{-1}C,\,\\[1.0ex]
I_{\mathcal{E}^{out}}-S(z)S^{\ast}(z)=(1-|z|^2)\,%
B(I-zA)^{-1}(I-\overline{z}{}A^{\ast})^{-1}B^{\ast}\,.
\label{SpRe2}%
\end{gather}
\end{subequations}
The inequalities \eqref{IForm1} follow from equalities \eqref{SpRe}, which hold
for all \(z\in\mathbb{D}\). The equalities \eqref{SpRe} furthermore hold for all
but finitely many \(z\in\mathbb{T}\). Thus, the rational function \(S(z)\)
is bounded in \(\mathbb{T}\), except on a finite set. \(S\) therefore has
no singularities in \(\mathbb{T}\) and takes unitary values there.\\

The following theorem serves as a `unitary' counterpart to Theorem \ref{RealTh}.
\begin{theorem}%
\label{ReThU}%
Let \(S(z)\) be a  \emph{rational inner} function
whose values are
operators acting from \(\mathcal{E}_1\) into \(\mathcal{E}_2\),
where \(\mathcal{E}_1\) and \(\mathcal{E}_2\) are finite-dimensional
Hilbert spaces.

Then there exists a finite-dimensional, \emph{minimal},
unitary operator colligation
\((\mathcal{E}_1^{in},\,\mathcal{E}_1^{out},\,\mathcal{H},\,U)\),
\eqref{UnCol}\,-\,\eqref{CODB}\,-\,\eqref{FWTW}, with
 \(\mathcal{E}^{in}=\mathcal{E}_1\) and \( \mathcal{E}^{out}=\mathcal{E}_2\),
whose characteristic function \(S_U(z)=A+zB(I-zD)^{-1}C\) coincides with
 the original function \(S(z)\). In other words, the function \(S\) is representable
 in the form \eqref{CharFunc}.
\end{theorem}%

\begin{definition}
\label{DeUnEq}%
 Let \((\mathcal{E}_1^{in},\,\mathcal{E}_1^{out},\,\mathcal{H}_1,\,U_1)\) and
\((\mathcal{E}_2^{in},\,\mathcal{E}_2^{out},\,\mathcal{H}_2,\,U_2)\) be operator
colligations, \eqref{DCO}\,-\,\eqref{fwtw}. If these colligations are equivalent
(i.e. if they satisfy the intertwining relation \eqref{LinEq}\,-\eqref{FWTWb})
and each of the operators \(E^{in},\,E^{out},\,V\) is unitary, we say that
\((\mathcal{E}_1^{in},\,\mathcal{E}_1^{out},\,\mathcal{H}_1,\,U_1)\) and
\((\mathcal{E}_2^{in},\,\mathcal{E}_2^{out},\,\mathcal{H}_2,\,U_2)\) are
\textsf{unitarily equivalent}.
\end{definition}

Clearly, if two operator colligations are unitarily equivalent
and one of these colligations is unitary, then the second colligation
is also unitary.

The following theorem provides us with a `unitary' counterpart to Theorem \ref{ChSO}.
\begin{theorem}
\label{CChF}%
 Let \((\mathcal{E}_1^{in},\,\mathcal{E}_1^{out},\,\mathcal{H}_1,\,U_1)\) and
\((\mathcal{E}_2^{in},\,\mathcal{E}_2^{out},\, \mathcal{H}_2,\,U_2)\) be
unitary colligations, \eqref{DCO}. Furthermore, let these colligations be unitarily
equivalent, i.e. suppose that the intertwining relation \eqref{LinEq} holds
for some \emph{unitary} operators \(E^{in},\,E^{out}\) and \(V\).

The respective characteristic functions \(S_1(z)\) and \(S_2(z)\) of these
colligations, \eqref{CFs}, then satisfy the intertwining relation \eqref{Intw}
with these very same \emph{unitary} operators \(E^{in}\) and \(E^{out}\).
\end{theorem}

If we, furthermore, assume that both unitary colligations are simple, we
can show that the converse to Theorem \ref{CChF} also holds.

The next theorem serves as a `unitary' counterpart to Theorem \ref{ChSOIn}.

\begin{theorem}
\label{ChSOInU}%
Let \((\mathcal{E}_1^{in},\,\mathcal{E}_1^{out},\,\mathcal{H}_1,\,U_1)\) and
\((\mathcal{E}_2^{in},\,\mathcal{E}_2^{out},\, \mathcal{H}_2,\,U_2)\) be
 finite-dimensional \emph{unitary} operator colligations, \eqref{DCO}.
 Let \(S_1(z)\) and \(S_2(z)\), \eqref{CFs}, be the characteristic functions of
\((\mathcal{E}_1^{in},\,\mathcal{E}_1^{out},\,\mathcal{H}_1,\,U_1)\) and
\((\mathcal{E}_2^{in},\,\mathcal{E}_2^{out},\, \mathcal{H}_2,\,U_2)\), respectively.
We now make the following assumptions:
\begin{enumerate}
\item[\textup{1.}]
The functions
\(S_1(z)\) and \(S_2(z)\) satisfy
the intertwining relation \eqref{Intw} for \(z\in\mathbb{D}\),
where \(E^{in}:\mathcal{E}_2^{in}\to{}\mathcal{E}_1^{in},\,
E^{out}:\mathcal{E}_2^{out}\to{}\mathcal{E}_1^{out}\) are some unitary operators.
\item[\textup{2.}]
The colligations \((\mathcal{E}_1^{in},\,\mathcal{E}_1^{out},\,\mathcal{H}_1,\,U_1)\) and
\((\mathcal{E}_2^{in},\,\mathcal{E}_2^{out},\, \mathcal{H}_2,\,U_2)\) are simple.
\end{enumerate}

The colligations \((\mathcal{E}_1^{in},\,\mathcal{E}_1^{out},\,\mathcal{H}_1,\,U_1)\) and
\((\mathcal{E}_2^{in},\,\mathcal{E}_2^{out},\, \mathcal{H}_2,\,U_2)\)
are then unitarily equivalent, i.e. there exists a unitarily
operator \(V:\mathcal{H}_2\to\mathcal{H}_1\) such that the intertwining
relation \eqref{LinEq} holds.
\end{theorem}

Let us compare the assumptions of \textup{Theorems \ref{ChSOIn}}
and \textup{\ref{ChSOInU}}. In \textup{Theorem~\ref{ChSOIn}}
we assume that the colligations
\((\mathcal{E}_i^{in},\,\mathcal{E}_i^{out},\,\mathcal{H}_i,\,U_i),\,i=1,2,\)
are minimal, however it is not assumed that these colligations are
unitary. In \textup{Theorem \ref{ChSOInU}} we assume that the
colligations are unitary and simple, but we do not explicitly assume
that these colligations are minimal, because they \emph{are}, in fact, already
\emph{minimal}.

\begin{theorem}
\label{SimMin}%
 Let \((\mathcal{E}^{in},\,\mathcal{E}^{out},\,\mathcal{H},\,U)\)  be
a finite-dimensional, \emph{unitary} operator colligation.
The following statements are then equivalent:
\begin{enumerate}
\item[\textup{1.}] The colligation is simple.
\item[\textup{2.}] The colligation is minimal.
\item[\textup{3.}] The colligation is controllable.
\item[\textup{4.}] The colligation is observable.
\end{enumerate}
\end{theorem}

\vspace{1.0ex} \emph{In what follows we deal only with scalar-valued
inner functions} \(S(z)\), i.e. with functions whose
values are complex numbers. The input space
\(\mathcal{E}^{in}\) and the output space \(\mathcal{E}^{out}\) of
the unitary colligation
\((\mathcal{E}^{in},\mathcal{E}^{out},\,\mathcal{H},\,U)\)
representing this \(S(z)\) can be identified with the space
\(\mathbb{C}\): \(\mathcal{E}^{in}=\mathcal{E}^{out}=\mathbb{C}\).
The finite-dimensional state space \(\mathcal{H}\), with, say
\( \ \dim \mathcal{H} = n\), can be identified with the space
\(\mathbb{C}^n\) (with the standard scalar product):
\(\mathcal{H}=\mathbb{C}^n\). With these conventions in place, the orthogonal sums
\(\mathcal{E}^{in}\oplus\mathcal{H}\) and \(\mathcal{E}^{out}\oplus\mathcal{H}\)
can be identified naturally with the space \(\mathbb{C}\oplus\mathbb{C}^n\).

We note that \(\mathbb{C}\oplus\mathbb{C}^n\) represents a canonical
decomposition of the space \(\mathbb{C}^{n+1}\) into an orthogonal sum.
We consider the space \(\mathbb{C}^{n+1}\) as the set
\(\mathfrak{M}_{(n+1)\times{}1}\) of all \((n+1)\)-column-vectors,
along with the standard linear operations and scalar product:
\begin{equation}
\label{StScPr} \langle{}f\,,\,g{}\rangle=g^\ast{}f,\quad f,\,g\in
\mathfrak{M}_{(n+1)\times{}1},
\end{equation}
where the asterisk \({}^\ast\) denotes Hermitian conjugation.

 A unitary operator, \(U\), acting in
\(\mathbb{C}^{n+1}\) is described by a unitary \((1+n)\times(1+n)\)-matrix,
which will also be denoted by \(U\). \(U\) maps the column-vector \(f\) to
the column-vector \(Uf\), where \(Uf\) is the usual matrix product.
The decomposition \(\mathbb{C}^{n+1}=\mathbb{C}\oplus\mathbb{C}^n\) of
the space \(\mathbb{C}^{n+1}\) suggest that we consider the following
block-matrix decomposition of \(U\):
\begin{subequations}
\label{CBMD}
\begin{gather}
\label{CBMDa}
 U=
\begin{bmatrix}
A&B\\[0.7ex]
C&D
\end{bmatrix}\,,\\[0.5ex]
\label{CBMDb} A\in\mathfrak{M}_{1\times{}1},\,
B\in\mathfrak{M}_{1\times{}n},\, C\in\mathfrak{M}_{n\times{}1},\,
D\in\mathfrak{M}_{n\times{}n}\,.
\end{gather}
\end{subequations}
The matrix entries are considered as operators: %
\begin{equation}
\label{Act} %
A: \mathcal{E}\to\mathcal{E}\,, \
B:\mathcal{H}\to\mathcal{E}\,, \
C: \mathcal{E}\to\mathcal{H}\,, \  %
D:  \mathcal{H}\to\mathcal{H}\,.
\end{equation}
where
\begin{equation}
\label{Spa} %
\mathcal{E}=\mathfrak{M}_{1\times{}1} \ \ (=\mathbb{C}), \ \
\mathcal{H}=\mathfrak{M}_{n\times{}1} \ \ (=\mathbb{C}^n)\,.
\end{equation}

\begin{definition}
\label{UCCUM}%
 Given a unitary matrix \(U\in\mathfrak{M}_{(n+1)\times{}(n+1)}\) with block decomposition
 \eqref{CBMD}, we associate the unitary
colligation \((\mathcal{E},\,\mathcal{H},\,U)\) with \(U\).
The exterior space \(\mathcal{E}\) and the state space
\(\mathcal{H}\) of this colligation are as in \eqref{Spa}, where
the spaces \(\mathbb{C}\) and \(\mathbb{C}^n\) have the
standard scalar products. The exterior, principal and channel
operators \(A,\,D,\,B,\,C\) correspond to the block-matrix entries
in \eqref{CBMD} and satsify \eqref{Act}.

We call this colligation the \textsf{unitary
colligation associated with the unitary matrix} \(U\).
\end{definition}

Given two unitary colligations associated with unitary matrices
\(U^{\prime}\) and \(U^{\prime\prime}\), how do we express that
these colligations are unitarily equivalent? The exterior spaces of
both colligations are `copies' of the same space \(\mathbb{C}\).
To identify the exterior spaces
\(\mathbb{C}\) of two different colligations, we should  specify
 the unitary operators \(E^{in}\) and \(E^{out}\) for
 the two copies of \(\mathbb{C}\)
 (These operators, \(E^{in}\) and \(E^{out}\), appear in \eqref{FWTWb}
and in the intertwining relations \eqref{LinEq} and
\eqref{Intw}.) We can naturally choose these
identification operators as the identity operators, i.e. such that each of
operators \(E^{in}\) and \(E^{out}\) is represented by the
\(1\times{}1\)-matrix whose (unique) entry is the number \(1\)
(Such operators can be represented by \(1\times{}1\)-matrices, where the
matrices corresponding to \(E^{in}\) and \(E^{out}\) consist, respectively, of an
arbitrary number \(\nu^{in}\) and \(\nu^{out}\) with
\(|\nu^{in}|=1\) and \(|\nu^{out}|=1\).)

With this convention in place, the unitary equivalence of the colligations
associated with the block-matrices
 \begin{equation}%
 \label{TUM}
 U^{\prime}=\begin{bmatrix}
A^{\prime}&B^{\prime}\\[1.0ex]
C^{\prime}&D^{\prime}
\end{bmatrix}%
\in\mathfrak{M}_{(n+1)\times{}(n+1)}%
 \quad\text{and} \quad
U^{\prime\prime}=\begin{bmatrix}
A^{\prime\prime}&B^{\prime\prime}\\[1.0ex]
C^{\prime\prime}&D^{\prime\prime}
\end{bmatrix}%
\in\mathfrak{M}_{(n+1)\times{}(n+1)}%
\end{equation}
 means that these matrices
 satisfy the intertwining relation:
\begin{equation}
\label{LinEqM}
\begin{bmatrix}
1&0\\[1.0ex]
0&V
\end{bmatrix}
\begin{bmatrix}
A^{\prime\prime}&B^{\prime\prime}\\[1.0ex]
C^{\prime\prime}&D^{\prime\prime}
\end{bmatrix}=
\begin{bmatrix}
A^{\prime}&B^{\prime}\\[1.0ex]
C^{\prime}&D^{\prime}
\end{bmatrix}
\begin{bmatrix}
1&0\\[1.0ex]
0&V
\end{bmatrix}
\end{equation}
where \(V\in\mathfrak{M}_{\,n\times{}n}\) is a unitary matrix.
The equality \eqref{Intw} then becomes:
\begin{equation*}
    S_1(z)=S_2(z).
\end{equation*}

\begin{definition}
\label{DEUM}%
We say that the unitary matrices
\(U^{\prime}\in\mathfrak{M}_{(n+1)\times{}(n+1)}\) and
\(U^{\prime\prime}\in\mathfrak{M}_{(n+1)\times{}(n+1)}\),
\eqref{TUM}, are \textsf{equivalent} if there exists a
 unitary matrix \(V\in\mathfrak{M}_{\,n\times{}n}\) such that the
 intertwining relation \eqref{LinEqM} holds.
\end{definition}

Let
 \(U=\begin{bmatrix}
A&B\\[1.0ex]
C&D
\end{bmatrix}\in\mathfrak{M}_{(n+1)\times{}(n+1)}\)\,.
We now consider the following matrices associated with the unitary matrix \(U\):
\begin{subequations}
\label{COM}
\begin{alignat}{2}
\label{COMa}%
 &\mathscr{C}(U)=
[C\,\,\,\,\,DC\,\,\,\,\,\,\,\,\,\ldots\,\,\,\,\,\,\,\,D^{n-1}C],&\quad
&\mathscr{C}(U)\in\mathfrak{M}_{\,n\times{}n}\,,\\[1.5ex]
\label{COMb}
&\mathscr{B}(U)=[B^{\ast}\,\,\,D^{\ast}B^{\ast}\,\,\,\ldots\,\,\,(D^{\ast})^{n-1}B^{\ast}],&
&\mathscr{B}(U)\in\mathfrak{M}_{\,n\times{}n}\,,
\end{alignat}
and
\begin{multline}
\label{COMc}%
\mathscr{S}(U)=
[C\,\,\,DC\,\,\,\,\,\ldots\,\,\,\,\,D^{n-1}C,\,B^{\ast}\,\,\,
D^{\ast}B^{\ast}\,\,\,\ldots\,\,\,(D^{\ast})^{n-1}B^{\ast}]\,,\\[0.5ex]
\mathscr{S}(U)\in\mathfrak{M}_{\,n\times{}2n}\,.
\end{multline}
\end{subequations}
If the unitary colligation associated with the matrix
\(U\) is controllable, observable or simple, this means that the
matrix \eqref{COMa},  \eqref{COMb} or \eqref{COMc} is, respectively,
of rank \(n\).
\begin{remark}%
\label{RemRan}%
If one of the matrices \eqref{COM} has rank \(n\), then
    its columns (considered as vectors in
 \(\mathbb{C}^n=\mathfrak{M}_{n\times{}1}\)) generate the whole
 space. The columns of these matrices are of the form \(D^{k}C\)
 or \((D^{\ast})^kB^{\ast}\), where \(k\) takes values in the
 interval \([0,\,\ldots\,,\,(n-1)]\). It is possible
 to consider matrices of this kind for \(k\) over a
 larger interval. Extending the interval \([0,\,\ldots\,,\,(n-1)]\)
 does not, however, lead to an increase in rank for these matrices:
 The Cayley-Hamilton Theorem tells us that the column-vectors, \(D^{p}C\)
 and \((D^{\ast})^kB^{\ast}\) with \(p\geq{}n\), are, respectively,
 linear combinations of the column-vectors
 \(D^{k}C\) and \((D^{\ast})^kB^{\ast}\)  with \(k\in[0,\,\ldots\,,\,(n-1)]\).
\end{remark}%
\begin{definition}
\label{DMM}%
    We say that a unitary matrix \(U\in\mathfrak{M}_{(n+1)\times{}(n+1)}\), \linebreak
    expressed using the block-decomposition in \eqref{CBMD},
  is \textsf{controllable} if \linebreak \( \ \mathrm{rank} ~ \mathscr{C}(U) = n\),
  \textsf{observable} if \( \ \mathrm{rank} ~ \mathscr{B}(U) = n\)
  and \textsf{simple} if \( \ \mathrm{rank} ~ \mathscr{S}(U) = n\).
    If the matrix \(U\) is both controllable and observable, we say that
    it is \textsf{minimal}.
\end{definition}
(We note that any one of the matrices \eqref{COM} is of rank \(n\) if and
only if the other two have rank \(n\). See Theorem \ref{SimMin}.)

The results of this section on the state space representation of scalar
(i.e. complex-valued) rational inner functions can be
summarized in the following way:
\begin{theorem}  \label{SyR}
\textup{(\textsf{{\small\,Rational Inner
Functions \(\boldsymbol\Longleftrightarrow\)\,Equivalence Classes of Unitary Matrices}})}
\newline
\begin{enumerate}
\item[\textup{1.}]
 Let \(S(z)\) be an inner rational function
 of degree \(n\).
Then \(S(z)\) can be represented in the form:
 \begin{equation}
 \label{SSR}%
 S(z)=A+zB(I_n-zD)^{-1}C,\,
 \end{equation}
 where \(A,B,C,D\) are blocks of some  unitary minimal matrix
 \(U\),\\
  \(U\in\mathfrak{M}_{(n+1)\times{}(n+1)}\), \eqref{CBMD}\,.
\item[\textup{2.}] Let \(U\in\mathfrak{M}_{(n+1)\times{}(n+1)}\)
be a unitary  matrix with block-decomposition \eqref{CBMD}
and let the function \(S(z)\) be defined in terms of \(U\) by
\eqref{SSR}.%
Then the function \(S(z)\) is a rational inner function with \(\deg{}%
S\leq n\). If the matrix \(U\) is minimal, then \(\deg{} S=n\).
\item[\textup{3.}] Let \(U^{\prime}\in\mathfrak{M}_{(n+1)\times{}(n+1)}\)
 and \(U^{\prime\prime}\in\mathfrak{M}_{(n+1)\times{}(n+1)}\) be
 unitary matrices with block-decomposition
\eqref{TUM} and let \(S^{\prime}(z)\) and \(S^{\prime\prime}(z)\)
be the functions defined in terms of \(U^{\prime}\) and
\(U^{\prime\prime}\) by:
\begin{equation}
 \label{SSRE}%
 S^{\prime}(z)=A^{\prime}+zB^{\prime}(I_n-zD^{\prime})^{-1}C^{\prime},\,\qquad
 S^{\prime\prime}(z)=A^{\prime\prime}+
 zB^{\prime}(I_n-zD^{\prime\prime})^{-1}C^{\prime\prime},\,
 \end{equation}
 If the matrices \(U^{\prime}\) and
 \(U^{\prime\prime}\) are equivalent, then \(S^{\prime}(z)\equiv{}
 S^{\prime\prime}(z)\). \linebreak
 If \(S^{\prime}(z)\equiv{}
 S^{\prime\prime}(z)\) and the matrices \(U^{\prime}\) and
 \(U^{\prime\prime}\) are minimal, then \(U^{\prime}\) and
 \(U^{\prime\prime}\) are equivalent.
 \end{enumerate}
\end{theorem}

The substance of this theorem can be summarized as follows:
\begin{itemize}
\item[\(\circ\)]
There exists a one-to-one correspondence between the set of all
rational inner functions of degree \(\leq{}n\) and the set of all
equivalence classes of unitary matrices in
\(\mathfrak{M}_{(n+1)\times{}(n+1)}\).
\item[\(\circ\)]
This correspondence can be expressed as a mapping from the set of
all rational inner functions of degree \(n\) onto the set of all equivalence
classes of minimal unitary matrices in \(\mathfrak{M}_{(n+1)\times{}(n+1)}\).
\end{itemize}

For a proof of Theorem \ref{SyR}, see the Appendix at the end of this
paper.\\

\centerline{\textsf{The Main Objective of This Paper.}}
 Applying the Schur algorithm to a given rational inner function \(s(z)\)
 of degree \(n\) produces the sequence \(s_k(z),\,k=0,\,1,\,\dots n,\) of rational inner
functions with \(s_0(z)=s(z)\) and \(\deg{}s_k(z)=n-i\). In particular,
\(s_n(z)\equiv{}s_n\) is a unitary constant. According to what
was stated in Section \ref{SRIF}, each of the functions \(s_k(z)\)
admits a system representation,
\begin{equation}
\label{SRUSF} s_k(z)=A_k+zB_k(I-zD_k)^{-1}C_k\,,
\end{equation}
in terms of the blocks of some minimal unitary matrix
\(U_k\in\mathfrak{M}_{(1+n-k)\times{}(1+n-k)}\):
\begin{equation}
U_k=
\begin{bmatrix}
A_k&B_k\\
C_k&D_k
\end{bmatrix}\,.
\end{equation}
We assume that from the very beginning, the given inner rational
function \(s(z)=s_0(z)\) is \textsf{determined} in terms of its
state space representation, so that the matrix \(U_0\) is given. The
goal is to recursively produce the sequence of matrices \(U_k\)
representing the functions \(s_k(z)\), \(k=1,\,2,\,\dots\,,n\).
The matrix \(U_{k+1}\), representing the
function \(s_{k+1}(z)\), is thus constructed from the matrix \(U_{k}\),
representing the function \(s_{k}(z)\). In other words, the
steps \eqref{kas} of the Schur algorithm must be described in terms of the
state space representation \eqref{SRUSF}.

It should be noted that the unitary matrices in the system representation
of a rational inner function are determined only up to the
equivalence
\begin{equation}
\label{UnEq}
\begin{bmatrix}
A_k&B_k\\
C_k&D_k
\end{bmatrix}\sim
\begin{bmatrix}
1&0\\
0&V_k^{-1}
\end{bmatrix}\,\cdot\,
\begin{bmatrix}
A_k&B_k\\
C_k&D_k
\end{bmatrix}\,\cdot\,
\begin{bmatrix}
1&0\\
0&V_k
\end{bmatrix}\,,
\end{equation}
where \(V_k\) is an arbitrary unitary \(k\times{}k\) matrix.
 So we have to find a rule for
constructing a matrix \(U_{k+1}\), which belongs to the
equivalence class of matrices representing the function
\(s_{k+1}(z)\), from an arbitrary element \(U_{k}\) of the
equivalence class of matrices representing the function
\(s_{k}(z)\).

\begin{center}
\textrm{The Schur algorithm in the framework of system representations\\
is described in Section \ref{SchA}.}
\end{center}


\textsf{\large{Historical Remarks.}} %
The definition of a characteristic function was developed gradually,
starting from the pioneering works of M.S.Livshitz. The first definition
appeared in \cite{Liv1} (for operators for which \(I-T^{\ast}T\)
and \(I-TT^{\ast}\) have rank one) and in \cite{Liv2} (for the
case that these operators have finite rank). M.S.Livshitz and
those working in the same field, subsequently turned their attention to
bounded operators \(T\) for which \(T-T^{\ast}\) is of finite rank
or at least of finite trace. For these operators \(T\), a
characteristic function was defined in an analogous way and by
means of this function, a wide-reaching theory for these operators developed.
In particular, triangular models of non-self-adjoint operators
were introduced. See \cite{Liv3}, \cite{BrLi}, \cite{Br}. In the
course of the evolution of the concept of characteristic functions,
it became clear that it was advantageous to consider, not just
non-self-adjoint operators, but also more general
objects: operator nodes (or operator colligations). The
notion of an operator colligation was prompted by physical
applications of the Livshitz theory of non-selfadjoint operators.
(See \cite{BrLi}, \cite{Liv9} and references there.)

B.\,Sz.\,Nagy and C.\,Foias used a different approach to characteristic functions in
1962. Their work involved harmonic analysis of the unitary dilation of the contractive operator
\(T\). Moreover, they simultaneously obtained a functional model
of \(T\) depending explicitly and exclusively on the characteristic
function of \(T\). See \cite[especially Chapter VI]{SzNFo} and
references therein.

The version of operator colligations, which appears in
Definition
 \ref{DOC} goes back to a remark of M.G.\,Krein to
the work \cite{BrSv1}.
 In \cite{BrSv1}, the notion of a contractive
operator colligation (node) was defined as the collection of
Hilbert spaces \(\mathcal{H},\mathcal{F},\ \mathcal{G}\) and
operators
\begin{equation}%
\label{HS}
 T_0:\,\mathcal{G}\to{}\mathcal{F},\
F:\,\mathcal{F}\to{}\mathcal{H}, \
G:\,\mathcal{G}\to{}\mathcal{H}, \
T:\,\mathcal{H}\to{}\mathcal{H}\,,
\end{equation}
satisfying the conditions
\begin{multline}%
\label{OI}
 I-TT^{\ast}=FF^{\ast},\,I-T^{\ast}T=GG^{\ast},\,\\
 I-T_0T_0^{\ast}=F^{\ast}F,\ I-T_0^{\ast}T_0=G^{\ast}G,\,
 TG=FT_0\,,
\end{multline}
The results presented in the paper \cite{BrSv1} were reported on
in a seminar of Krein's in Odessa. In the remark to this talk,
M.G.Krein noticed that the conditions \eqref{HS}-\eqref{OI} mean
that the block-operator
\begin{equation}
\begin{bmatrix}
T_0^{\ast}& G^{\ast}\\
-F& T
\end{bmatrix}:\ \ \
\begin{bmatrix}
\mathcal{F}\\ \mathcal{H}
\end{bmatrix}
\to
\begin{bmatrix}
\mathcal{G}\\ \mathcal{H}
\end{bmatrix}\,,
\end{equation}
acting in the appropriate orthogonal sums of Hilbert spaces, is a
unitary operator.   Starting from this remark of M.G.\,Krein's,
mathematicians belonging to the Odessa school as well as other
mathematicians, defined the operator colligation as the block
operator acting from the direct sum
\begin{math}
\begin{bmatrix}
\textit{input space}\\ \textit{state space}
\end{bmatrix}
\end{math}
into the direct sum
\begin{math}
\begin{bmatrix}
\textit{output space}\\ \textit{state space}
\end{bmatrix}
\end{math}\,.
If the spaces have scalar products and the block
operator is a unitary operator with respect to this product, then
the operator colligation is called an \textit{unitary colligation}.

It should be mentioned that the paper \cite{BrSv1} has connections
to the theory of functional models of contractive operators
developed in \cite{SzNFo}. The definition \eqref{CharFunc} of the
characteristic function of the colligation
\eqref{CODB}\,-\,\eqref{FWTW} agrees with the definition of the
characteristic function in \cite{BrSv1}.

The notions of controllability and observability (and minimality)
in the setting of State Space Theory were introduced by R.Kalman
in \cite{Kal1}. The study of controllability and observability of
composite systems was first dealt with in \cite{Gil}. Under other
names, the notion of controllability also appears in the Livshitz
theory of open systems. See the notions of the simple system and
of the complementary component in section 1.3 of \cite{Liv9}. (See
pages  36\,-\,37 of the Russian original, or pages 27-29 of the
English translation.)

The fact that every rational matrix-function \(S\) can be realized
as the transfer function of some \textit{minimal} stationary
linear system (which here appears as Theorem \ref{RealTh}),
the uniqueness of the state space representation (Theorem \eqref{ChSOIn})
and the equality \(\dim \mathcal{H}=\deg S\) were all
established by R.\,Kalman in a very general setting.
These results, as well as many other results, can be found in Chapter 10
of \cite{KFA}. See also Chapter 1 of \cite{Fuh}.

Some algorithms for the system realization of a given rational
function were proposed by R.\,Kalman and his collaborators. (See
Chapter 10 of the monograph \cite{KFA} and references there.)
R.\,Kalman did not consider questions related to the realization
of contractive or inner matrix-functions: He developed system
theory over arbitrary fields rather than over the field of complex
numbers.

An excellent (and short!) presentation of the state space approach
to the problems of minimal realization and factorization of
rational functions can be found in \cite{Kaa}.

Realizations of contractive or inner rational
matrix-functions (rational and more general) were later considered
in the framework of the SzNagy-Foias model for contractive operators.
These and also more general results can be found in many
publications now. For convenience, we
present some basic facts on system realizations for inner rational
functions (scalar) in the Appendix to the present paper.

The state space description of the composite system, which is
formed by the cascade (or Redheffer) coupling of several state space
systems, was dealt with in \cite{HeBa} in more generality. We
make use of these results, but prefer to derive them
independently of \cite{HeBa} in the form and in the generality
which is most suitable for our goal.


\section{Coupled Systems and The Schur Transformation
: Input-Output Mappings.\label{SASRio}}
\setcounter{equation}{0}

To describe the Schur algorithm using system
representations, we must first consider how the
linear-fractional Schur transformation
\begin{equation}
\label{FLT}
\omega(z)\rightarrow{}s(z),\quad s(z)=\frac{s_0+z\omega(z)}{1+\,z\omega(z)\overline{s_0}}\
\qquad \text{(\(s_0\) is a complex number, \(|s_0|<1\))}
\end{equation}
 can be described in terms of the \textit{input-output mappings of linear systems.}
The linear-fractional transform \eqref{FLT} is of the form
\begin{equation}
\label{tFLT}%
 s(z)=\frac{w_{11}(z)\omega(z)+w_{12}(z)}{w_{21}(z)\omega(z)+w_{22}(z)}\,.
\end{equation}
This form of a linear-fractional transform is the most familiar to the classical analyst.
In the theory of unitary operator colligations, the
\textsf{Redheffer%
\footnote%
{Raymond Redheffer (1921-2005) was a US mathematician working at UCLA.} %
form} for linear-fractional transforms, i.e.
\begin{equation}
\label{Redh}
 s(z)=s_{11}(z)+s_{12}(z)\omega(z)(I-s_{22}(z)\omega(z))^{-1}s_{21}(z)\,.
\end{equation}
is often more convenient. Every linear-fractional transformation of the form
\eqref{tFLT} can be rewritten in the Redheffer form \eqref{Redh},
but not every transformation in Redheffer form can be expressed in
linear-fractional form.

The matrix %
\begin{math}
         \ W(z)=\Bigl[\begin{smallmatrix}
        w_{11}(z)&w_{12}(z)\\[0.5ex]
        w_{21}(z)&w_{21}(z)
        \end{smallmatrix}\Bigr] \
\end{math}%
 for the transformation \eqref{FLT} and \eqref{tFLT}
(under the appropriate normalization%
\footnote%
{The matrix of the linear-fractional
transform \eqref{tFLT} is determined only up to the proportionality
\(W(z)\to\lambda(z)W(z)\), where \(\lambda\in\mathbb{C}\setminus{} \{0\}.\)}%
) is
\begin{equation}
\label{FLT1}%
 W(z)=(1-|s_0|^2)^{-1/2}
\begin{bmatrix}
z&s_0\\[0.5ex]
z\overline{s_0}&1
\end{bmatrix}\,.
\end{equation}
\(W(z)\) in \eqref{FLT1} is not an inner matrix but it is
 a \(j\)-inner matrix:
\[j-W^{\ast}(z)jW(z)\geq{}0, \ \ z\in\mathbb{D},\qquad j-W^{\ast}(t)jW(t)=0, \ \
 t\in\mathbb{T}\,,\]
 where
 \begin{math}j=\begin{bmatrix}
1&0\\[0.5ex]
0&-1
\end{bmatrix}\,.
\end{math}


Let us express the fractional-linear transformation \eqref{FLT}
in the Redheffer form \eqref{Redh},
where the \(2\times{}2\)-matrix-function %
\(S(z)=%
\Bigl[\begin{smallmatrix}
s_{11}(z)&s_{12}(z)\\
s_{21}(z)&s_{22}(z)
\end{smallmatrix}\Bigr]\) is:
\begin{equation}
\label{TSM}%
 S(z)=%
\begin{bmatrix}
s_0&\ \ z\,(1-|s_0|^2)^{1/2}\\[2.0ex]
(1-|s_0|^2)^{1/2}&-z\,\overline{s_0}
\end{bmatrix}
\end{equation}
Unlike \(W(z)\), \eqref{FLT1}, the matrix-function \(S(z)\),
\eqref{TSM}, is an \textit{inner} function.

 The transformation in the Redheffer form \eqref{Redh}
 admits an interpretation in System Theory.
 We discuss this in more generality than is needed for our
 considerations, which are centered on the linear-fractional
 Schur transformation.

 Suppose that \(\text{LSDS}^{\,\text{I}}\) and
\(\text{LSDS}^{\,\text{II}}\) are two linear stationary dynamical
systems. In this section, we focus on the input-output mapping and
do not touch on considerations related to state spaces.

Let \(S(z):\,\mathcal{E}^{\text{I}}\to\mathcal{E}^{\text{I}}\) be
the transfer matrix-function of the system
\(\text{LSDS}^{\,\text{I}}\). Furthermore, let
\[\psi(z)=S(z)\varphi(z)\]
be the input-output mapping corresponding to the system
\(\text{LSDS}^{\,\text{I}}\), where
\(\varphi(z):\,\mathbb{D}\to\mathcal{E}^{\text{I}}\) is the input
signal and \(\psi(z):\,\mathbb{D}\to\mathcal{E}^{\text{I}}\) is
the output signal. Suppose now that the exterior space
\(\mathcal{E}^{\text{I}}\) of the system
\(\text{LSDS}^{\,\text{I}}\) is the orthogonal sum of the
subspaces \(\mathcal{E}^{\text{I}}_1\) and
\(\mathcal{E}^{\text{I}}_2\):
\begin{equation}%
\label{OrDe}
\mathcal{E}^{\text{I}}=\mathcal{E}^{\text{I}}_1\oplus{}\mathcal{E}^{\text{I}}_2\,.
\end{equation}%
Equation \eqref{OrDe} suggests that the input and output signals be decomposed as follows:
\begin{equation}
\varphi(z)=
\begin{bmatrix}
\varphi_1(z)\\[1.0ex]
\varphi_2(z)
\end{bmatrix}\,,\qquad
\psi(z)=
\begin{bmatrix}
\psi_1(z)\\[1.0ex]
\psi_2(z)
\end{bmatrix}\,,
\end{equation}
And furthermore that the matrix \(S(z)\) be decomposed accordingly:
\begin{equation}
S(z)=
\begin{bmatrix}
s_{11}(z)&s_{12}(z)\\[1.0ex]
s_{21}(z)&s_{22}(z)
\end{bmatrix}\,,
\end{equation}
So that
\begin{equation}
\label{EqS_I}
\begin{bmatrix}
\psi_1(z)\\[1.0ex]
\psi_2(z)
\end{bmatrix}=
\begin{bmatrix}
s_{11}(z)&s_{12}(z)\\[1.0ex]
s_{21}(z)&s_{22}(z)
\end{bmatrix}\,
\begin{bmatrix}
\varphi_1(z)\\[1.0ex]
\varphi_2(z)
\end{bmatrix}\,.
\end{equation}
 The system  \(\text{LSDS}^{\,\text{I}}\) can be considered as a
linear stationary dynamical system with two input channels,
corresponding to the input signals \(\varphi_1(z)\) and
\(\varphi_2(z)\), and two output channels, corresponding to the
output signals \(\psi_1(z)\) and \(\psi_2(z)\):\\ %

\noindent%
\centerline{\includegraphics*[scale=0.4]{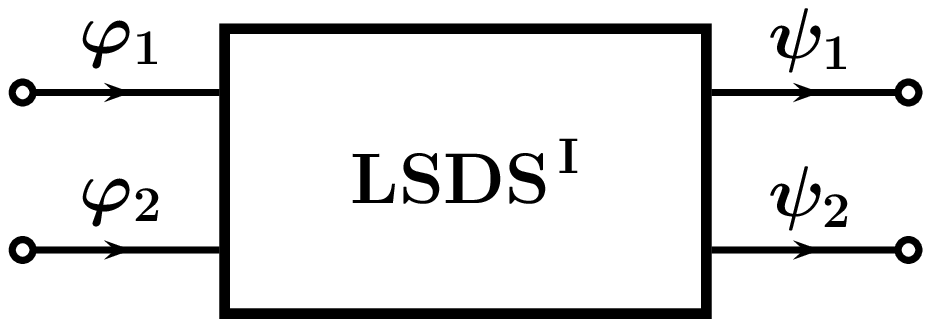}}
\centerline{Figure 1}%

\vspace{2.0ex}%
\noindent Let
\begin{equation}%
\label{EqS_II} \tau(z)=\omega(z)\sigma(z)
\end{equation}
  be the
input-output mapping corresponding to the system
\(\text{LSDS}^{\,\text{II}}\), where
\(\sigma(z):\,\mathbb{D}\to\mathcal{E}^{\text{II}}\) is the input
signal and \(\tau(z):\,\mathbb{D}\to\mathcal{E}^{\text{II}}\) is the
output signal.
 The system  \(\text{LSDS}^{\,\text{II}}\) can be considered as a
linear stationary dynamical system with one input channel,
corresponding to the input signal \(\sigma(z)\), and one output
channel, corresponding to the
output signal \(\tau(z)\): %

\vspace{2.0ex}%
\noindent%
\centerline{\includegraphics*[scale=0.4]{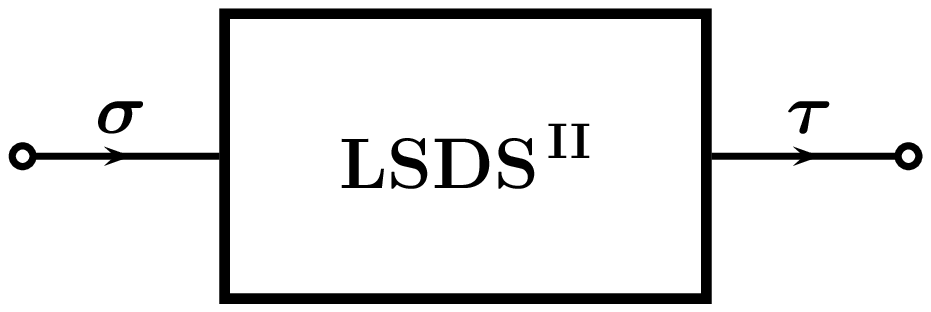}}
\centerline{Figure 2}%
\vspace{2.0ex}%

Suppose now that
\begin{equation}%
\label{CD}%
 \mathcal{E}^{\text{I}}_2=\mathcal{E}^{\text{II}}
\end{equation}%
This allows us to `link' the systems
 \(\text{LSDS}^{\,\text{I}}\) and  \(\text{LSDS}^{\,\text{II}}\).
 We connect the output channel of the system \(\text{LSDS}^{\,\text{II}}\)
 with the second \(\text{LSDS}^{\,\text{I}}\) input channel
 and the \(\text{LSDS}^{\,\text{II}}\) input channel
 with the second \(\text{LSDS}^{\,\text{I}}\) output channel, as shown in Figure~3.

 \noindent
 \centerline{\includegraphics*[scale=0.4]{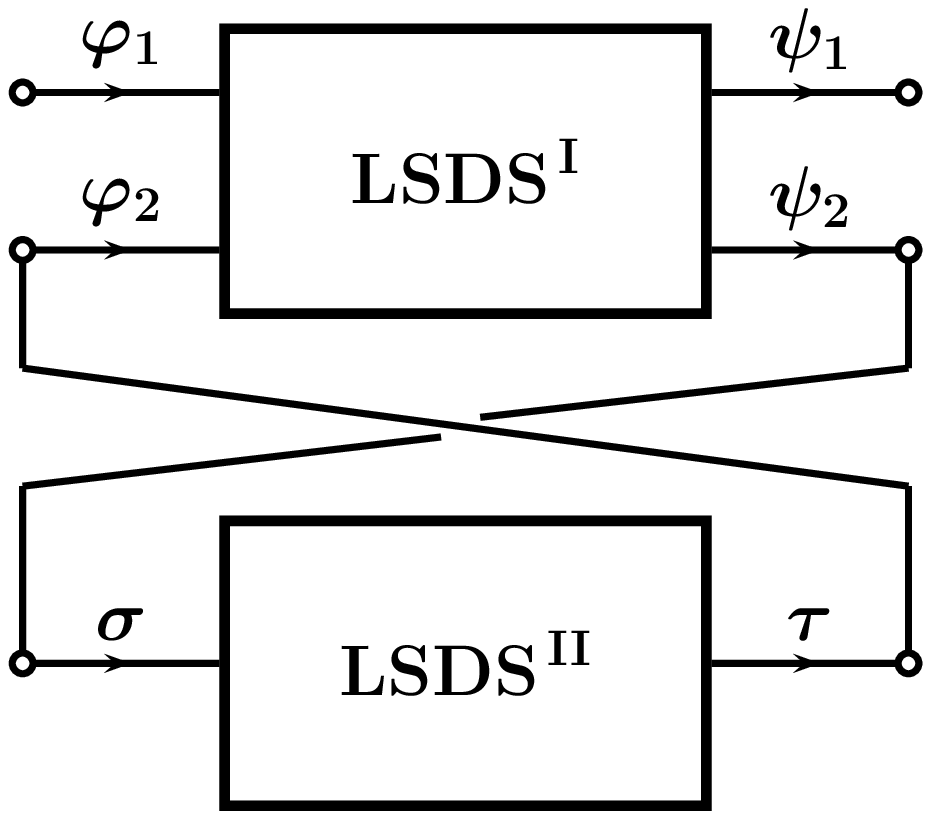}}
\centerline{Figure 3}%


 The resulting linear stationary dynamical system \(\text{LSDS}\)
 has exterior space \(\mathcal{E}^{\text{I}}_1\), input signal \(\varphi_1(z)\)
 and output signal \(\psi_1(z)\). The output signal \(\psi_1(z)\)
 is linearly dependent on the input signal \(\varphi_1(z)\):
\begin{equation}
\label{TMC}%
 \psi_1(z)=s(z)\,\varphi_1(z)\,,
\end{equation}
where \(s(z)\) is the transfer function for \(\text{LSDS}\).\\

We call \(\text{LSDS}\) the \textsf{Redheffer coupling} of
the systems \(\text{LSDS}^{\,\text{I}}\) and \(\text{LSDS}^{\,\text{II}}\).

\noindent%
We now look to express \(s(z)\) in terms of \(S(z)\) and \(\omega(z)\).
The above-described connection between the systems
\(\text{LSDS}^{\,\text{I}}\) and \(\text{LSDS}^{\,\text{II}}\) can
be formally expressed by means of the constraints
\begin{equation}
\label{Con}
 \varphi_2(z)=\tau(z),\qquad \psi_2(z)=\sigma(z)\,.
\end{equation}

Eliminating \(\varphi_2(z),\,\psi_2(z),\,\sigma(z),\,\tau(z)\)
from the system of linear equations \eqref{EqS_I}, \eqref{EqS_II}
and \eqref{Con}, we obtain the equation \eqref{TMC}, where \(s(z)\)
has the form
\begin{equation}
\label{RFLT}
s(z)=s_{11}(z)+s_{12}(z)\omega(z)\big(I-s_{22}(z)\omega(z)\big)^{-1}s_{21}(z)\,.
\end{equation}
We now turn our attention to the `\textsf{energy relation}' associated with
the linear fractional transformation \eqref{RFLT}:
\(\omega(z)\to{}s(z)\).

Equation \eqref{EqS_I} yields,
\[\varphi_1^\ast\varphi_1+\varphi_2^\ast\varphi_2-
\psi_1^\ast\psi_1-\psi_2^\ast\psi_2=
\begin{bmatrix}\varphi_1^\ast&\varphi_2^\ast\end{bmatrix}
(I-S^\ast{}S)\begin{bmatrix}\varphi_1\\\varphi_2\end{bmatrix}\,.
\]
Making the substitutions
\(\psi_1=s\varphi_1,\,\psi_2=\sigma\) and \(\varphi_2=\omega\sigma\), we
obtain
\begin{equation}
\label{ChD}
\varphi_1^\ast(1-s^\ast{}s)\varphi_1=\begin{bmatrix}\varphi_1^\ast&\varphi_2^\ast\end{bmatrix}
(I-S^\ast{}S)\begin{bmatrix}\varphi_1\\\varphi_2\end{bmatrix}+
\sigma^\ast(1-\omega^\ast\omega)\sigma\,,
\end{equation}
where
\begin{equation}
\label{ChDc}%
 \sigma=(1-s_{22}\omega)^{-1}s_{21}\varphi_1,\qquad
 \varphi_2=\omega(1-s_{22}\omega)^{-1}s_{21}\varphi_1\,.
\end{equation}
It follows from equation \eqref{ChD} that if \(I-S^\ast{}S\geq{}0\) and
 \(1-\omega^\ast\omega\geq{}0\), then \(1-s^\ast{}s\geq{}0\).
 If \(I-S^\ast{}S=0\) and \(1-\omega^\ast\omega=0\), then
 \(1-s^\ast{}s=0\). In particular, this brings us to:
 \begin{theorem}
 \label{PRT}
 Let \(S(z)\) and \(\omega(z)\) be rational inner matrix-functions.
 Furthermore, let \(s(z)\) be given by the Redheffer linear-fractional
 transform \eqref{RFLT}. Then \(s(z)\) is a rational inner matrix-function.
 \end{theorem}

 We note that the linear-fractional transform, in its
 classical form \eqref{tFLT}, is related to another kind of
 coupling. The relevant connection is shown in \textup{Figure 3}.

 \vspace*{2.0ex}%
\noindent%
\centerline{\includegraphics*[scale=0.4]{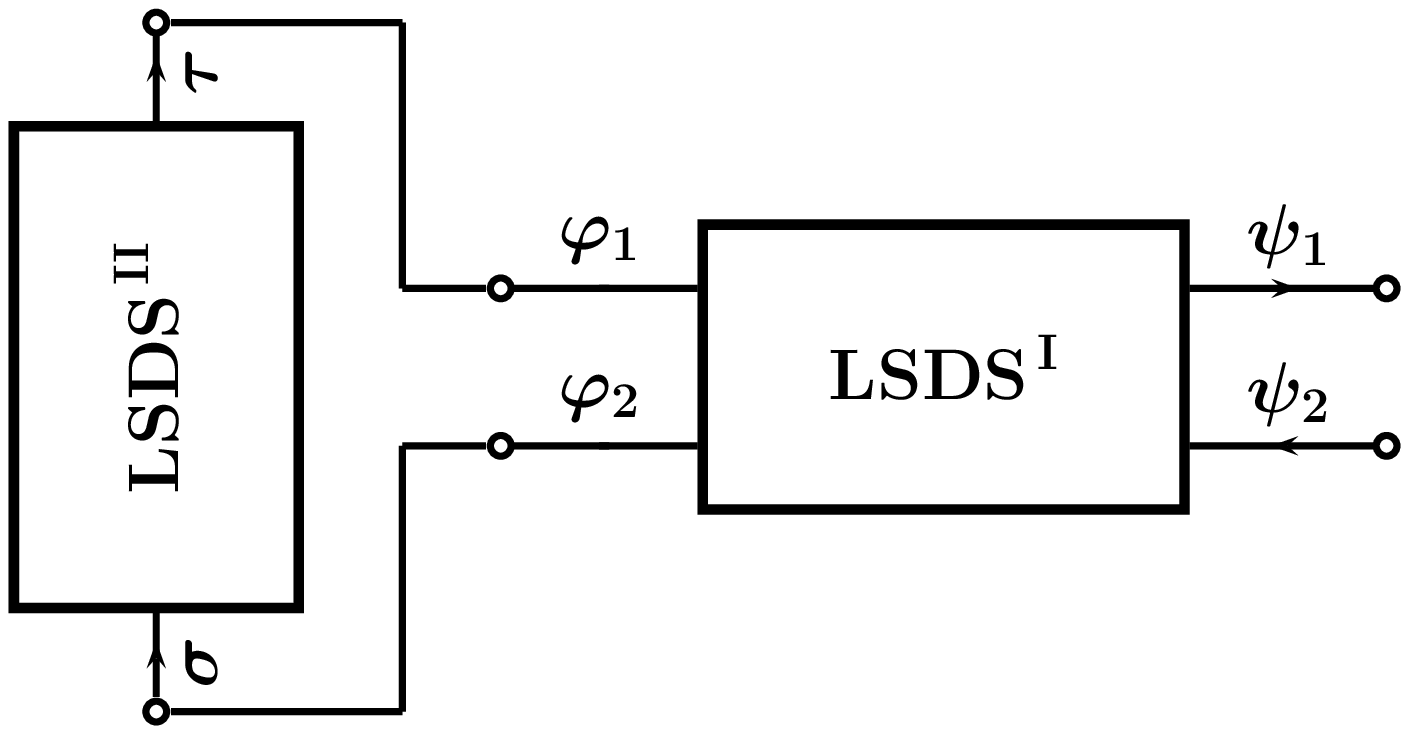}}
\centerline{\textup{Figure 3}}%
\vspace{1.0ex}%

\noindent \(\textup{LSDS}^{\textup{I}}\) has two input
channels with input signals \(\varphi_1(z)\) and
\(\varphi_2(z)\). \(\textup{LSDS}^{\textup{I}}\) also has two
output channels with output signals \(\psi_1(z)\) and \(\psi_2(z)\)
(with frequency representation).
The system \(\textup{LSDS}^{\textup{II}}\) has
one input channel with input signal \(\sigma(z)\) and
output signal \(\tau(z)\).
We connect the \(\textup{LSDS}^{\textup{II}}\) output channel
with the first input channel of the system \(\textup{LSDS}^{\textup{I}}\)
as well as the \(\textup{LSDS}^{\textup{II}}\) input channel with the
second input channel of the system \(\textup{LSDS}^{\textup{I}}\)
(We assume that the systems are compatible with respect to these connections,
i.e. that the appropriate subspaces coincide.)
 We consider the second output channel of the system
\(\textup{LSDS}^{\textup{I}}\) as the input channel of the new
coupled system \(\textup{LSDS}\) and the first output channel
of  \(\textup{LSDS}^{\textup{I}}\) as the output channel of
\(\textup{LSDS}\) (Shown in \textup{Figure 3}.)
Let %
\begin{math}W(z)=\Bigl[\begin{smallmatrix}
w_{11}(z)&w_{12}(z)\\[0.5ex]
w_{21}(z)&w_{21}(z)
\end{smallmatrix}\Bigr]\end{math} %
be the transfer matrix for \(\textup{LSDS}^{\textup{I}}\)
and \(\omega(z)\) be the transfer matrix for \(\textup{LSDS}^{\textup{II}}\):
\begin{equation*}
\begin{bmatrix}
\psi_1(z)\\[0.5ex]
\psi_2(z)
\end{bmatrix}=
\begin{bmatrix}
w_{11}(z)&w_{12}(z)\\[0.5ex]
w_{21}(z)& w_{22}(z)
\end{bmatrix}
\begin{bmatrix}
\varphi_1(z)\\[0.5ex]
\varphi_2(z)
\end{bmatrix}\,,\qquad
\tau(z)=\omega(z)\sigma(z)\,.
\end{equation*}
The link between the systems \(\textup{LSDS}^{\textup{I}}\) and
\(\textup{LSDS}^{\textup{II}}\), shown in \textup{Figure 3}, is
described by the constraints
\[\sigma(z)=\varphi_2(z), \qquad \tau(z)=\varphi_1(z)\,.\]
In which the input and the output signals of the system
\(\textup{LSDS}\) are denoted by \(\varphi(z)\) and \(\psi(z)\), respectively:
\[\varphi(z)=\psi_2(z)\,,\qquad \psi(z)=\psi_1(z)\,,\]
so that:
\[\psi(z)=s(z)\varphi(z)\,,\]
where %
\[s(z)=\big(w_{11}(z)\omega(z)+w_{12}(z)\big)\cdot%
\big(w_{21}(z)\omega(z)+w_{22}(z)\big)^{-1}\,.\]

\textsf{\large{Historical Remark.}} %
The coupling of input-output systems having four terminals,
considered in this section (See Figures 1-3), is sometimes called
\textit{cascade coupling}. This kind of coupling (as well as
related mathematical questions) was investigated by R.\,Redheffer
in \cite{Red1}\,-\,\cite{Red5}. Because of this, we use the name
\textit{Redheffer coupling}. Redheffer did not consider
questions related to cascade coupling of state space linear
systems. These questions were later addressed in \cite{HeBa} (Without any
reference to Redheffer.)

The results presented in \cite{HeBa} are more general than here needed.
We have tailored our approach to the theory of
Redheffer coupling in the next two sections to fit our needs.

\section{The Redheffer Coupling
 of Unitary Colligations. \label{RCOC}} \setcounter{equation}{0}

 As rational inner functions, \(S(z),\,\omega(z)\) and
 \(s(z)\) admit system representations as characteristic
 functions of the unitary operator colligations with
 colligation operators \(U^{\text{\,I}},\,U^{\text{\,II}}\) and \(U\),
 respectively. We now turn to the question of %
 \textsf{how we might express \(U\) in terms of the operators
  \(U^{\,I}\) and \(U^{\,II}\).} %

Our approach to this problem will be more general than is here
called for, our goal being to describe the colligations related to
Schur transformations. We assume that the unitary colligations
corresponding to the systems \(\text{LSDS}^{\,\text{I}}\) and
\(\text{LSDS}^{\,\text{II}}\) are given. We look to obtain the
unitary colligation corresponding to the system \(\text{LSDS} \),
the Redheffer coupling of the systems \(\text{LSDS}^{\,\text{I}}\)
and \(\text{LSDS}^{\,\text{II}}\). The system
\(\text{LSDS}^{\,\text{I}}\) is not assumed to be related to the
Schur transformation. \(\text{LSDS}^{\,\text{I}}\) and
\(\text{LSDS}^{\,\text{II}}\) can be generic systems. The only
condition imposed on these systems is that the exterior space
\(\mathcal{E}^{\textup{II}}\) of the system
\(\text{LSDS}^{\,\text{II}}\) is identified with the subspace
\(\mathcal{E}^{\textup{I}}_1\) of the exterior space
\(\mathcal{E}^{\textup{I}}\) belonging to
\(\text{LSDS}^{\,\text{I}}\). To avoid technical complications we
assume that the exterior and state spaces of the systems
\(\text{LSDS}^{\,\text{I}}\) and \(\text{LSDS}^{\,\text{II}}\) are
finite-dimensional.

 To simplify the notation, we denote the matrix entries
of the colligation operator \(U^{\text{I}}\), corresponding to the
system \(\text{LSDS}^{\,\text{I}}\), as follows
\begin{gather}
\label{U1}
U^{\text{\,I}}=
\begin{bmatrix}
a_{11}&a_{12}&\  b_{1}\\[0.5ex]
a_{21}&a_{22}&\  b_{2}\\[0.8ex]
c_{1}&c_{2}&\  d
\end{bmatrix}\,,
\intertext{where}
a_{p,q}:\mathcal{E}_p^{\,\text{I}}\to\mathcal{E}_q^{\,\text{I}},\ \quad
b_q: \mathcal{H}^{\,\text{I}}\to\mathcal{E}_q^{\,\text{I}},\ \quad
c_p:\mathcal{E}_p^{\,\text{I}}\to\mathcal{H}^{\,\text{I}},\ \quad
d: \mathcal{H}^{\,\text{I}}\to\mathcal{H}^{\,\text{I}}\,.\notag
\end{gather}
The matrix entries for the colligation operator \(U^{\text{II}}\),
corresponding to the system \(\text{LSDS}^{\,\text{II}}\),
are denoted as follows:
\begin{gather}
\label{U2}
U^{\text{\,II}}=
\begin{bmatrix}
\alpha &\  \beta\\[0.5ex]
\gamma&\delta
\end{bmatrix}\,,
\intertext{where}
\alpha:\mathcal{E}^{\,\text{II}}\to\mathcal{E}^{\,\text{II}},\ \quad
\beta: \mathcal{H}^{\,\text{II}}\to\mathcal{E}^{\,\text{II}},\ \quad
\gamma:\mathcal{E}^{\,\text{II}}\to\mathcal{H}^{\,\text{II}},\ \quad
\delta: \mathcal{H}^{\,\text{II}}\to\mathcal{H}^{\,\text{II}}\,.\notag
\end{gather}
The linear equations describing the dynamics of the
system \(\text{LSDS}^{\,\text{I}}\) are
\begin{gather}
\label{L1}
\begin{bmatrix}
\psi_1(z)\\[0.5ex]%
\psi_2(z)\\[0.8ex]%
z^{-1}h(z)
\end{bmatrix}
=
\begin{bmatrix}
a_{11}&a_{12}&\  b_{1}\\[0.5ex]
a_{21}&a_{22}&\  b_{2}\\[0.8ex]
c_{1}&c_{2}&\  d
\end{bmatrix}
\begin{bmatrix}
\varphi_1(z)\\[0.5ex]%
\varphi_2(z)\\[0.8ex]%
h(z)
\end{bmatrix}
\,,
\intertext{where}
\varphi(z)=
\begin{bmatrix}
\varphi_1(z)
\\[0.5ex]
\varphi_2(z)
\end{bmatrix}\,,\
\ \psi(z)=
\begin{bmatrix}
\psi_1(z)
\\[0.5ex]
\psi_2(z)
\end{bmatrix}\, \ \ \text{and} \  \ %
h(z)
\notag
\end{gather}
are, respectively, the input signal, the output signal and the
inner state signal corresponding to the system \(\text{LSDS}^{\,\text{I}}\).

The linear equations describing the dynamics of the system
\(\text{LSDS}^{\,\text{II}}\) are
\begin{equation}
\label{L2}
\begin{bmatrix}
\tau(z)\\[0.5ex]%
z^{-1}l(z)
\end{bmatrix}
=
\begin{bmatrix}
\alpha&\  \beta \\[0.5ex]
\gamma&\  \delta
\end{bmatrix}
\begin{bmatrix}
\sigma(z)\\[0.5ex]%
l(z)
\end{bmatrix}
\,,
\end{equation}
where
\(\sigma(z),\
\tau(z)\ \text{and} \  %
l(z)\) are, respectively, the input signal, output signal and the
interior state signal corresponding to the system
\(\text{LSDS}^{\,\text{II}}\).

The constraints
\begin{equation}
\label{Res}
\tau(z)=\varphi_2(z)\,,\ \ \sigma(z)=\psi_2(z)
\end{equation}
correspond to the Redheffer coupling of the systems \(\text{LSDS}^{\,\text{I}}\) and
\(\text{LSDS}^{\,\text{II}}\).

We now aim to eliminate the variables \(\varphi_2(z),\,\psi_2(z),\,\sigma(z),\,\tau(z)\)
from the systems \eqref{L1}, \eqref{L2}, \eqref{Res}. To this end,
we substitute the expressions \(\alpha\sigma(z)+\beta{}l(z)\) and \(\sigma(z)\)
for the variables \(\varphi_2(z)\) and \(\psi_2(z)\) into the equation %
\[\psi_{2}(z)=a_{21}\varphi_1(z)+a_{22}\varphi_2(z)+\beta{}h(z)\,.\]
With this we can express \(\sigma(z)\) in terms of \(\varphi_1(z), \,h(z)\) and \(l(z)\):
\begin{multline}
\label{El}
\sigma(z)=\\
=(1-a_{22}\,\alpha)^{-1}\,a_{21}\,\varphi_{1}(z)+(1-a_{22}\,\alpha)^{-1}\,b_2\,h(z)
+(1-a_{22}\,\alpha)^{-1}\,\beta{}\,l(z)\,.
\end{multline}
Substituting this expressions for \(\sigma\) into \eqref{L1}, \eqref{L2}, \eqref{Res},
we obtain
\begin{equation}
\label{L}
\begin{bmatrix}
\psi_1(z)\\[0.5ex]
z^{-1}h(z)\\[0.5ex]
z^{-1}l(z)
\end{bmatrix}%
=%
\begin{bmatrix}
\ \ A \ &B_1&B_2\\[0.8ex]
C_1&D_{11}&D_{12}\\
C_2&D_{21}&D_{21}
\end{bmatrix}%
\begin{bmatrix}
\varphi_1(z)\\[0.5ex]
h(z)\\[0.5ex]
l(z)
\end{bmatrix}\,,%
\end{equation}
where %
 \begin{multline*}
 A:\,\mathcal{E}^{\text{\,I}}_1\to\mathcal{E}^{\text{\,I}}_1,
 \ B_1:\,\mathcal{H}^{\,\text{I}}\to\mathcal{E}^{\text{\,I}}_1,\,
 B_2:\,\mathcal{H}^{\,\text{II}}\to\mathcal{E}^{\text{\,I}}_1,\,\\
\ C_1:\,\mathcal{E}^{\text{\,I}}_1\to\mathcal{H}^{\,\text{I}},\,
 C_2:\,\mathcal{H}^{\,\text{II}}\to\mathcal{E}^{\,\text{I}}_1,\,\hspace{13.0ex}\\
 D_{11}:\mathcal{H}^{\,\text{I}}\to\mathcal{H}^{\,\text{I}}\,,
 \ \
 D_{12}:\mathcal{H}^{\,\text{II}}\to\mathcal{H}^{\,\text{I}}\,,\ \
 D_{21}:\mathcal{H}^{\,\text{I}}\to\mathcal{H}^{\,\text{II}}\,,
 \ \
 D_{22}:\mathcal{H}^{\,\text{II}}\to\mathcal{H}^{\,\text{II}}\,.
 \end{multline*}
The matrix
\begin{equation}%
\label{U}%
U=\begin{bmatrix}
\ \ A \ &B_1&B_2\\[0.8ex]
C_1&D_{11}&D_{12}\\
C_2&D_{21}&D_{21}
\end{bmatrix}%
\end{equation}
can be expressed using the entries of the matrices
\(U^{\,\text{\,I}}\), \eqref{U1}, and  \(U^{\,\text{\,II}}\),
\eqref{U2}, as follows:
\begin{equation}
\label{RedhC} U=%
\begin{bmatrix}
a_{11}&b_1&a_{12}\,\beta\\[0.5ex]
c_1&d&c_{2}\,\beta\\[0.5ex]
0&0&\delta
\end{bmatrix}
+
\begin{bmatrix}
a_{12}\,\alpha\\[0.5ex]
c_2\,\alpha\\[0.5ex]
\gamma
\end{bmatrix}
\cdot\big(1-a_{22}\,\alpha\big)^{-1}\cdot
\begin{bmatrix}
\,a_{21}\,&b_2\,&a_{22}\,\beta\,
\end{bmatrix}
\,.
\end{equation}
The operator \(U\) is called the \textsf{Redheffer product} of
the operators \(U_1\) and \(U_2\).

We again turn our attention to the `\textsf{energy relation}' associated
with the operators \(U^{\,\text{I}}\), \(U^{\,\text{II}}\) and \(\,U\).
Suppose that \(U^{\,\text{I}}\) and \(U^{\,\text{II}}\) are unitary. Let \linebreak
\(\varphi_1\in\mathcal{E}^{\text{I}}_1,\ \varphi_2\in\mathcal{E}^{\text{I}}_2, \
h\in\mathcal{H}^{\text{I}}\) \
\(\sigma\in{}\mathcal{E}^{\text{II}}\) and \(l\in\mathcal{H}^{\text{II}}\)
be arbitrary vectors. If
\(\psi_1\in\mathcal{E}^{\text{I}}_1,\ \psi_2\in\mathcal{E}^{\text{I}}_2, \
k\in\mathcal{H}^{\text{I}}\), \
\(\tau\in\mathcal{E}^{\text{II}}\) and \(m\in\mathcal{H}^{\text{II}}\)
are defined by the equalities
\begin{equation*}
\begin{bmatrix}
\psi_1\\[0.5ex]
\psi_2\\[0.5ex]
k
\end{bmatrix}=\, \text{\large \(U^{\text{I}}\)}
\begin{bmatrix}
\varphi_1\\[0.5ex]
\varphi_2\\[0.5ex]
h
\end{bmatrix}%
\,,\quad%
\begin{bmatrix}
\tau\\[0.5ex]
 m
\end{bmatrix}=\, \text{\large \(U^{\text{II}}\)}
\begin{bmatrix}
\sigma\\[0.5ex]
 l
\end{bmatrix}\,,
\end{equation*}
then %
\begin{equation}
\label{ER1}%
 ||\psi_1||^2+||\psi_2||^2+||k||^2=||\varphi_1||^2+||\varphi_2||^2+||h||^2\,,
\end{equation}
and
\begin{equation}
\label{ER2}%
 ||\tau||^2+||m||^2=||\sigma||^2+||l||^2\,.
\end{equation}
For arbitrary \(\varphi_1,\,h,\,l\) and
\begin{equation}
\label{ressig}
\sigma=(1-a_{22}\,\alpha)^{-1}\,a_{21}\,\varphi_{1}+(1-a_{22}\,\alpha)^{-1}\,b_2\,
+(1-a_{22}\,\alpha)^{-1}\,\beta{}\,,
\end{equation}
\begin{equation}
\label{resphi2}%
\varphi_2=\alpha\big((1-a_{22}\,\alpha)^{-1}\,a_{21}\,\varphi_{1}+
(1-a_{22}\,\alpha)^{-1}\,b_2\,
+(1-a_{22}\,\alpha)^{-1}\,\beta\big)+\beta{}l\, ,
\end{equation}
it follows that
\[\psi_2=\sigma,\qquad\tau=\varphi_2\,,\]
and
\begin{equation}
\label{ER}%
 ||\psi_1||^2+||l||^2+||m||^2= ||\varphi_1||^2+||h||^2+||l||^2\,.
\end{equation}
According to the definition of the operator \(U\),
\begin{equation}
\label{UDet}%
\begin{bmatrix}
\psi_1\\[0.5ex]
k\\[0.5ex]
m
\end{bmatrix}=\, \text{\large \(U\)}
\begin{bmatrix}
\varphi_1\\[0.5ex]
h\\[0.5ex]
l
\end{bmatrix}\,.%
\end{equation}
Since \(\varphi_1,\,h,\,l\) are arbitrary, equality \eqref{ER}
means that \(U\) is unitary. This operator, partitioned into
blocks according to \eqref{U}, is related to the unitary
colligation \((\mathcal{E},\mathcal{H},U)\), where
\(\mathcal{E}=\mathcal{E}^\text{I}\),
 \(\mathcal{H}=\mathcal{H}^\textup{I}\oplus\mathcal{H}^\textup{II}\).
\vspace{0.5ex}
\begin{definition}
\textit{The colligation \((\mathcal{E},\mathcal{H},U)\) is called the
\textsf{Redheffer coupling of the colligations}
\((\mathcal{E}^\textup{I},\mathcal{H}^\textup{I},U^\textup{I})\)
and \((\mathcal{E}^\textup{II}, \ \mathcal{H}^\textup{II},U^\textup{II})\).}
\end{definition}
\vspace{0.5ex}
\begin{theorem}
\label{CFCC}%
 Let
\begin{math}%
S(z)=\Bigl[\begin{smallmatrix}s_{11}(z)&s_{12}(z)\\
s_{21}(z)&s_{22}(z)\end{smallmatrix}\Bigr]%
\end{math},  \(\omega(z)\) and \(s(z)\) be the characteristic functions of
the colligations
\((\mathcal{E}^\text{I},\mathcal{H}^\text{I},U^\text{I})\),
\((\mathcal{E}^\text{II},\mathcal{H}^\text{II},U^\text{II})\) and
their Redheffer coupling \((\mathcal{E},\mathcal{H},U)\),
respectively:
\begin{equation}
\begin{bmatrix}s_{11}(z)&s_{12}(z)\\
s_{21}(z)&s_{22}(z)\end{bmatrix}=\\
\begin{bmatrix}a_{11}&a_{12}\\
a_{21}&a_{22} \end{bmatrix}+z
\begin{bmatrix}b_1\\ b_2\end{bmatrix}
(I-zd)^{-1}
\begin{bmatrix}c_1& c_2\end{bmatrix}\,,
\end{equation}
\begin{equation}
\omega(z)=\alpha+z\beta(1-z\delta)^{-1}\gamma\,,
\end{equation}
\begin{equation}
s(z)=A+z\begin{bmatrix}B_1\\ B_2
\end{bmatrix}
\left(\begin{bmatrix}I_{\mathcal{H}^\textup{I}}&0\\
0&I_{\mathcal{H}^\textup{II}}
\end{bmatrix}-
z\begin{bmatrix}D_{11}&D_{12}\\
D_{21}&D_{22}
\end{bmatrix}
\right)^{-1}
\begin{bmatrix}C_1&C_2
\end{bmatrix}\,.
\end{equation}
(The notation for the entries of the matrices \(U^{\text{I}}\),
\(U^{\text{II}}\) and \(U\) is taken from \eqref{U1}, \eqref{U2}
 and \eqref{U}, respectively.)

Then
\begin{equation}
s(z)=s_{11}(z)+s_{12}(z)\omega(z)(I-s_{22}(z)\omega(z))^{-1}s_{21}(z)\,.
\end{equation}
\end{theorem}
\section{The Inverse Schur Transformation and\\
Redheffer Couplings of Colligations.\label{SASRoc}}
\setcounter{equation}{0}
We now focus again on the linear-fractional transformation \eqref{FLT}
in the Redheffer form \eqref{Redh}, where \(\omega(z)\) is a rational inner
matrix-function of degree \(n-1\), so that \(s(z)\) is a rational inner
matrix-function of degree \(n\).

 The function \(S(z)\), which appears in \eqref{TSM} is a rational inner function.
It is a characteristic matrix-function for the
operator colligation \((\mathcal{E}^{\text{I}},\
\mathcal{H}^\text{I},\,U^{\text{I}})\), which we now describe.

The outer space \(\mathcal{E}^{\text{I}}\) is two-dimensional. We
identify  \(\mathcal{E}^{\text{I}}\) with
 \(\mathbb{C}^2\). The space \(\mathcal{E}^{\text{I}}\) is
 considered as the orthogonal sum \(\mathcal{E}^{\text{I}}=
 \mathcal{E}^{\text{I}}_1\oplus\mathcal{E}^{\text{I}}_2\), where
\(\mathcal{E}^{\text{I}}_1\) is identified with \(\mathbb{C}\)
and \(\mathcal{E}^{\text{I}}_2\) is identified with
\(\mathbb{C}\). The orthogonal decomposition
\(\mathcal{E}^{\text{I}}=
 \mathcal{E}^{\text{I}}_1\oplus\mathcal{E}^{\text{I}}_2\)
is thus the canonical decomposition
\(\mathbb{C}^2=\mathbb{C}\oplus\mathbb{C}\). The inner space
\(\mathcal{H}^{\text{I}}\) is one-dimensional. We identify
\(\mathcal{H}^{\text{I}}\) with
 \(\mathbb{C}^1\).
 The colligation operator \(U^{\text{I}}\) is defined by the
unitary \(3\times{}3=(2+1)\times(2+1)\)-matrix considered as an
operator acting in
\(\mathbb{C}^3=\mathbb{C}^2\oplus\mathbb{C}^1\):
\begin{equation}
\label{Uu1}
U^\text{\,I}=%
\begin{bmatrix}
\ A^{\text{\,I}}&B^{\text{\,I}}\ \\[0.5ex]
\ C^{\text{\,I}}&D^{\text{\,I}}\
\end{bmatrix}
\end{equation}
with
\begin{multline*}
A^{\text{\,I}}=
\begin{bmatrix}
\ s_0\ & \ 0 \  \\[1.0ex]
\ (1-|s_0|^2)^{1/2} \ & \ 0 \
\end{bmatrix},\
B^{\text{I}}=
\begin{bmatrix}
\ (1-|s_0|^2)^{1/2} \ \\[1.0ex] %
\ -\overline{s_0},\
\end{bmatrix},\ \\[0.7ex]
C^{\text{I}}=
\begin{bmatrix}
\ 0 &1 \
\end{bmatrix},\ \
 D^{\text{I}}=
\begin{bmatrix}
\ 0\
\end{bmatrix}\,.
\end{multline*}
The characteristic function of the colligation
\((\mathcal{E}^{\text{I}},\mathcal{H}^{\text{I}},U^{\text{I}})\)
is the matrix-function \(S(z)\) of the form \eqref{TSM}:
\begin{equation}
\label{SRS1}
\begin{bmatrix}
s_0&\ \ z\,(1-|s_0|^2)^{1/2}\\[2.0ex]
(1-|s_0|^2)^{1/2}&-z\,\overline{s_0}
\end{bmatrix}=
A^{\,\text{I}}+zB^{\text{\,I}}(I-zD^{\,\text{I}})^{-1}C^{\text{\,I}}\,.
\end{equation}

 The rational inner function
\(\omega(z)\) of degree \(n-1\) is the characteristic function of
the colligation
\((\mathcal{E}^{\text{II}},\mathcal{H}^{\text{II}},U^{\text{II}})\).
The outer space \(\mathcal{E}^{\text{II}}\) is one-dimensional and
is identified with \(\mathbb{C}\) and the inner space
\(\mathcal{H}^{\text{II}}\) is \((n-1)\)-dimensional and is
identified with \(\mathbb{C}^{n-1}\). The colligation operator
\(U^{\text{II}}\) thus acts in
\(\mathbb{C}^{n}=\mathbb{C}\oplus\mathbb{C}^{n-1}\). We identify
the operator \(U^{\text{II}}\) with its matrix in the canonical
basis of \(\mathbb{C}^{n}\):
\begin{equation}
\label{CROm}%
 U^{\text{II}}=%
 \begin{bmatrix}
 \,\alpha&\beta\,\\[0.5ex]
 \,\gamma&\delta\,
 \end{bmatrix}\,,
\end{equation}
where \[\alpha\in\mathfrak{M}_{1\times{}1}, \
\beta\in\mathfrak{M}_{1\times{}(n-1)}, \
\gamma\in\mathfrak{M}_{(n-1)\times{}1}, \
\delta\in\mathfrak{M}_{(n-1)\times{}(n-1)}.\]
\(\alpha\) is simply a complex number. The matrix \(U^{\text{\,II}}\) is unitary.
The system rep\-resentation of the function \(\omega(z)\) is given by:
\begin{equation}
\label{SROm}%
 \omega(z)=\alpha+\beta(1-z\delta)^{-1}\gamma\,.
\end{equation}
In particular,
\begin{equation}
\label{SPOm}%
 \omega(0)=\alpha\,.
\end{equation}
The function
\begin{equation}
\label{sSA}
s(z)=\frac{s_0+z\omega(z)}{1+\,z\omega(z)\overline{s_0}}\,,
\end{equation}
written as a Redheffer fractional-linear transform, takes the form:
\begin{equation}
\label{sSAR} s(z)=s_0+z(\,1-|s_0|^{2}\,)^{1/2}\,
\omega(z)\,(\,1+z\omega(z)\overline{s_0}\,)^{-1}\,(1-|s_0|^{2})^{1/2}\,,
\end{equation}
and admits a system realization by means of the operator colligation
\((\mathcal{E},\mathcal{H},U)\), where \((\mathcal{E},\mathcal{H},U)\)
 is the Redheffer coupling of the colligations
 \((\mathcal{E}^{\text{I}},\mathcal{H}^{\text{I}},U^{\text{I}})\),
representing the function \(S(z)\), and
\((\mathcal{E}^{\text{II}},\mathcal{H}^{\text{II}},U^{\text{II}})\),
representing the function \(\omega(z)\).

 Clearly, \(\mathcal{E}=\mathbb{C}\) and \(\mathcal{H}=\mathbb{C}^n\).
 \(U\) is the Redheffer coupling of the matrices
 \(U^{\text{I}}\) and \(U^{\text{II}}\).
 Applying formula \eqref{RedhC} to \(U^{\,\text{I}}\) and
\(U^{\,\text{II}}\), we obtain:
\begin{equation}
\label{SCo}%
 U=
\begin{bmatrix}
s_0\ & \ (1-|s_0|^2)^{1/2}  & \ 0_{1\times(n-1)} \ \\[0.7ex]
\alpha\,(1-|s_0|^2)^{1/2} \ & \ -\alpha\,\overline{s_0} \ & \
\beta \ \\[0.7ex]
\gamma\,(1-|s_0|^2)^{1/2} \ & \ -\gamma\,\overline{s_0} \ & \
\delta \
\end{bmatrix}
\,,
\end{equation}
so that \(U\) takes the form:
\begin{equation}
\label{SpStU}%
 U=%
 \begin{bmatrix}
 A &B\ \\[0.5ex]
C    & D
\end{bmatrix},
\end{equation}
where%
\begin{multline}        \label{SRsE}
\qquad A=s_0,  \qquad
B=\begin{bmatrix}
\ (1-|s_0|^2)^{1/2}& 0_{1\times(n-1)} \
\end{bmatrix},\,\\[0.7ex]
C=\begin{bmatrix}
 \ \alpha\,(1-|s_0|^2)^{1/2}\ \\[0.5ex]
\gamma\,(1-|s_0|^2)^{1/2} \
\end{bmatrix}\,, \qquad
D=\begin{bmatrix}%
 -\alpha\,\overline{s_0}&\beta \ \\[0.5ex]
 -\gamma\,\overline{s_0}&\delta
\end{bmatrix}\,,\\[0.7ex]
A\in\mathfrak{M}_{1\times{}1}\,,\
B\in\mathfrak{M}_{1\times{}n}\,,\ B\in\mathfrak{M}_{n\times{}1}, \
 D\in\mathfrak{M}_{n\times{}n}\,.
\end{multline}
Clearly, \(U\) in \eqref{SCo}-\eqref{SpStU} can be expressed as follows:
\begin{equation}
\label{SCoo}%
 U=
\begin{bmatrix}
1  & \ 0 \ & 0_{1\times(n-1)} \\[0.7ex]
0   & \ \alpha & \beta\\[0.5ex]
0_{(n-1)\times{}1}   & \ \gamma & \delta
\end{bmatrix}
\begin{bmatrix}
s_0 & \ (1-|s_0|^2)^{1/2} \ &\  0_{1\times{}(n-1)} \\[0.5ex]
(1-|s_0|^2)^{1/2}& \!\!-\overline{s_0} \ &  0_{1\times(n-1)}\\[0.8ex]
0_{(n-1)\times{}1}   &  0_{(n-1)\times{}1} &  1_{(n-1)\times(n-1)}
\end{bmatrix}
\!.
\end{equation}

\vspace{2.0ex}
Applying Theorem \ref{CFCC} to the Redheffer coupling of the
colligations \(U^{\textup{I}}\), \eqref{Uu1}, and
\(U^{\textup{II}}\), \eqref{Uu1}, yields:
\begin{theorem}
\label{SRIST}%
Let \(\omega(z)\) be an rational inner matrix-function of degree
\(n-1\) and let
\begin{equation}
\label{SROmT}
 \begin{bmatrix}
 \,\alpha&\beta\,\\[0.5ex]
 \,\gamma&\delta\,
 \end{bmatrix}, \
 \alpha\in\mathfrak{M}_{1\times{}1},\,
\beta\in\mathfrak{M}_{1\times{}(n-1)},\,
\gamma\in\mathfrak{M}_{(n-1)\times{}1},\,
\delta\in\mathfrak{M}_{(n-1)\times{}(n-1)}\,,
\end{equation}
be a unitary matrix so that the system representation
\eqref{SROm} for \(\omega(z)\) holds. Let \(s_0\) be a
complex number with \(|s_0|<1\). Let the function \(s(z)\) be defined
as the inverse Schur transform \eqref{sSA} (using \(s_0\) and \(\omega(z)\))
and let the matrix \(U\),
\begin{equation}
\label{SRUuu} U=
 \begin{bmatrix}
 \,A&B\,\\[0.5ex]
 \,C&D\,
 \end{bmatrix}, \
 A\in\mathfrak{M}_{1\times{}1},\,
B\in\mathfrak{M}_{1\times{}(n)},\,
C\in\mathfrak{M}_{(n)\times{}1},\,
D\in\mathfrak{M}_{(n)\times{}(n)}\,,
\end{equation}
be defined by equation \eqref{SCoo}.

\(U\) is then unitary and yields the system representation of \(s(z)\):
\begin{equation}
\label{SRs}%
s(z)=A+zB(I-zD)^{-1}C\,.
\end{equation}
\end{theorem}

\centerline{\textsf{Unitary Equivalence Freedom.\\[1.0ex]}}

The same function \(s(z)\), for which we earlier found a representation
using the matrix \(U\) in \eqref{SCoo}, can also be represented with the
help of a matrix having the form:
\begin{equation}%
\label{DeFr}
U^{\,V}=\begin{bmatrix}
 1 &0_{n\times{}n}\ \\[0.5ex]
0_{n\times{}n}    & V^\ast\end{bmatrix}\,%
U\,%
\begin{bmatrix}
 1 &0_{n\times{}n}\ \\[0.5ex]
0_{n\times{}n}   & V
\end{bmatrix}\,,
\end{equation}
where \(V\in\mathfrak{M}_{n\times{}n}\) is a unitary matrix.

 The matrix representing \(s(z)\) and which, furthermore, \textit{appears as the
Redheffer coupling matrix} for the matrices representing \(S(z)\) and \(\omega(z)\),
can be considered to have fewer `degrees of freedom' than matrices of the form \eqref{DeFr}.
The degree of freedom for the Redheffer coupling matrix is derived from this same property
in the Redheffer coupled matrices. The more general form of the matrix, which represents the
\(2\times{}2\)-matrix-function \(S(z)\), is the `transformed' matrix:
\begin{gather}
\label{Uuu1}
U^{\,\text{I},\,\varepsilon}=%
\begin{bmatrix}
1_{\,2\times{}2}&0_{\,2\times{}1}\\
0_{\,1\times{}2}&\overline{\varepsilon}
\end{bmatrix}
U^\text{\,I}%
\begin{bmatrix}
1_{\,2\times{}2}&0_{\,2\times{}1}\\
0_{\,1\times{}2}&\varepsilon
\end{bmatrix},%
\intertext{i.e.}%
U^{\,\text{I},\,\varepsilon}=
\begin{bmatrix}
\ A^{\,\text{I}}&B^{\text{\,I}}\varepsilon{}\ \\[0.5ex]
\ \overline{\varepsilon} \,C^{\text{\,I}}&D^{\text{\,I}}\
\end{bmatrix}\,,
\end{gather}
where \(U^{\,\text{I}}\) is the matrix from \eqref{Uu1} and \(\varepsilon\)
is an arbitrary unimodular complex number.
A more general form of the colligation matrix representing the function
\(\omega(z)\) is given by:
\begin{gather}
\label{Uu2} %
U^{{\,\text{II}},\,v}=
\begin{bmatrix}
1&0_{1\times{}(n-1)}\\
0_{(n-1)\times{}1}&v^\ast
\end{bmatrix}U^{\text{\,II}}
\begin{bmatrix}
1&0_{1\times{}(n-1)}\\
0_{(n-1)\times{}1}&v
\end{bmatrix}\,,\\
\intertext{i.e.}
U^{{\,\text{II}},\,v}=
\begin{bmatrix}
\alpha &\  \beta^{\,v}\\[0.5ex]
\gamma^{\,v}&\delta^{\,v}
\end{bmatrix}\,,
\end{gather}
where \(U^{\text{\,II}}\), \ \eqref{CROm}\,, is some \(n\times{}n\) unitary
colligation matrix representing the function \(\omega(z)\),
\begin{equation}
\beta^{\,v}=\beta\,v,\qquad \gamma^{\,v}=v^{\,\ast}\gamma,
\qquad \delta^{\,v}=v^{\,\ast}\delta\,v\,,
\end{equation}
 and
\(v\) is an arbitrary unitary \((n-1)\times(n-1)\)-matrix.
Applying formula \eqref{RedhC} to the matrices \(U^{\,\text{I},\,\varepsilon}\)
 and  \(U^{\,\text{II},\,v}\), we obtain the Redheffer coupling matrix:
\begin{equation}
\label{SCoMg}%
 U^{\,\varepsilon,\,v}=
\begin{bmatrix}
s_0\ & \ \varepsilon\,(1-|s_0|^2)^{1/2}  & \ 0_{1\times(n-1)} \ \\[0.7ex]
\alpha\,\overline{\varepsilon}\,(1-|s_0|^2)^{1/2} \ & \ -\alpha\,\overline{s_0} \ & \
\,\overline{\varepsilon}\,\beta^{\,v} \ \\[0.7ex]
\gamma^{\,v}\,(1-|s_0|^2)^{1/2} \ & \ -\gamma\,\varepsilon\,\overline{s_0} \ & \
\delta^{\,v} \
\end{bmatrix}
\,.
\end{equation}
Clearly,
\begin{multline}
\label{SCoot}%
 U^{\,\varepsilon,v}=
\begin{bmatrix}
1  & \ \ \ \ 0 \ & 0_{1\times(n-1)} \\[1.4ex]
0   & \ \ \ \alpha &\overline{\varepsilon}\,\beta^{\,v}\\[0.5ex]
0_{(n-1)\times{}1}   & \ \ \ \ \gamma^{\,v}\varepsilon&\delta^{\,v}
\end{bmatrix}\times %
\\
\times%
\begin{bmatrix}
s_0 & \ \varepsilon\,(1-|s_0|^2)^{1/2} \ &\  0_{1\times{}(n-1)} \\[0.5ex]
\overline{\varepsilon}\,(1-|s_0|^2)^{1/2}& \!\!-\overline{s_0} \ &  0_{1\times(n-1)}\\[0.8ex]
0_{(n-1)\times{}1}   &  0_{(n-1)\times{}1} &  1_{(n-1)\times(n-1)}
\end{bmatrix}
\!,
\end{multline}
and finally
\begin{equation}%
\label{DeFrv}
U^{\,\varepsilon,\,v}=\begin{bmatrix}
 1 &0_{n\times{}n}\ \\[0.5ex]
0_{n\times{}n}    & V^\ast_{\,\varepsilon,\,v}\end{bmatrix}\,%
U\,%
\begin{bmatrix}
 1 &0_{n\times{}n}\ \\[0.5ex]
0_{n\times{}n}   & V_{\,\varepsilon,\,v}
\end{bmatrix}\,,
\end{equation}
where
\begin{equation}
\label{MDMv}
V_{\,\varepsilon,v}=%
\begin{bmatrix}
\varepsilon&0_{1\times{}(n-1)}\\[0.5ex]
0_{(n-1)\times{}1}&v
\end{bmatrix}\,,
\end{equation}%
\(\varepsilon\) is an arbitrary unimodular
complex number and \(v\) is an arbitrary\\
unitary \((n-1)\times(n-1)\)-matrix (\(\varepsilon\) and \(v\) are
the same as in \eqref{Uu2}).

Comparing formulas \eqref{DeFr} and
\eqref{DeFrv}-\eqref{MDMv}, we see that the matrices
\(U^{\,\varepsilon,\,v}\) which come from Redheffer coupling of
the  matrices representing \(S(z)\) and \(\omega(z)\) are special.
The distinguishing feature of the matrices \(U^{\,\varepsilon,\,v}\)
can be summarized as follows:\\

\noindent
\hspace{2.0ex}
\begin{minipage}[h]{0.95\linewidth}
 \textit{Among all of the \((n+1)\times(n+1)\)-matrices
\begin{math}%
U^V=%
\begin{bmatrix}
\label{GenU}
s_0&\,B^V\\[0.5ex]
C^V&\,D^V
\end{bmatrix}
\end{math}
of the form \eqref{DeFr}, it is precisely those for which the block-matrix
entry \(B^V\) takes the form
 \begin{equation*}
B^V=%
\begin{bmatrix}
\varepsilon\,(1-|s_0|^2)^{1/2}  & \ 0_{1\times(n-1)}
\end{bmatrix}\,,
 \end{equation*}
where \(\varepsilon\) is an arbitrary unimodular complex number,
that can be expressed as in \eqref {DeFrv}-\eqref{MDMv}.}
\end{minipage}

\section{One Step of the Schur Algorithm, Expressed\\
in the Language of Colligations.\label{SSA}}
\setcounter{equation}{0}
The results from Section \ref{SASRoc} can be summarized as follows:
Starting from the unitary \(n\times{}n\)-matrix %
\begin{math}
 \begin{bmatrix}
 \alpha&\beta\,\\[0.5ex]
 \gamma&\delta\,
 \end{bmatrix}\,
\end{math}
representing a given inner rational matrix-function \(\omega(z)\) of degree \(n\),
\begin{equation*}
\omega(z)=\alpha+z\beta(I-z\delta)^{-1}\gamma\,,
\end{equation*}
we constructed the unitary \((n+1)\times(n+1)\)-matrix
\begin{math}
 \begin{bmatrix}
 s_0&B\,\\[0.5ex]
 C&D\,
 \end{bmatrix}\,
\end{math}
representing the function \(s(z)\):
\begin{equation*}
s(z)=s_0+zB(I-zD)^{-1}C\,,
\end{equation*}
where \(s(z)\) is the inverse Schur transform \eqref{sSA}.

Our goal is not, however, to determine  \(s(z)\) from \(\omega(z)\),
but instead to start with \(s(z)\) and determine \(\omega(z)\).
We look to describe a step of the Schur algorithm when
applied to a rational inner function \(s(z)\),
\begin{equation*}
s(z)\longrightarrow \omega(z)=\frac{1}{z}\frac{s(z)-s_0}{1-s(z)\overline{s_0}}\,,
\qquad s_0=s(0)\,,
\end{equation*}
in terms of system representations. In other words, we would like to find
the unitary matrix \(\begin{bmatrix}\alpha&\beta\\
\gamma&\delta\end{bmatrix}\) representing \(\omega(z)\), starting from
the matrix \(U\) representing the function \(s(z)\).

Equation \eqref{SCoo} serves as a heuristic argument. Until
now,
\(\begin{bmatrix}
 \alpha & \beta\\[0.5ex]
 \gamma & \delta
\end{bmatrix}\)
was given and \(U\) was unknown. Now we assume that
the unitary matrix \(U\) is given and that the matrix
 \(\begin{bmatrix}
 \alpha & \beta\\[0.5ex]
 \gamma & \delta
\end{bmatrix}\)
 is unknown. We consider \eqref{SCoo} as an equation with respect to
 the matrix
\(\begin{bmatrix}
 \alpha & \beta\\[0.5ex]
 \gamma & \delta
\end{bmatrix} \)
 and \(U\) as given. Because the second factor on the right-hand
 side of \eqref{SCoo} is a unitary matrix, the solution matrix
\(\begin{bmatrix}
 \alpha & \beta\\[0.5ex]
 \gamma & \delta
\end{bmatrix}\)
(if it exists) for equation \eqref{SCoo} is also a unitary matrix.

For a \textit{general}
unitary matrix \(U\), equation \eqref{SCoo} has \textit{no solution} \linebreak
with respect to the matrix
\(\begin{bmatrix}
 \ \alpha & \beta\\[0.5ex]
 \ \gamma & \delta
\end{bmatrix}\): We know that the block-matrix \linebreak entry \(B\) in
 \(U=\begin{bmatrix}s_0&B\\C&D\end{bmatrix} \ \)
 (\(U\) as in \eqref{SCoo})
 is necessarily of the form \linebreak
\(B=\begin{bmatrix}(1-|s_0|^2)^{1/2}  & \ 0_{1\times(n-1)}\end{bmatrix}\,.\)

 Since the characteristic functions of unitarily equivalent colligations coincide,
 it is enough to find a solution for \eqref{SCoo} with \(U\) replaced by some
 matrix \(U^V\) of the form \eqref{DeFr}:
\begin{gather}
U^V=
\begin{bmatrix}
1  & \ 0 \ & 0_{1\times(n-1)} \\[0.7ex]
0   & \ \alpha & \beta\\[0.5ex]
0_{(n-1)\times{}1}   & \ \gamma & \delta
\end{bmatrix}
\begin{bmatrix}
s_0 & \ (1-|s_0|^2)^{1/2} \ &\  0_{1\times{}(n-1)} \\[0.5ex]
(1-|s_0|^2)^{1/2}& \!\!-\overline{s_0} \ &  0_{1\times(n-1)}\\[0.8ex]
0_{(n-1)\times{}1}   &  0_{(n-1)\times{}1} &  1_{(n-1)\times(n-1)}
\end{bmatrix}
\!,\notag \\
\label{SCooM}%
\end{gather}

\begin{lemma}
\label{sSa} Given a unitary \((n+1)\times(n+1)\)-matrix \(U\),
the unitary \(n\times{}n\)-matrix \(V\) can be chosen such
that equation \eqref{SCooM} has a solution.
\end{lemma}
\begin{lemma}
\label{RMCF}
Given a unitary \((n+1)\times{}(n+1)\)-matrix \(U\):
\begin{equation}
\label{GivU}%
 U=%
 \begin{bmatrix}
\ s_0\ &B\\[0.5ex]
 C&D
 \end{bmatrix}\,,
\end{equation}
we can find a unitary \(n\times{}n\)-matrix \(V_0\) such that \(U^{V_0}\),
given by
\begin{equation*}
U^{\,V_0}=\begin{bmatrix}
 1 &0_{n\times{}n}\ \\[0.5ex]
0_{n\times{}n}    & V_0^\ast\end{bmatrix}\,%
U\,%
\begin{bmatrix}
 1 &0_{n\times{}n}\ \\[0.5ex]
0_{n\times{}n}   & V_0
\end{bmatrix}\,,
\end{equation*}
 takes the form \(U^{V_0}=U^0\), where
\begin{equation}
\label{GivUr}%
 U^0=%
 \begin{bmatrix}
\ s_0\ &B_0\\[0.5ex]
 C_0&D_0
 \end{bmatrix}\,,
\end{equation}
and the block-matrix entry \(B_0\in\mathfrak{M}_{1\times{}n}\)
is
\begin{equation}
\label{b0}%
 B_0=\big[
 (1-|s_0|^2)^{1/2} \ \ \cdots \ \ 0_{1\times(n-1)}\big]\,.
\end{equation}
\end{lemma}
\noindent
\textsf{PROOF}. The row-vectors \(B\) and \(B_0\) satisfy the condition
\[BB^{\ast}=B_0B_0^{\ast}\,\ \ \left( \, =(1-|s_0|^2) \, \right) \,\]
The equality \(B_0B_0^{\ast}=1-|s_0|^2\) follows from the definition of \(B_0\), \eqref{b0}.
The equality \(s_0\overline{s_0}+BB^{\ast}=1\) holds, since the matrix \(U\), \eqref{GivU}, is unitary.
Applying Lemma \ref{Hh} to the row-vectors \(B\) and \(B_0\), we find the unitary
\(n\times{}n\)-matrix \(V_0\) such that \(BV_0=B_0\). For every such choice of \(V_0\), the matrix
\(U^{V_0}\) has the form \eqref{GivUr}-\eqref{b0}.
\begin{remark}
\label{NUNM}
If \(n>1\), the matrices \(U^0\) and \(V_0\) are not uniquely defined.
The row-vector
\(B\) of any matrix \(U^V\) with \(V\) of the form
\begin{math}V=V_0
\bigl[\begin{smallmatrix}
1&0\\[0.4ex]
0&v
\end{smallmatrix}\bigr]\,,
\end{math}
where \(v\) is an arbitrary unitary \((n-1)\times(n-1)\)-matrix
is also of the form \eqref{b0}.
\end{remark}
\begin{theorem}
\label{fnt}%
 Given a unitary \((n+1)\times(n+1)\)-matrix of the form
\begin{equation}
\label{GivU-0}%
 U^0=%
 \begin{bmatrix}
\ s_0\ &(1-|s_0|^2)^{1/2} & 0_{1\times{}(n-1)}\\[0.7ex]
 c_1&d_{11}&d_{12}\\[0.7ex]
 c_2&d_{21}&d_{22}
 \end{bmatrix},
\end{equation}
where \(s_0\in\mathbb{C},\,|s_0|\leq{}1\), \(n\geq{}2\),
\begin{equation*}
\begin{matrix}
c_1\in\mathfrak{M}_{1\times{}1}\ \ \ \ \ ,&d_{11}\in\mathfrak{M}_{1\times{}1}\ \ \ \ , &
d_{12}\in\mathfrak{M}_{1\times{}(n-1)} \ \ \ \ ,\\
c_2\in\mathfrak{M}_{(n-1)\times{}1}\,,&d_{21}\in\mathfrak{M}_{(n-1)\times{}1}, &d_{22}\in\mathfrak{M}_{(n-1)\times{}(n-1)} \,,
\end{matrix}
\end{equation*}
the equation
\begin{equation}
\label{SE}%
 U^0=
\begin{bmatrix}
1  & \ 0 \ & 0_{1\times(n-1)} \\[0.7ex]
0   & \ \alpha & \beta\\[0.5ex]
0_{(n-1)\times{}1}   & \ \gamma & \delta
\end{bmatrix}
\begin{bmatrix}
s_0 & \ (1-|s_0|^2)^{1/2} \ &\  0_{1\times{}(n-1)} \\[0.5ex]
(1-|s_0|^2)^{1/2}& \!\!-\overline{s_0} \ &  0_{1\times(n-1)}\\[0.8ex]
0_{(n-1)\times{}1}   &  0_{(n-1)\times{}1} &  1_{(n-1)\times(n-1)}
\end{bmatrix}
\!,
\end{equation}
where
\begin{equation*}
\begin{matrix}
\alpha\in\mathfrak{M}_{1\times{}1}\ \ \ \ , &
\beta\in\mathfrak{M}_{1\times{}(n-1)} \ \ \ \ ,\\
\gamma\in\mathfrak{M}_{(n-1)\times{}1}, &\delta\in\mathfrak{M}_{(n-1)\times{}(n-1)} \,,
\end{matrix}
\end{equation*}
has a solution with respect to the matrix
\begin{equation}%
\label{U-1}
U^1= %
\begin{bmatrix}
 \alpha & \beta\\[0.5ex]
 \gamma & \delta
\end{bmatrix}\,.
\end{equation}
The solution of this equation can be expressed as
\begin{equation}
\label{EES}
\begin{bmatrix}
 \alpha & \beta\\[0.57ex]
 \gamma & \delta
\end{bmatrix}=
\begin{bmatrix}
-d_{11}\,s_0+c_1\,(1-|s_0|^2)^{1/2}  & \ \ d_{12}\\[0.7ex]
 -d_{21}\,s_0+c_2\,(1-|s_0|^2)^{1/2} & \ \ d_{22}
\end{bmatrix}
\end{equation}
\end{theorem}
\textsf{PROOF of THEOREM \ref{fnt}.}
We consider equation \eqref{SE} in further detail.
If this equation is solvable, then
\begin{multline}
\label{SEr}%
 \begin{bmatrix}
\ s_0\ &(1-|s_0|^2)^{1/2} & 0_{1\times{}(n-1)}\\[0.7ex]
 c_1&d_{11}&d_{12}\\[0.7ex]
 c_2&d_{21}&d_{22}
 \end{bmatrix}\times\\
\times\begin{bmatrix}
\overline{s_0} & \ (1-|s_0|^2)^{1/2} \ &\  0_{1\times{}(n-1)} \\[0.5ex]
(1-|s_0|^2)^{1/2}& \!\!-s_0 \ &  0_{1\times(n-1)}\\[0.8ex]
0_{(n-1)\times{}1}   &  0_{(n-1)\times{}1} &  1_{(n-1)\times(n-1)}
\end{bmatrix}
=\\
=\begin{bmatrix}
1  & \ 0 \ & 0_{1\times(n-1)} \\[0.7ex]
0   & \ \alpha & \beta\\[0.5ex]
0_{(n-1)\times{}1}   & \ \gamma & \delta
\end{bmatrix}
\end{multline}
Multiplying the matrices on the left-hand side
of \eqref{SEr}, we see that their product is of the form
\begin{math}
\Bigl[\begin{smallmatrix}
1  &  0  & 0 \\[0.3ex]
\ast   & \ast& \ast\\[0.3ex]
\ast   & \ast & \ast
\end{smallmatrix}\Bigr]\,.
\end{math}
 Since the matrix \(U^0\), \eqref{GivU-0},
is unitary, the scalar product of its different rows vanishes.
 The fact that the first row of this matrix is orthogonal to each
 other row can be expressed as
 \begin{equation}%
 \label{ld}
 \begin{bmatrix}
 c_1\\
 c_2
 \end{bmatrix}%
 \overline{s_0}+%
 \begin{bmatrix}
 d_{11}\\
 d_{21}
 \end{bmatrix}%
 (1-|s_0|^2)^{1/2}=0\,.
 \end{equation}
 The latter equalities mean that the product of the
 matrices on the left-hand side of \eqref{SEr} takes the form
\begin{math}
\Bigl[
\begin{smallmatrix}
\ast  &  \ast  & \ast \\[0.3ex]
0   & \ast & \ast \\[0.3ex]
0   & \ast & \ast
\end{smallmatrix}\Bigr]\,.
\end{math}
Thus, the product of the matrices on the left-hand side of
\eqref{SEr} has the desired form
\begin{math}
\Bigl[
\begin{smallmatrix}
1  &  0  & 0 \\[0.3ex]
0   & \ast & \ast \\[0.3ex]
0   & \ast & \ast
\end{smallmatrix}\Bigr]\,.
\end{math}
Multiplying out the matrices in \eqref{SEr}, we obtain \eqref{EES}.
\hfill{Q.E.D.}
\begin{remark}
In view of \eqref{ld}, the solution
\begin{math}
\Bigl[
 \begin{smallmatrix}
  \ \alpha & \beta\\[0.5ex]
 \ \gamma & \delta
 \end{smallmatrix}
 \Bigr]
 \end{math}
 of equation \eqref{SE} can also be written as:
 \begin{equation}
 \label{yaf1}
 \begin{bmatrix}
 \ \alpha & \beta\\[0.5ex]
 \ \gamma & \delta
 \end{bmatrix}
 =%
 \begin{bmatrix}
 (1-|s_0|^2)^{-1/2} c_1 & d_{12}\\[0.5ex]
  (1-|s_0|^2)^{-1/2} c_2 & d_{22}
 \end{bmatrix}\,\quad \text{if} \ \ |s_0|<1,
\end{equation}
and
 \begin{equation}
 \label{yaf2}
 \begin{bmatrix}
 \ \alpha & \beta\\[0.5ex]
 \ \gamma & \delta
 \end{bmatrix}
 =%
 \begin{bmatrix}
 -(\overline{s_0})^{-1}d_{11} & d_{12}\\[0.5ex]
 -(\overline{s_0})^{-1}d_{21} & d_{22}
 \end{bmatrix}\, \quad \text{if} \ \  s_0\not=0\,.
 \end{equation}
\end{remark}
\vspace{2.0ex}
\noindent
\begin{remark}
\label{LShP}
 If \(n=1\), then there is no room for the matrices
\(d_{12},d_{21},d_{22}\) and \(\beta,\gamma,\delta\). In this
case \(U^0\), \eqref{GivU-0}, should be replaced by the matrix:
\begin{equation}
\label{GivU-0_r}%
 U^0=%
 \begin{bmatrix}
\ s_0\ &(1-|s_0|^2)^{1/2}\\[0.7ex]
 c_1&d_{11}
 \end{bmatrix}\,,
\end{equation}
where \(s_0\in\mathbb{C}\) with \(|s_0|\leq{}1\),
\begin{equation*}
\begin{matrix}
c_1\in\mathfrak{M}_{1\times{}1}\ \ \ \ \
,&d_{11}\in\mathfrak{M}_{1\times{}1}
\end{matrix}\,,
\end{equation*}
and the matrix \(U^1\), \eqref{U-1}, should be replaced with:
matrix \(U^1\)
\begin{equation}%
\label{U-1-r}
U^1= %
\begin{bmatrix}
 \alpha
\end{bmatrix}\,,
\end{equation}%
where
\begin{equation*}
\begin{matrix}
\alpha\in\mathfrak{M}_{1\times{}1}
\end{matrix}\,.
\end{equation*}
Equation \eqref{SE} takes the form
\begin{equation}
\label{SE-r}%
 U^0=
\begin{bmatrix}
1  & \ 0 \  \\[0.7ex]
0   & \ \alpha
\end{bmatrix}
\begin{bmatrix}
s_0 & \ (1-|s_0|^2)^{1/2} \  \\[0.5ex]
(1-|s_0|^2)^{1/2}& \!\!-\overline{s_0}
\end{bmatrix}
\!.
\end{equation}
The solution of this equation can be expressed as
\begin{equation}
\label{EES-r}
\begin{bmatrix}
 \alpha
\end{bmatrix}=
\begin{bmatrix}
-d_{11}\,s_0+c_1\,(1-|s_0|^2)^{1/2}
\end{bmatrix}\,,
\end{equation}
as well as in the forms:
 \begin{equation}
 \label{yaf1-r}
 \begin{bmatrix}
 \alpha
 \end{bmatrix}
 =%
 \begin{bmatrix}
 (1-|s_0|^2)^{-1/2} c_1
 \end{bmatrix}\,\quad \text{if} \ \ |s_0|<1,
\end{equation}
and
 \begin{equation}
 \label{yaf2-r}
 \begin{bmatrix}
 \alpha
 \end{bmatrix}
 =%
 \begin{bmatrix}
 -(\overline{s_0})^{-1}d_{11}
 \end{bmatrix}\, \quad \text{if} \ \  s_0\not=0\,.
 \end{equation}
\end{remark}%

\vspace{1.0ex}
Since both factors on the right-hand side of \eqref{SE} are unitary matrices,
we have that \(U^1\) is also a unitary matrix.
The matrix \(U^0\) in \eqref{GivU-0} can be considered as a matrix of
the unitary colligation \((\mathcal{E}^0,\mathcal{H}^0,U^0)\) with outer space
\(\mathcal{E}^0=\mathbb{C}\) and with inner
space \(\mathcal{H}^0=\mathbb{C}^n\). The matrix \(U^1\) in \eqref{U-1}
can, in turn, be considered as a matrix of the
the unitary colligation \((\mathcal{E}^1,\mathcal{H}^1,U^1)\)
with outer space \(\mathcal{E}^1=\mathbb{C}\) and with inner
space \(\mathcal{H}^1=\mathbb{C}^{n-1}\).
\begin{lemma} {\ } \\
\label{CO} \hspace*{4.0ex} \textup{I.} If the  colligation
\((\mathcal{E}^0,\mathcal{H}^0,U^0)\) is controllable, then the
colligation \\ %
\((\mathcal{E}^1,\mathcal{H}^1,U^1)\)
is also controllable.\\
\hspace*{4.0ex} \textup{II.} If the  colligation
\((\mathcal{E}^0,\mathcal{H}^0,U^0)\) is observable, \,then the
colligation \\ %
\((\mathcal{E}^1,\mathcal{H}^1,U^1)\) is also observable.
\end{lemma}
\textsf{PROOF.} Without loss of generality, we assume that \(|s_0|<1\).
Otherwise the colligation \((\mathcal{E}^0,\mathcal{H}^0,U^0)\) can not be neither controllable nor
observable. Our reasoning is based on the equalities
\begin{equation}
\label{UE1}
\begin{bmatrix}
 d_{11} & d_{12}\\[0.5ex]
 d_{21} & d_{22}
\end{bmatrix}
 =%
\begin{bmatrix}
  (-\overline{s_0})\alpha & \beta\\[0.5ex]
  (-\overline{s_0})\gamma & \delta
 \end{bmatrix}
 \, .
\end{equation}
and
\begin{equation}
\label{UE2}
\begin{bmatrix}
c_1\\[0.5ex]c_2
\end{bmatrix}
=%
(1-|s_0|^2)^{1/2}%
\begin{bmatrix}
\alpha\\[0.5ex]\gamma
\end{bmatrix}
\,.
\end{equation}
\textsf{Proof of the Statement I.} The condition that
the colligation \((\mathcal{E}^0,\mathcal{H}^0,U^0)\) be controllable means that
\begin{equation}
\label{CoCo}
\Vvee\limits_{0\leq{}k} \ %
\begin{bmatrix}
 d_{11} & d_{12}\\[0.5ex]
 d_{21} & d_{22}
\end{bmatrix}^k
\begin{bmatrix}
 c_{1}\\[0.5ex]
 c_{2}
\end{bmatrix}=
\mathfrak{M}_{n\times{}1}\,\quad\left(\ \ =
\begin{bmatrix}
\mathfrak{M}_{1\times{}1}\\[0.5ex]
\mathfrak{M}_{(n-1)\times{}1}
\end{bmatrix}\ \
\right)\,.
\end{equation}
And the controllability of the colligation
\((\mathcal{E}^1,\mathcal{H}^1,U^1)\) can be expressed as:
\begin{equation}
\label{coco}
\Vee\limits_{0\leq{}k}\,\delta^k\gamma=\mathfrak{M}_{(n-1)\times{}1}\,.
\end{equation}
We look to show that \eqref{coco} follows from \eqref{CoCo}.
In view of \eqref{UE1} and \eqref{UE2}, we can express \eqref{CoCo} as:
\begin{equation}
\label{CoCoR}
\Vvee\limits_{0\leq{}k} \ %
\begin{bmatrix}
(-\overline{s_0}) \alpha & \beta\\[0.5ex]
 (-\overline{s_0}) \beta & \delta
\end{bmatrix}^k
\begin{bmatrix}
 \alpha\\[0.5ex]
 \gamma
\end{bmatrix}=
\mathfrak{M}_{n\times{}1}\,\quad\left(\ \ =
\begin{bmatrix}
\mathfrak{M}_{1\times{}1}\\[0.5ex]
\mathfrak{M}_{(n-1)\times{}1}
\end{bmatrix}\ \
\right)\,.
\end{equation}
Let
\begin{equation}
\label{No}
\begin{bmatrix}
(-\overline{s_0}) \alpha & \beta\\[0.5ex]
 (-\overline{s_0}) \gamma & \delta
\end{bmatrix}^k
\begin{bmatrix}
 \alpha\\[0.5ex]
 \gamma
\end{bmatrix}=\begin{bmatrix}
 \ast\\[0.5ex]
 f_k
\end{bmatrix}\,,\quad k=0,\,1,\,2\,\,\ldots\,,
\end{equation}
where %
\begin{math}
f_k\in\mathfrak{M}_{(n-1)\times{}1}\,.
\end{math}
In view of \eqref{CoCoR},
\begin{equation}
\label{CoCor}
\Vee\limits_{0\leq{}k}
\begin{bmatrix}
f_k
\end{bmatrix}=\mathfrak{M}_{(n-1)\times{}1}\,.
\end{equation}
Clearly, we have that for every \(k=0,\,1,\,2,\,\ldots\)
\begin{equation}
\label{tm}
f_k=\xi_{0,k}\delta^0\gamma+\,\cdots\,+\xi_{k-1,k}\delta^{k-1}\gamma+
\delta^k\gamma\,,
\end{equation}
where \(\xi_{j,k},\,\,0\leq{}j\leq{}k-1\,,\) are some complex numbers.
Therefore,
\begin{equation*}
\Vee\limits_{0\leq{}k}
\begin{bmatrix}
f_k
\end{bmatrix}=\Vee\limits_{0\leq{}k}\,\delta^k\gamma\,.
\end{equation*}
We have thus proved Statement I.

\textsf{Proof of the Statement II.} The condition that the
colligation be observable \((\mathcal{E}^0,\mathcal{H}^0,U^0)\) can be written as:
\begin{equation}
\label{CoCo2}
\Vvee\limits_{0\leq{}k} \ %
\begin{bmatrix}
 (1-|s_0|^2)^{1/2}&0_{1\times(n-1)}
\end{bmatrix}
\begin{bmatrix}
 d_{11} & d_{12}\\[0.5ex]
 d_{21} & d_{22}
\end{bmatrix}^k
= \mathfrak{M}_{1\times{}n}\,\ \left(  =
\begin{bmatrix}
\mathfrak{M}_{1\times{}1}& \mathfrak{M}_{1\times{}(n-1)}
\end{bmatrix}
\right)\,.
\end{equation}
And the observability of the colligation
\((\mathcal{E}^1,\mathcal{H}^1,U^1)\) can be expressed as:
\begin{equation}
\label{coco2}
\Vee\limits_{0\leq{}k}\,\beta\delta^k=\mathfrak{M}_{1\times{}(n-1)}\,.
\end{equation}
We aim to show that \eqref{coco2} follows from the formulas \eqref{CoCo2},
\eqref{EES} and \eqref{ld}. In view of \eqref{UE1}, we can
express \eqref{CoCo2} as follows:
\begin{equation}
\label{CoCoR2}
\Vvee\limits_{0\leq{}k} \ %
\begin{bmatrix}
 1&
 0_{1\times(n-1)}
\end{bmatrix}
\begin{bmatrix}
(-\overline{s_0}) \alpha & \beta\\[0.5ex]
 (-\overline{s_0}) \beta & \delta
\end{bmatrix}^k
=
\mathfrak{M}_{1\times{}n}\,\quad\left(\ \ =
\begin{bmatrix}
\mathfrak{M}_{1\times{}1}&
\mathfrak{M}_{1\times{}(n-1)}
\end{bmatrix}\ \
\right)\,.
\end{equation}
\begin{math}
\end{math}
Let
\begin{equation}
\label{No2}
\begin{bmatrix}
 1& 0_{1\times{}(n-1)}
\end{bmatrix}
\begin{bmatrix}
(-\overline{s_0}) \alpha & \beta\\[0.5ex]
 (-\overline{s_0}) \gamma & \delta
\end{bmatrix}^k
=\begin{bmatrix}
 \ast&
 g_{k}
\end{bmatrix}\,, \quad k=0,\,1,\,2\,\,\ldots\,,
\end{equation}
where %
\begin{math}
g_k\in\mathfrak{M}_{1\times{}(n-1)}\,.
\end{math}
In view of \eqref{CoCoR2}, we have that
\begin{equation}
\label{CoCor2}
\Vee\limits_{0\leq{}k}
\begin{bmatrix}
g_k
\end{bmatrix}=\mathfrak{M}_{1\times{}(n-1)}\,.
\end{equation}
Clearly, \(g_0=0_{1\times(n-1)}\) and for every \(k=0,\,1,\,2,\,\ldots\)
\begin{equation}
\label{tm2}
g_{k+1}=\eta_{0,k}\,\beta\delta^0+\,\cdots\,+\eta_{k-1,k}\,\beta\delta^{k-1}+
\beta\delta^{k}\,,
\end{equation}
where \(\eta_{j,k},\,\,0\leq{}j\leq{}k-2\,,\) are some complex numbers.
Therefore,
\begin{equation*}
\Vee\limits_{0\leq{}k}
\begin{bmatrix}
g_k
\end{bmatrix}=\Vee\limits_{0\leq{}k}\,\beta\delta^k\,.
\end{equation*}
We have thus proved Statement II.\hfill Q.E.D.

\vspace{2.0ex} \noindent
 The following Lemma is an immediate consequence of Lemma \ref{CO} \\[-4.0ex]
\begin{lemma}
\label{COF} Let \(s_0\in\mathbb{C}\) with \(|s_0|<1\) and \(U^0\) be a
unitary \((n+1)\times(n+1)\)-matrix of the form \eqref{GivU-0}.
Suppose that the \(n\times{}n\)-matrix  \(U^1\), \eqref{U-1}, is
related to the matrix \(U^0\) by equation \eqref{SE}. Let
\((\mathcal{E}^0,\mathcal{H}^0,U^0)\) and
\((\mathcal{E}^1,\mathcal{H}^1,U^1)\) be the above-described
operator colligations related to the matrices \(U^0\)
and \(U^1\). If the colligation \((\mathcal{E}^0,\mathcal{H}^0,U^0)\)
is minimal, then the colligation \((\mathcal{E}^1,\mathcal{H}^1,U^1)\)
is also minimal.
\end{lemma}
\begin{theorem}
\label{StScAl}
 Let \(s(z)\) be a  rational inner
matrix-function of degree \(n>1\) (\(s(z)\) is thus non-constant
and \(|s_0|<1\), where \(s_0=s(0)\)) and let
\begin{equation}
\label{om}
\omega(z)=\frac{1}{z}\cdot\frac{s(z)-s_0}{1-s(z)\overline{s_0}}\,,\qquad s_0=s(0),
\end{equation}
be the Schur transformation of the function \(s(z)\).\\
Let the unitary matrix \(U\),%
\begin{equation}%
\label{OSSi}
U=
\begin{bmatrix}
s_0&B_0\\[0.5ex]
C_0&D_0
\end{bmatrix}\,,
\quad B_0\in\mathfrak{M}_{1\times{}n},\,C_0\in\mathfrak{M}_{n\times{}1},\,
D_0\in\mathfrak{M}_{n\times{}n},
\end{equation}%
which yields the minimal system representation
\begin{equation}
\label{SR-s}%
 s(z)=s_0+zB_0(I-zD_0)^{-1}C_0\,,
\end{equation}
have row \(B_0\) of the special form
\begin{equation}
\label{SF}
B_0=[b\ \ \ 0_{1\times{(n-1)}}]\,, \qquad b>0\,.
\end{equation}

Then the function \(\omega(z)\) admits the system representation
\begin{equation}
\label{SR-o}%
 \omega(z)=\alpha+z\beta(I-z\delta)^{-1}\gamma\,,
\end{equation}
where the unitary \(n\times{}n\)-matrix
\begin{math}%
\begin{bmatrix}
 \ \alpha & \beta\\[0.5ex]
 \ \gamma & \delta
\end{bmatrix}\,,
\end{math}
\[
\alpha\in\mathfrak{M}_{1\times{}1}, \ \
\beta\in\mathfrak{M}_{1\times{}(n-1)}, \ \
\gamma\in\mathfrak{M}_{(n-1)\times{}1}, \ \ \delta\in\mathfrak{M}_{(n-1)\times{}(n-1)} \,,
\]
can be determined from the matrix \(U^0\) using
 \begin{equation}
 \label{Rrc}
 \begin{bmatrix}
 \ \alpha & \beta\\[0.5ex]
 \ \gamma & \delta
 \end{bmatrix}
 =%
 \begin{bmatrix}
 (1-|s_0|^2)^{-1/2} c_1 & d_{12}\\[0.5ex]
  (1-|s_0|^2)^{-1/2} c_2 & d_{22}
 \end{bmatrix}\,,
\end{equation}
where \(c_j\) and \(d_{jk}\) are the block-matrix entries of
the block-matrix decompositions
\begin{equation}%
\label{FBD}
C_0=
\begin{bmatrix}
c_1\\[0.5ex]
c_2
\end{bmatrix}\,
\qquad
D_0=
\begin{bmatrix}
d_{11}&d_{12}\\[0.5ex]
d_{21}&d_{22}
\end{bmatrix}\,,
\end{equation}
\[c_1\in\mathfrak{M}_{1\times{}1},\ c_2\in\mathfrak{M}_{(n-1)\times{}1},\]
\[ D_{11}\in\mathfrak{M}_{1\times{}1},\ D_{12}\in\mathfrak{M}_{1\times{}(n-1)},\
D_{21}\in\mathfrak{M}_{(n-1)\times{}1},\
D_{22}\in\mathfrak{M}_{(n-1)\times{}(n-1)}\,.\] The unitary
colligation associated with the matrix \(\tilde{U}^1\) is minimal.
\end{theorem}
\textsf{PROOF}.  The matrix \(U^0\), \eqref{OSSi}, the matrix
\begin{math}
 \begin{bmatrix}
 \ \alpha & \beta\\[0.5ex]
 \ \gamma & \delta
 \end{bmatrix}
\end{math}, \eqref{Rrc}, and the number \(s_0\) are related by
equation \eqref{SCoo}. According to Theorem \ref{SRIST}, the
function \(\omega(z)\), defined by \eqref{SR-o}, and the function
\(s(z)\) are related by the equality \eqref{sSA}. \hfill
Q.E.D.\\[1.0ex]

Theorem \ref{StScAl} together with Lemma \ref{RMCF} describe a step
of the Schur algorithm in terms of system representations.
Before applying the direct Schur transform \eqref{om}, which is a
step of the Schur algorithm, we should first `normalize' the colligation
matrix \(U\) representing the `initial' function \(s(z)\).
This normalization starts with the matrix \(U\), from which we determine
the unitarily equivalent matrix \(U^0\), \eqref{GivU-0},
whose row \(B^0\) is of the special form \eqref{GivUr}. We then aim
to solve the equation \eqref{SE} with respect to the
matrix
\begin{math}
\bigl[\begin{smallmatrix}
 \ \alpha & \beta\\[0.5ex]
 \ \gamma & \delta
\end{smallmatrix}\bigr]\,.
\end{math}
The solution of \eqref{SE} is given by \eqref{Rrc}.
The unitary matrix %
\begin{math}
\bigl[\begin{smallmatrix}
 \ \alpha & \beta\\[0.5ex]
 \ \gamma & \delta
\end{smallmatrix}\bigr]
\end{math} %
yields the system representation of the function \(\omega(z)\).

It should be emphasized that, in general, \textit{the matrix %
\begin{math}
\bigl[\begin{smallmatrix}
 \ \alpha & \beta\\[0.5ex]
 \ \gamma & \delta
\end{smallmatrix}\bigr]
\end{math} %
is not normalized, i.e. its row \(\beta\) is not of the form
\(\beta=[\ast\ \ 0_{1\times(n-2)}]\).}
To perform the next step of the Schur algorithm, we must therefore `normalize'
the matrix \begin{math}
\bigl[\begin{smallmatrix}
 \ \alpha & \beta\\[0.5ex]
 \ \gamma & \delta
\end{smallmatrix}\bigr]
\end{math}, %
obtaining the `normalized' form \(U^1\).
We then have to solve the equation of the form \eqref{SE}, where \(U^0\) is
replaced by \(U^1\), etc. The normalization procedure must therefore be performed
\emph{at every step of the Schur algorithm.} This normalization
procedure is, however, not quite unique. It has some degrees of freedom
(See Remark \ref{NUNM}.)
It turns out that we can use these degrees of freedom to make the normalization procedure a
\emph{one-time} procedure, so that it might be dealt with during \emph{preprocessing}
for the further step-by-step recurrence. In further processing there is then no need
for normalization and one only has to solve the recurrent chain of equations of the form \eqref{SE}.
A \emph{one-time normalization} of this kind is related to the reduction
of the `initial' colligation matrix to the \emph{lower Hessenberg form}.

\section{Hessenberg Matrices. \\ The Householder Algorithm.
\label{HaHe}} \setcounter{equation}{0}%

 Roughly speaking, the lower (upper) Hessenberg matrix,
 is a matrix which is almost lower (upper) triangular. The precise definition is:
\begin{definition}\ \\
\label{HeML}%
\hspace*{2.0ex}\textup{\textsf{I.}} We say that a square matrix \(H\) is
a \textsf{lower Hessenberg matrix} if  it has zero-entries
above the first superdiagonal. If
\(H=||h_{jk}||_{0\leq{}j,k\leq{}n}\), then \(H\) is lower
Hessenberg matrix if \(h_{jk}=0\) for \(k>j+1,\,0\leq{}j\leq{}n-1\).\\
\hspace*{2.0ex}\textup{\textsf{II.}} We say that a lower Hessenberg matrix
\(H=||h_{jk}||_{0\leq{}j,k\leq{}n}\) is \textsf{special}
if all entries of its first superdiagonal are non-negative:
\(h_{j,j+1}\geq{}0\,,0\leq{}j\leq{}n-1.\)\\
\hspace*{2.0ex}\textup{\textsf{III.}} We say that a Hessenberg matrix
\(H=||h_{jk}||_{0\leq{}j,k\leq{}n}\) is
\textsf{HL-non-singular} if all entries of its first superdiagonal
are non-zero: \(h_{j,j+1}\not=0\,,0\leq{}j\leq{}n-1.\)
\end{definition}
The definition of an upper Hessenberg matrix, special
upper Hessenberg matrix and non-singular upper Hessenberg matrix
is similar to Definition \ref{HeML}:
\begin{definition}\ \\
\label{HeMLU}%
\hspace*{2.0ex}\textup{\textsf{I.}} We say that a square matrix \(H\) is
an \textsf{upper Hessenberg matrix} if it has zero-entries
below the first subdiagonal. If
\(H=||h_{jk}||_{\,0\leq{}j,k\leq{}n}\), then \(H\) is upper
Hessenberg matrix if \(h_{jk}=0\) for \(k<j-1,\,0\leq{}j\leq{}n-1\).\\
\hspace*{2.0ex}\textup{\textsf{II.}} We say that an upper Hessenberg matrix
\(H=||h_{jk}||_{0\leq{}j,k\leq{}n}\) is \textsf{special}
if all entries of its first subdiagonal are non-negative:
\(h_{j,j-1}\geq{}0\,,1\leq{}j\leq{}n.\)\\
\hspace*{2.0ex}\textup{\textsf{III.}} We say that an upper Hessenberg matrix
\(H=||h_{jk}||_{0\leq{}j,k\leq{}n-1}\) is
\textsf{HU-non-singular} if all entries of its first subdiagonal are
non-zero: \(h_{j,j-1}\not=0\,,1\leq{}j\leq{}n.\)
\end{definition}

\vspace{2.0ex}%
Hessenberg matrices were investigated by Karl Hessenberg
(1904-1959), a German engineer whose dissertation dealt with the
computation of eigenvalues and eigenvectors of linear operators.

\begin{theorem}\ \\
\label{RUHF}
 \hspace*{2.0ex}\textup{\textsf{I}}. Given an
\((n+1)\times{}(n+1)\)-matrix
 \(M=||M_{j,k}||_{0\leq{}j,k\leq{}n}\),
there exists a unitary \(n\times{}n\)-matrix \(V\) such that
 the matrix \(H^L\),
 \begin{equation}
 \label{RLH}
 H^L=
 \begin{bmatrix}
 1&0_{1\times{}n}\\
 0_{n\times{}1}&V^{\ast}
 \end{bmatrix}
 M%
\begin{bmatrix}
 1&0_{1\times{}n}\\
 0_{n\times{}1}&V
 \end{bmatrix}
 \end{equation}
  is a special lower Hessenberg matrix.\\[0.5ex]
\hspace*{2.0ex}\textup{\textsf{II}}. If the matrix \(M\) is
HL-non-singular, then both matrices \(H^L\) and \(V\) are
uniquely determined. From the equalities
\begin{equation}
 \label{RLHUn}
 H_j^L=
 \begin{bmatrix}
 1&0_{1\times{}n}\\
 0_{n\times{}1}&V^{\ast}_j
 \end{bmatrix}
 M%
\begin{bmatrix}
 1&0_{1\times{}n}\\
 0_{n\times{}1}&V_j
 \end{bmatrix}\,,\quad j=1,\,2\,,
 \end{equation}
 where \(H_1^L\) and \(H_2^L\) are special upper Hessenberg matrices,
 \(V_1\) and \(V_2\) are unitary matrices
 and the Hessenberg matrix \(H_1^L\) is HL-non-singular, it follows
 that \(H_2^L=H_1^L\) and \(V_2=V_1\).
\end{theorem}
\begin{definition}
\label{UHF0}%
 Given a square matrix \(M\), a lower Hessenberg
matrix \(H^L\) to which \(M\) can be reduced, \eqref{RLH}, is
called a lower Hessenberg form of the matrix \(M\).
\end{definition}
\begin{theorem}\ \\
\label{RLHF}
 \hspace*{2.0ex}\textup{\textsf{I}}. Given an
\((n+1)\times{}(n+1)\)-matrix
 \(M=||M_{j,k}||_{\,0\leq{}j,k\leq{}n}\),
there exists a unitary \(n\times{}n\)-matrix \(V\) such that
 the matrix \(H^U\),
 \begin{equation}
 \label{RUH}
 H^U=
 \begin{bmatrix}
 1&0_{1\times{}n}\\
 0_{n\times{}1}&V^{\ast}
 \end{bmatrix}
 M%
\begin{bmatrix}
 1&0_{1\times{}n}\\
 0_{n\times{}1}&V
 \end{bmatrix}
 \end{equation}
  is a special upper Hessenberg matrix.\\[0.5ex]
\hspace*{2.0ex}\textup{\textsf{II}}. If the matrix \(M\) is
HU-non-singular, then both matrices \(H^U\) and \(V\) are
uniquely determined. From the equalities
\begin{equation}
 \label{RUHUn}
 H_j^U=
 \begin{bmatrix}
 1&0_{1\times{}n}\\
 0_{n\times{}1}&V^{\ast}
 \end{bmatrix}
 M%
\begin{bmatrix}
 1&0_{1\times{}n}\\
 0_{n\times{}1}&V_j
 \end{bmatrix}\,,\quad j=1,\,2\,,
 \end{equation}
 where \(H_1^U\) and \(H_2^U\) are upper Hessenberg matrices,
 \(V_1\) and \(V_2\) are unitary matrices and the Hessenberg matrix
 \(H_1^U\) is HU-non-singular, it follows that \(H_2^U=H_1^U\) and \(V_2=V_1\).
\end{theorem}
\begin{definition}
\label{LHF}%
 Given a square matrix \(M\), an upper Hessenberg
matrix \(H^U\) to which \(M\) can be reduced, \eqref{RUH}, is
called an upper Hessenberg form of the matrix \(M\).
\end{definition}
\begin{theorem}
\label{HCO}
Let \(U\) be an \((n+1)\times{}(n+1)\)-unitary matrix.\\
\hspace*{2.0ex}\textup{\textsf{I}}. The unitary colligation
associated with the matrix \(U\) is observable if and only if
lower Hessenberg form of \(U\) is HL-non-singular.\\[0.5ex]
\hspace*{2.0ex}\textup{\textsf{II}}. The unitary colligation
associated with the matrix \(U\) is controllable if and only if the
upper Hessenberg form of \(U\) is HU-non-singular.
\end{theorem}

\begin{corollary}\label{LAE}
\textup{According to Theorem \ref{SimMin}, the finite-dimensional
unitary colligation is observable if and only if it is
controllable. Thus,} for a unitary matrix \(U\), the lower
Hessenberg form of \(U\) is \(HL\)-nonsingular if and only if the
upper Hessenberg form of \(U\) is \(HU\)-nonsingular.
\end{corollary}

\noindent%
\begin{lemma}
\label{Hh} \ \ \ Given two row-vectors
\(B^{\,\prime}=[b_1^{\,\prime} \, \, b_2^{\,\prime} \ . \ . \ . \
.\,\, b_n^{\,\prime}]\in \mathfrak{M}_{1 \times n} \)  and
\(B^{\,\prime\prime}=[b_1^{\,\prime\prime} \, \,
b_2^{\,\prime\prime} \ . \ . \ . \ .\,\,
b_n^{\,\prime\prime}]\in \mathfrak{M}_{1 \times n} \) having same norm:
\begin{equation}%
\label{EqN}%
B^{\,\prime}{B^{\,\prime}}^{\,\ast}=B^{\,\prime\prime}{B^{\,\prime\prime}}^{\,\ast}\,,
\end{equation}
there exists a unitary \(n\times{}n\)-matrix \(V\) such that
\begin{equation}%
\label{OtoT}%
B^{\,\prime}V=B^{\,\prime\prime}.
\end{equation}%
\end{lemma}
\noindent \textsf{PROOF of LEMMA \ref{Hh}.} We first consider the
question in a more general setting.
 Assume that \(\mathfrak{H}\) is a \textit{complex} Hilbert
 space with scalar product \(\langle\,u\,,\,v\rangle\),
 where \(\langle\,u\,,\,v\rangle\) is linear with respect to the
 argument \(u\) and antilinear with respect to \(v\).
 Let \(x\in\mathfrak{H}\) and \(y\in\mathfrak{H}\) be two
 vectors such that \(\langle\,x\,,\,x\rangle=\langle\,y\,,\,y\rangle\not=0.\)
 Let \(||u||\) denote the norm of the vector \(u\):
 \(||u||={\langle\,u\,,\,u\rangle}^{1/2}\). Given two vectors
 \(x\in\mathfrak{H}\),  \(y\in\mathfrak{H}\) such that
 \(||x||=||y||\not={}0\),
 our goal is to construct a unitary operator
 \(V:\mathfrak{H}\to\mathfrak{H}\) such that \(Vx=y\). If the
 vector \(y\) is proportional to the vector \(x\):
 \(x=\lambda{}y\) for some \(\lambda\in\mathbb{C}\), we put
 \(Vz\stackrel{}{}=\lambda{}z\ \ \forall z\in\mathbb{C}\). This
 operator is unitary: \(|\lambda|=1\), because
 \(||x||=||y||\not={}0\). If the vectors \(x\) and \(y\) are not
 proportional, we choose \(\lambda\in\mathbb{C},\,|\lambda|=1\) such that
  \(\lambda\langle\,x\,,\,y\rangle\geq{}0\). (If \(\langle\,x\,,\,y\rangle\not={}0\),
  then this \(\lambda\) is unique. If \(\langle\,x\,,\,y\rangle={}0\), we can choose
  arbitrary \(\lambda\) with \(|\lambda|=1\).) Let
  \begin{equation}
  \label{Chl}
  Vz=\lambda{}z-2\frac{\langle{}z,x-
  \overline{\lambda}y\rangle}{||x-\overline{\lambda}y
  ||^2}(\lambda{}x-y)\,\ \quad \forall \ z\in\mathfrak{H}\,.
  \end{equation}
  The vectors \[e_1=x+\overline{\lambda}y \text{ \ \ and \ \ }
  e_2=x-\overline{\lambda}y\]
  are non-zero (\(x\) and \(y\) are not proportional to one another) and
  orthogonal:
  \begin{equation}
  \label{Oe}
\langle{}e_1,\,e_2\rangle=0\,,
  \end{equation}
  because
  \[\langle{}x+\overline{\lambda}y\,,\,x-\overline{\lambda}y\rangle{}=
 \langle{}x,\,x\rangle-\overline{\lambda}\lambda\langle{}y,\,y\rangle+
 \overline{\lambda}\langle{}y,\,x\rangle -\lambda\langle{}x,\,y\rangle\]
 and \(\langle{}x,\,x\rangle=\langle{}y,\,y\rangle\),
 \(\overline{\lambda}\lambda=1\), \
  \(\lambda\langle{}x,\,y\rangle=\overline{\lambda\langle{}x,\,y\rangle}=
  \overline{\lambda}\langle{}y,\,x\rangle\). From \eqref{Chl} and \eqref{Oe} it
  follows that
 \begin{equation}
  \label{Ev1}
Ve_1=\lambda{}e_1\,.
  \end{equation}
  From \eqref{Chl} it follows that
\begin{equation}
  \label{Ev2}
Ve_2=-\lambda{}e_2\,,
  \end{equation}
  and
\begin{equation}
  \label{Ev3}
Vz=\lambda{}z\,
  \end{equation}
   \(\forall\,z\in\mathfrak{H}: \langle{}z,\,e_1\rangle=0, \
  \langle{}z,\,e_2\rangle=0\). Therefore the operator \(V\) is
  unitary. Since
  \[x=\frac{1}{2}(e_1+e_2), \quad y=\frac{\lambda}{2}(e_1-e_2)\,\]
  from \eqref{Ev1} and \eqref{Ev2} it follows that \(Vx=y\).

  Let us turn to the proof of the statement of Lemma \ref{Hh}.
  Let \(\mathfrak{H}\) be the set of all \(n\)-row-vectors with
  complex entries (in other words,
  \(\mathfrak{H}=\mathfrak{M}_{1\times{}n}\))
  and with the following scalar product:
  if \(u=[u_1,\,\ldots\,,\,u_n]\) and \(v=[v_1,\,\ldots\,,\,v_n]\)
  are vectors in  \(\mathfrak{H}\), then their scalar product
  \(\langle{}u\,,\,v\rangle\) is defined as
  \[\langle{}u\,,\,v\rangle=u\,v^{\ast}\,\]
  where \(v^{\ast}\) is the Hermitian conjugate of the row-vector \(v\).
  If \(H\) is some \(n\times{}n\)-matrix, then it generates an operator
  in \(\mathfrak{H}\). This operator maps the row-vector \(u\) to the
  row-vector \(uH\), where \(uH\) is the product of the
  \textit{matrices} \(u\) and \(H\). This operator is unitary if and only if
  \(H\) is unitary.

  In the notation of Lemma \ref{Hh}: \(x=B^{\prime}=[b_1^{\,\prime} \, \,
b_2^{\,\prime} \ . \ . \ . \ .\,\, b_n^{\,\prime}]\),
\(y=B^{\prime\prime}=[b_1^{\,\prime\prime} \, \,
b_2^{\,\prime\prime} \ . \ . \ . \ .\,\, b_n^{\,\prime\prime}]\).
Thus the matrix \(V\) corresponding to the operator \eqref{Chl}
takes the form
\begin{gather}
\label{HeR1}
V=||v_{jk}||_{1\leq{}j,k\leq{}n},\\
\intertext{where} %
\label{HeR2}
v_{jk}=\lambda\delta_{jk}-2(\overline{b_j^{\prime}}-\lambda\overline{
b_j^{\prime\prime}})\langle\,B^{\prime}-\overline{\lambda}B^{\prime\prime},
\,B^{\prime}-\overline{\lambda}B^{\prime\prime} \rangle^{-1}
(\lambda{}b_k^{\prime}-b_k^{\prime\prime})\,.
\end{gather}
and \(\lambda\) is such that
\[\lambda{}B^{\prime}(B^{\prime\prime})^{\ast}\geq{}0,\ \ |\lambda|=1\,.\]
\(\delta_{jk}\) is the Kronecker symbol. \hfill Q.E.D.
\begin{remark}%
\label{HaRo}%
In the case when the rows \(B^{\prime}\) and \(B^{\prime\prime}\)
are real, the matrix \(V\), \eqref{HeR1}-\eqref{HeR2}, is also real. In this
 case matrices of the form \eqref{HeR1}-\eqref{HeR2} are known as
\textsf{Householder reflection matrices}. Householder
 reflection matrices and the Householder Algorithm (which is based
 on matrices of this type) are widely used in numerical linear
 algebra. See \textup{\cite{Wil}, \cite{Str}, \cite{GolV} and
  \cite{Hou}}.
\end{remark}%
\begin{remark}
\label{Si}%
A unitary matrix \(V\) satisfying the condition
\eqref{OtoT} is not unique. The process of constructing such
matrices \eqref{HeR1}-\eqref{HeR2} is constructive.

We will apply \textup{Lemma \ref{Hh}} to the following special situation:
Let \(B^{\prime}\not=0\) be an arbitrary \(1\times{}n\)-column and
\(B^{\prime\prime}\) be of the special form
\(B^{\prime\prime}=[b^{\prime\prime}\,\,0_{1\times{}(n-1)}]\),
where \(b^{\prime\prime}>0\) and thus
\(b^{\prime\prime}=(B^{\prime}(B^{\prime})^{\ast})^{1/2}\).  For
these \(B^{\prime}, B^{\prime\prime}\), the first column of the
unitary matrix \(V\) satisfying \eqref{OtoT} is uniquely determined:
\[v_{j1}=\overline{b_{j}^{\prime}}(B^{\prime}(B^{\prime})^{\ast})^{-1/2},\,1\leq{}j\leq{}n.\]
The construction of the desired matrix \(V\) is thus reduced to
the following problem: Given the first column of an \(n\times{}n\)-matrix,
one needs to extend this column to a full \emph{unitary}
matrix. The Householder reflection procedure is one way of doing
this.
\end{remark}

We use the Householder reflection matrices to reduce an arbitrary
matrix to a Hessenberg matrix.\\[1.0ex]

\noindent%
 \textsf{PROOF of THEOREM \ref{RUHF}.} Let \(M=M^0\) and let
\(m_{j,k}^0\) be entries of the matrix \(M^0\):
\begin{equation}
\label{EnM0}%
 M^0=||m^0_{j,k}||_{0\leq{}j,k\leq{}n}
\end{equation}
Applying Lemma \ref{Hh}, we choose the unitary matrix
\(V_1\in\mathfrak{M}_{n,n}\) such that
\begin{equation}
\label{zr0}
 [m_{0,1}^0,\,m_{0,2}^0,\,\ldots\,,\,m_{0,n}^0]V_1=
[m_{0,1}^1,\,m_{0,2}^1,\,\ldots\,,\,m_{0,n}^1]\,,
\end{equation}
where
\begin{equation}%
\label{pr0}
 m_{0,1}^1\geq{}0,\,m_{0,k}^1=0,\,2\leq{}k\leq{}n.
\end{equation}%
(So that
\(m_{0,1}=\bigl([m_{0,1}^0,\,m_{0,2}^0,\,\ldots\,,\,m_{0,n}^0]
\cdot[m_{0,1}^0,\,m_{0,2}^0,\,\ldots\,,\,m_{0,n}^0]^\ast\bigr)^{1/2}\).)

\noindent%
\(V_1\) can be considered as an appropriate Householder
rotation, for instance. Let us consider the matrix
\begin{equation}%
\label{M1}
M^1=%
\begin{bmatrix}
 1&0_{1\times{}n}\\[0.7ex]
 0_{n\times{}1}&V_1^{\ast}
 \end{bmatrix}
 M^0%
 \begin{bmatrix}
 1&0_{1\times{}n}\\[0.7ex]
 0_{n\times{}1}&V_1
 \end{bmatrix}\,,
\end{equation}%
and let \(m_{j,k}^1\) denote the entries of the matrix \(M^1\):
\begin{equation}
\label{EnM1}%
 M^1=||m^1_{j,k}||_{0\leq{}j,k\leq{}n}
\end{equation}
Clearly,
\begin{equation}
\label{St1}%
 m_{0,0}^1= m_{0,0}^0\,.
\end{equation}
We continue this procedure inductively. We next turn to the
inductive step from \(l\) to \(l+1\).

Suppose that the matrices \(M^p\in\mathfrak{M}_{n+1,n+1}\) and
\(V_p\in\mathfrak{M}_{n-p+1,n-p+1}\) with \( 0\leq{}p\leq{}l\) are
already known and that the following condition for the entries of
the matrix \(M^p\),
\begin{equation}
\label{EnMp}%
 M^p=||m^p_{j,k}||_{0\leq{}j,k\leq{}n}\,,
\end{equation}
are satisfied:
\begin{equation}
\label{prp}%
 m^p_{j,j+1}\geq{} 0 ,\quad m^p_{j,k}=0, \ \
j+2\leq{}k\leq{}n\,, \ \ j=0,\,1,\,\ldots\,,\,p-1\,,
\end{equation}
The matrices \(V_p,\,1\leq{}p\leq{}l\) are unitary and we have
\begin{equation}%
\label{Mp}
M^p=%
\begin{bmatrix}
 I_p&0_{p\times{}(n-p+1)}\\[0.7ex]
 0_{(n-p+1)\times{}p}&V_p^{\ast}
 \end{bmatrix}
 M^{p-1}%
 \begin{bmatrix}
 I_p
&0_{p\times{}(n-1p+1)}\\[0.7ex]
 0_{(n-p+1)\times{}p}&V_p
 \end{bmatrix}
\end{equation}%
for every \(p: p\leq{}l\,.\)

 We choose the unitary \((n-l)\times(n-l)\)-matrix \(V_{l+1}\) such that
\begin{equation}
\label{zrl}
 [m_{l,l+1}^{l},\,m_{l,l+2}^{l},\,\ldots\,,\,m_{l,n}^{l}]V_{l+1}=
[m_{l,l+1}^{l+1},\,m_{l,l+2}^{l+1},\,\ldots\,,\,m_{l,n}^{l+1}]\,,
\end{equation}
where %
\begin{equation}%
\label{prl}
 m_{l,l+1}^{l+1}\geq{}0,\quad m_{l,k}^{l+1}=0,\ \ l+2\leq{}k\leq{}n.
\end{equation}%
Lemma \ref{Hh} ensures that this choice is possible. We then
define the matrix \(M^{l+1}\),
\begin{equation}
\label{EnMl}%
 M^{l+1}=||m^{l+1}_{j,k}||_{0\leq{}j,k\leq{}n}
\end{equation}
as
\begin{equation}%
\label{Ml}
M^{l+1}=%
\begin{bmatrix}
 I_{l+1}&0_{l+1\times{}(n-l)}\\[0.7ex]
 0_{(n-l)\times{}(l+1)}&V_{l+1}^{\ast}
 \end{bmatrix}
 M^{l}%
 \begin{bmatrix}
 I_{l+1}
&0_{(l+1)\times{}(n-l)}\\[0.7ex]
 0_{(n-l)\times{}(l+1)}&V_{l+1}
 \end{bmatrix}\,.
\end{equation}%
The entries  of the matrix \(M^{l+1}\) satisfy the condition
\begin{equation}
\label{prlM}%
 m^{l+1}_{j,j+1}\geq{} 0 ,\quad m^{l+1}_{j,k}=0, \ \
j=0,\,1,\,\ldots\,,\,l\,, \ \ j+2\leq{}k\leq{}n\,.
\end{equation}
For \(j=l\), condition \eqref{prl} holds in view of
\eqref{zrl} (Ensuring this was our goal in choosing the matrix
\(V_{l+1}\) as we did.)\\
For \(0\leq{}j\leq{l-1}\),
condition \eqref{prl} holds, because going from the matrix
\(M^{l}\) to the matrix \(M^{l+1}\) we do not change the rows with
indices \(j:\ 0\leq{}j\leq{}l-1\):
\begin{equation}%
\label{UCh}%
 m^{l+1}_{j,k}=m^{l}_{j,k},\quad 0\leq{}j\leq{}l-1,\
0\leq{}k\leq{}n\,.
\end{equation}%
The equality \eqref{UCh} holds, firstly because the identity
matrix of size \(l+1\) is the left upper corner of the block-matrix
\begin{math}
\Bigl[\begin{smallmatrix}
 I_{l+1}&0_{l+1\times{}(n-l)}\\[0.7ex]
 0_{(n-l)\times{}(l+1)}&V_{l+1}
 \end{smallmatrix}
\Bigr]\,
\end{math}
and secondly, because
\[ m^{l}_{j,k}=0\quad\forall \,j,k:\ 0\leq{}j\leq{}l-1,\,l+1\leq{}k\leq{}n
 \, \]
(The latter is a consequence of the induction hypothesis
\eqref{prp} for \(p=l-1\).) \\
The inductive process finishes when we construct the matrix
\(M_n=M^{l+1}\) for \(l=n-1\).

 The matrix \(V\) satisfying \eqref{RLH} appears as the
product
\begin{multline}%
\label{FC}%
V=V_1\cdot\\ %
\cdot\begin{bmatrix}
 I_1&0_{1\times{}(n-1)}\\[0.7ex]
 0_{(n-2)\times{}2}&V_2
 \end{bmatrix}\cdot
 \begin{bmatrix}
 I_2&0_{2\times{}(n-2)}\\[0.7ex]
 0_{(n-2)\times{}2}&V_3
 \end{bmatrix}\cdot\,\,\cdots\,\,\,\cdot
 \begin{bmatrix}
 I_{n-2}&0_{(n-2)\times{}2}\\[0.7ex]
 0_{2\times{}(n-2)}&V_{n-1}
 \end{bmatrix}\,.
 \end{multline}%
According to the above construction, the entries of the matrix
\(H=\) \\ \(=||h_{j,k}||_{0\leq{}j,k\leq{}n}\), \eqref{RLH}, satisfy:
\begin{equation}
\label{fe}%
 h_{j,k}=m_{j,k}^{j+1},\ \ j\leq{}k\leq{}n\,,
\end{equation}
and thus we have:
\begin{equation}
\label{feco}%
 h_{j,j+1}=m_{j,j+1}^{j+1}\geq{0},\ \ h_{j,k}=0,
 \ j+2\leq{}k\leq{}n\,,
\end{equation}
\ \hfill Q.E.D.

The reduction of matrices to the Hessenberg form is a tool often applied in
numerical linear algebra as a preliminary step for further numerical
algorithms. See \cite{Wil}, \cite{Str}, \cite{GolV} and other
sources in numerical linear algebra.

The Householder algorithm is implemented in the programming system
\texttt{MATLAB}. The \texttt{MATLAB} command  \texttt{H=hess(A)}
reduces the matrix \texttt{A} to the upper Hessenberg form
\texttt{H}.

In the next section we discuss the Schur algorithm for rational
inner functions in terms of the unitary colligation for the system
representation of this function. Reducing the colligation matrix
to the upper Hessenberg form is a preliminary step for further
developing the Schur algorithm in terms of system representations.
\begin{remark}
\label{Ref}%
 In \cite{KiNe}, the Householder algorithm and the Hessenberg form
 for unitary matrices are used to study the probability measures
 associated with finite Blaschke products via Cayley transform.
\end{remark}

\section{The Schur Algorithm in Terms of System Representations.
\label{SchA}} \setcounter{equation}{0}%
 We have now finished all necessary preparations and we are well positioned to
 present the Schur algorithm in terms of unitary colligations
 representing the appropriate functions.

Let \(s(z)\) be a  rational inner matrix-function of degree
\(n>0\) (\(s(z)\) is thus non-constant and \(|s_0|<1\), where
\(s_0=s(0)\)) and let
\[s_0(z)=s(z),\ \ \ s_k(z), \ k=1,\,2,\,\ldots\,,\,n\,,\] be the sequence
of rational inner functions constructed according to \eqref{kas}
(\(\deg{}s_k(z)=n-k\), so that \(s_n(z)=s_n\) is a unitary constant.)

Let
\begin{equation}
\label{SRst} s(z)=A+zB(I_n-zD)^{-1}C
\end{equation}
be the system representation of \(s(z)\), where
\begin{equation}%
\label{OSS} U=
\begin{bmatrix}
A&B\\[0.5ex]
C&D
\end{bmatrix}\,
\end{equation}%
is the matrix of the minimal unitary colligation representing
\(s(z)\):
\[ A\in\mathfrak{M}_{1\times{}1},\ B\in\mathfrak{M}_{1\times{}n},\
C\in\mathfrak{M}_{n\times{}1},\
D\in\mathfrak{M}_{n\times{}n}\,.\quad (\text{So}, A=s_0.)\]
We first reduce \(U\) to the lower Hessenberg form.
Let \(V\) be a unitary \(n\times{}n\)-matrix such that the matrix
\(U^0\) (also unitary):
\begin{equation}
 \label{U-0}
 U^0=
 \begin{bmatrix}
 1&0_{1\times{}n}\\
 0_{n\times{}1}&V^{\ast}
 \end{bmatrix}
 U%
\begin{bmatrix}
 1&0_{1\times{}n}\\
 0_{n\times{}1}&V
 \end{bmatrix}\,,\quad j=1,\,2\,,
 \end{equation}
 is an upper Hessenberg matrix.
The block entries of the matrix
 \begin{equation}%
\label{OSS0} U^0=
\begin{bmatrix}
A^0&B^0\\[0.5ex]
C^0&D^0
\end{bmatrix}\,
\end{equation}%
are:
\[ A^0\in\mathfrak{M}_{1\times{}1},\ B^0\in\mathfrak{M}_{1\times{}n},\
C^0\in\mathfrak{M}_{n\times{}1},\
D^0\in\mathfrak{M}_{n\times{}n}\,.\quad (\text{So}, A^0=A=s_0.)\]
The unitary colligations associated with the matrices \(U\) and
\(U^0\) are unitarily equivalent. The unitary colligation
associated with the unitary matrix \(U^0\) is therefore minimal and
represents the function \(s_0(z)=s(z)\):
\begin{equation}
\label{SR0} s_0(z)=A^0+zB^0(I_n-zD^0)^{-1}C^0\,.
\end{equation}

\noindent%
Inductively, we construct the sequence \(U^p,\
p=0,\,1,\,\ldots\,,\,n-1\) of unitary upper Hessenberg matrices
such that the unitary colligation associated with the matrix
\(U^p\) is minimal and represents the function \(s_p(z)\), which
appears in the \(p\)-th step of the Schur algorithm.

For \(p=0\), the representation in \eqref{SR0} holds. We consider
the step from \(p\) to \(p+1\).

Suppose that \(U^p\,,\ 0\leq{}p<(n-1)\) is a unitary lower
\(HL\)-non-singular \((n-p+1)\times(n-p+1)\) Hessenberg matrix
 with the block-matrix decomposition:
\begin{equation}%
\label{OSSp} U^p=
\begin{bmatrix}
A^p&B^p\\[0.5ex]
C^p&D^p
\end{bmatrix}\,,
\end{equation}%
where
\[ A^p\in\mathfrak{M}_{1\times{}1},\ B^p\in\mathfrak{M}_{1\times{}(n-p)},\
C^p\in\mathfrak{M}_{(n-p)\times{}1},\
D^p\in\mathfrak{M}_{(n-p)\times{}(n-p)}\,.\]%
 The unitary colligation associated with the matrix \(U^p\) is
 minimal and represents the function \(s_p(z)\), which appears in the
 \(p\)-th step of the Schur algorithm:
\begin{equation}
\label{SRp} s_p(z)=A^p+zB^p(I_{n-p}-zD^p)^{-1}C^p\,.
\end{equation}
Let
\begin{equation}%
\label{FBDp} C^p=
\begin{bmatrix}
C_1^p\\[0.5ex]
C_2^p
\end{bmatrix}\,
\qquad D^p=
\begin{bmatrix}
D_{11}^p&D_{12}^p\\[0.5ex]
D_{21}^p&D_{22}^p
\end{bmatrix}\,,
\end{equation}
be the more refined block matrix decomposition of the block-matrix
entries \(C^p\) and \(D^p\):
\[C_1^p\in\mathfrak{M}_{1\times{}1},\ C_2^p\in\mathfrak{M}_{(n-1-p)\times{}1},\]
\[ D_{11}\in\mathfrak{M}_{1\times{}1},\ D_{12}\in\mathfrak{M}_{1\times{}(n-1-p)},\
D_{21}\in\mathfrak{M}_{(n-1-p)\times{}1},\
D_{22}\in\mathfrak{M}_{(n-1-p)\times{}(n-1-p)}\,.\] %
Since \(U^p\) is an upper Hessenberg matrix and also an
\(HU\)-non-singular matrix, we have that \(B^p\not=0\).
Because \(U^p\) is also unitary, it follows that \(|A^p|<1\), i.e. that
\begin{equation}
\label{pScp} |s_p|<1 \ \ \text{where} \ \ s_p=s_p(0)\,.
\end{equation}
The row  \(B^p\) is of the form
\begin{equation}%
\label{Bp}
B^p=[(1-|s_p|^2)^{1/2},\,0_{1\times(n-p-1)}]%
 \end{equation}

\noindent We construct the \((n-p)\times(n-p)\)-matrix \(U^{p+1}\):
\begin{equation}
 \label{Uppo}
 U^{p+1}=
 \begin{bmatrix}
 \ A^{p+1} & B^{p+1}\\[0.5ex]
 \ C^{p+1}& D^{p+1}
 \end{bmatrix}\,,
 \end{equation}
 \[ A^{p+1}\in\mathfrak{M}_{1\times{}1}, B^{p+1}\in\mathfrak{M}_{1\times{}(n-p-1)},\
C^{p+1}\in\mathfrak{M}_{(n-p-1)\times{}1},
D^{p+1}\in\mathfrak{M}_{(n-p-1)\times{}(n-p-1)},\]%
 where
 \begin{equation}
 \label{DDMp}
 \begin{bmatrix}
 \ A^{p+1} & B^{p+1}\\[0.5ex]
 \ C^{p+1}& D^{p+1}
 \end{bmatrix}\stackrel{\tiny\text{def}}{=}
 \begin{bmatrix}
 (1-|s_p|^2)^{-1/2} C_1^p & D_{12}^p\\[0.5ex]
  (1-|s_p|^2)^{-1/2} C_2^p & D_{22}^p
 \end{bmatrix}\,.
\end{equation}
 To obtain the matrix \(U^p\) from \(U^{p+1}\), one should
delete the left column and the upper row of the matrix \(U^p\)
 and then recalculate the first column of the resulting matrix.
 The matrix \(U^{p+1}\) is then an upper Hessenberg
 matrix. The matrix \(U^{p+1}\) is \(HL\)-non-degenerate,
 because \(U^p\) is \(HL\)-non-degenerate and because the first superdiagonal of the
 matrix \(U^{p+1}\) is a subset of the first superdiagonal of the
 matrix \(U^{p}\).
  According to Theorem \ref{StScAl}
(which can be applied to the matrix \(U^{p}\) in view of \eqref{Bp}),
the matrix \(U^{p+1}\) is unitary and the unitary colligation
associated with \(U^{p+1}\) represents the function \(s_{p+1}(z)\)
appearing in the \(p+1\)-th step of the Schur algorithm:
\begin{equation}
\label{SRpp1}
s_{p+1}(z)=A^{p+1}+zB^{p+1}(I_{n-p-1}-zD^{p+1})^{-1}C^{p+1}\,.
\end{equation}
These considerations do not directly apply when \(p=n-1\). In this case,
there is no room for \(B^{n},\,C^{n},\,D^{n}\). However, we can
construct `part' of the matrix \eqref{DDMp}:
\begin{equation}%
\label{LShPa}%
 A^n=(1-|s_{n-1}|^2)^{-1/2}C_1^{n-1}\,.
\end{equation}%
(See Remark \ref{LShP}.) The \(1\times{}1\)-matrix \(A^n\) is
unitary, hence it is a unitary constant. Clearly, \(A^n=s_n\),
where \(s_n\) is the \(n\)-th Schur parameter. This completes the
description of the Schur algorithm for inner rational matrix-functions
in terms of system representations. \hfill Q.E.D.
\begin{remark}
\label{DSf} It is particularly easy to determine the sequence
\(\{D^p\}_{p=0,\,1\,\ldots\,,\,n}\) of matrices representing
the inner operators of the unitary colligations associated with
the colligation matrices \(U^p\). The matrix \(D^p\) makes up the
\((n-p)\times(n-p)\) lower-right corner of the matrix \(D^0\).
The inner rational matrix-function \(s(z)\) is the ratio of two
polynomials:
\begin{equation}%
s_p(z)=c_p\frac{z^{n-p}\overline{\chi_p(1/\overline{z})}}{\chi_p(z)},\quad
\deg{}\chi_p(z)=n-p\,, \ \chi_p(0)=1\, \ \ |c_p|=1.
\end{equation}%
Clearly,%
\begin{subequations}
\label{Corn}
 \begin{equation}
 \label{Corn_1}
\chi_p(z)=\det(I_{n-p}-zD^p)\,,\ \
z^{n-p}\overline{\chi_p(1/\overline{z})}=\det\big(zI_{n-p}-(D^p)^{\ast}\big)\,,
\end{equation}thus
\begin{equation}
 \label{Corn_2}
s_p(z)=c_p\det\Big(\big(zI_{n-p}-(D^p)^{\ast}\big)\big(I_{n-p}-zD^p\big)^{-1}\big)\,.
\end{equation}
\end{subequations}
\end{remark}
\section{An Expression for the Colligation Matrix in Terms of the Schur Parameters.
\label{ESchP}} \setcounter{equation}{0}%
Let \(s(z)\) be a rational inner matrix-function of degree \(n\).
Let \(s_p(z),\,p=0,\,1,\,\ldots\,,\,n\) be the sequence of
rational inner functions produced by the Schur algorithm from
the function \(s(z)\), as described in \eqref{kas}, \(\deg
s_p(z)=n-p\). Let \(U^p\) , \eqref{OSSp}, be the colligation
matrix of the minimal unitary colligation, which yields the
system representation \eqref{SRp} of the function \(s_p\). Among
all unitary \((n-p+1)\times(n-p+1)\)-matrices representing the
function \(s_p\) we choose a lower Hessenberg matrix \(U^p\).
Such a matrix \(U^p\) exists and is unique.

The equality \eqref{SE}, where \(U^p\) is taken as the
matrix \(U^0\) and \(U^{p+1}\) is taken as the matrix
\begin{math}
\bigl[\begin{smallmatrix}
 \alpha & \beta\\[0.5ex]
 \gamma & \delta
\end{smallmatrix}\bigr]\,
\end{math}
takes the form
\begin{multline*}
U^p=
\begin{bmatrix}
1  & \ 0_{1\times{}1} \ & \hspace*{2.0ex}0_{1\times(n-1-p)} \\[0.7ex]
0_{1\times{}1}   &  &\\[-0.5ex]
& &\hspace*{-4.0ex}U^{p+1}\hspace*{6.0ex} \\[-1.0ex]
 0_{(n-1-p)\times{}1}   & &
\end{bmatrix}\cdot\\
\cdot\begin{bmatrix}
s_p & \ (1-|s_p|^2)^{1/2} \ &\  0_{1\times{}(n-1-p)} \\[0.5ex]
(1-|s_p|^2)^{1/2}& \!\!-\overline{s_p} \ &  0_{1\times(n-1-p)}\\[0.8ex]
0_{(n-1-p)\times{}1}   &  0_{(n-1-p)\times{}1} &
I_{(n-1-p)\times(n-1-p)}
\end{bmatrix}
\,,\\
p=0,\,1,\, \ldots\,,\,n-2.
\end{multline*}
The latter formula can be rewritten in the equivalent but more
convenient form:\\
\begin{multline}
\label{SEp}%
\begin{bmatrix}
I_{p}& 0_{p\times(n-p+1)}\\[0.7ex] 0_{(n-p+1)\times{}p}& U^p
\end{bmatrix}=
\begin{bmatrix}
I_{p+1}& 0_{(p+1)\times(n-p)}\\[0.7ex] 0_{(n-p)\times{}(p+1)}&
U^{p+1}
\end{bmatrix}\cdot\\[2.5ex]%
\cdot
\begin{bmatrix}
I_{p} &0_{p\times{2}} & \ \ \ 0_{p\times(n-1-p)} \\[3.0ex] %
 0_{2\times{}p} & \ \mbox{\(\begin{matrix}s_{p}&(1-|s_{p}|^2)^{1/2}\\[1.0ex] %
 (1-|s_{p}|^2)^{1/2}&-\overline{s_{p}}
 \end{matrix}\)} \ & 0_{2\times(n-p-1)}\\[5.0ex] %
 \ \ \  0_{(n-p-1)\times{}p} &0_{(n-p-1)\times{}2} & I_{n-1-p}
\end{bmatrix}\,.\\[1.0ex]
0\leq{}p\leq{}n-1\,.
\end{multline}
 For \(p=0\), the matrix on the left-hand side of \eqref{SEp} takes
 the form  \(\begin{bmatrix}U^0\end{bmatrix}\). For \(p=n-1\),
 the matrix \(U^{p+1}\) takes the form \(U^n=s_n\) and
 the second factor on the right-hand side of \eqref{SEp} takes
 the form
\begin{equation*}
\begin{bmatrix}
I_{n-2} &0_{(n-2)\times{2}}  \\[3.0ex] %
 0_{2\times(n-2)} & \ \mbox{\(\begin{matrix}s_{n-1}&(1-|s_{pn-1}|^2)^{1/2}\\[1.0ex] %
 (1-|s_{n-1}|^2)^{1/2}&-\overline{s_{n-1}}
 \end{matrix}\)} %
\end{bmatrix}\,.
\end{equation*}
From \eqref{SEp} it follows that
\begin{multline}
\label{PrE}%
 U^0=\\
\Pprod\limits_{0\leq{}p\leq{}n-1}^{\mbox{\Curvearrowleft}}
\begin{bmatrix}
I_{p} &0_{p\times{2}} & \ \ \ 0_{p\times(n-1-p)} \\[3.0ex] %
 0_{2\times{}p} & \ \mbox{\(\begin{matrix}s_{p}&(1-|s_{p}|^2)^{1/2}\\[1.0ex] %
 (1-|s_{p}|^2)^{1/2}&-\overline{s_{p}}
 \end{matrix}\)} \ & 0_{2\times(n-p-1)}\\[5.0ex] %
 \ \ \  0_{(n-p-1)\times{}p} &0_{(n-p-1)\times{}2} & I_{n-1-p}
\end{bmatrix}\,.
\end{multline}
Multiplying the matrices in \eqref{PrE}, we obtain an expression
for the entries of the  matrix \(U^0\), which gives us the system
representation of the function \(s(z)\) in terms of the
Schur parameters of \(s(z)\):
\begin{equation}%
\label{ESchPa}%
U^0=||u^0_{j,k}||_{0\leq{j,k}\leq{}n}\,,
\end{equation}%
where
\begin{equation}%
\label{EfE}%
 u^0_{j,k}=
\begin{cases}
\phantom{-}s_0,& j=0,\hspace*{10.0ex} k=0,\\
 \phantom{-}s_j\,\Delta_{j-1}\Delta_{j-2}\cdot\,\,\cdots\,\,\cdot{}\Delta_1\Delta_0\,,&
 1\leq{}j\leq{}n\,,\hspace*{5.0ex} k=0\,,\\
 -s_j\,\Delta_{j-1}\Delta_{j-2}\cdot\,\,\cdots\,\,\cdot{}\Delta_{k}\,%
 \overline{s_{k-1}},&
 1\leq{}j\leq{}n\,,\hspace*{5.0ex} 1\leq{}k\leq{j},\,\\
 \phantom{-}\Delta_j\,, &\,0\leq{}j\leq{}n-1\,,\,\,k=j+1\,,\\
\phantom{-}0\,,&\, 0\leq{}j<n-1,\,\,\,j+1<k \leq{}n\,,
\end{cases}
\end{equation}%
with
\begin{equation}%
\label{Def}%
 \Delta_j=(1-|s_j|^2)^{1/2}\,.
\end{equation}%
One can, in the same way, obtain expressions for the matrices
\(U^j\) of the unitary colligations representing the functions
\(s_j, \,1\leq{}j\leq{}n\).

It should be mentioned that a matrix of the form \eqref{ESchPa},
\eqref{EfE} appeared in the paper \cite[formula
\((66^\prime)\)]{Ger} and was then rediscovered a number times.
See \cite{Grg}, \cite{Con1}, \cite[Section 2.5]{Con2}, \cite{Tep},
\cite[Chapter 4]{Sim}, \cite[Theorem 2.17]{Dub}.

\section{On Work Related to System Theoretic Interpretations
of the Schur Algorithm\label{RelWork}}
\setcounter{equation}{0}%

In this section we discuss the connections between the present work
and other work relating to the Schur algorithm as expressed in terms of
system realizations. In particular, we discuss the results presented
in \cite{AADL} and in \cite{KiNe}.

The paper \cite{AADL} deals with functions of the class %
\(\boldsymbol{S_\kappa}\), i.e. with the functions \(s(z)\)
meromorphic in the unit disc and possessing the properties:
\begin{enumerate}
\item[1).]
For every \(N\) and for all points
\(z_1,\,\ldots\,,\,z_N\in\mathbb{D}\) which are holomorphicity points
for \(s\), the matrix \(\|K(z_p,z_q)\|_{1\leq{}p,q\leq{}N}\),
\(K(z,\zeta)=\frac{1-s(z)\overline{s(\zeta)}}{1-z\overline{\zeta}}\),
does not have more than \(\kappa\) negative squares.
\item[2).] There exists an \(N\) and points \(z_1,\,\ldots\,,\,z_N\in\mathbb{D}\)
such that this matrix has precisely \(\kappa\) negative squares.
\end{enumerate}
One of the goals of the paper \cite{AADL} is to discuss the Schur
algorithm for functions from the class \(\boldsymbol{S_\kappa}\)
in terms of system realizations. In particular, the results of
\cite{AADL} are applicable to the special case%
\footnote{
\(\boldsymbol{S_0}\) is the class of contractive functions
holomorphic in the unit disc.} %
\(\kappa=0\), in which they can be simplified.
In our considerations on the algebraic structure of a step of the Schur
algorithm we will, for the sake of simplicity, restrict ourselves to
finite-dimensional systems, which correspond to rational inner functions
(of, say, degree \(n\)). We now describe the relevant result from
\cite{AADL}, adopting the notation used there (to make
the comparison with the results presented in our paper easier). In
\cite{AADL} the function \(s(z)\) is given by
\begin{equation}
\label{1}%
 s(z)=s_0+zB(I-zD)^{-1}C,
\end{equation}
 where \(B\), \(C\), \(D\) are entries of a unitary matrix \(U\),
\begin{equation}
\label{2}%
 U=
\begin{bmatrix}
s_0&B\\[0.5ex]
C&D
\end{bmatrix}\,,
\quad
B\in\mathfrak{M}_{1\times{}n},\,C\in\mathfrak{M}_{n\times{}1},\,
D\in\mathfrak{M}_{n\times{}n}
\end{equation}
  It is not \textit{explicitly} assumed from the very beginning that
the entry \(B\) of the matrix \(U\) has the special form (7.36).
The matrix \(U\) appears as the matrix \(V\), (1.2), in \cite{AADL}.
Our notation corresponds to that of \cite{AADL} as follows:
The objects, which appear as
\(\gamma\), \(v\), \(u\), \(T\) in formula (1.2) of \cite{AADL} are
\(s_0\), \(B^\ast\), \(C\), \(D\) in our formulas (7.34)-(7.35).
The state space which is denoted by \(\mathcal{K}\) in (1.2) of
\cite{AADL} is the space \(\mathcal{H}=\mathfrak{M}_{n\times{}1} \ \
(=\mathbb{C}^n)\) in our paper.

Let \(s_1(z)\) be the Schur transform of the function \(s(z)\),
\begin{equation}
\label{3}
s_1(z)=\frac{1}{z}\cdot\frac{s(z)-s_0}{1-s(z)\overline{s_0}}\,,\qquad
s_0=s(0)
\end{equation}
(or (7.33) in our paper). According to \cite{AADL}, \(s_1(z)\) is
representable in the form
\begin{equation}
\label{4}%
 s_1(z)=\alpha+z\beta(I-z\delta)^{-1}\gamma\,,
\end{equation}
 with
\begin{alignat}{2}
\alpha&=\frac{1}{1-|s_0|^2}\,BC,\,\quad&%
\beta&=\frac{1}{\sqrt{1-|s_0|^2}}\,BDP,\notag\\
 \gamma&=\frac{1}{\sqrt{1-|s_0|^2}}\,PC,\quad&\delta&=\,PDP\,,
 \label{5}
\end{alignat}
where \(P\) is the matrix of the orthogonal projector onto the
orthogonal complement of the vector \(B^\ast\) in \(\mathcal{H}\),
i.e. %
 \begin{equation*}
 P\in\mathfrak{M}_{n\times{}n},\,%
\textup{rank}\,P=n-1,\,P^2=P,\,P=P^\ast,\, BP=0\,.
\end{equation*}
(Formulas \eqref{5} are the formulas for the entries of the matrix
\(V_1\) which appear on page 11 of \cite{AADL}.)
 If we would like to represent the image space \(P\mathcal{H}\)
as
the space %
\(\mathcal{H}_1=\mathfrak{M}_{(n-1)\times{}1} \ \
(=\mathbb{C}^{n-1})\),
\(\mathcal{H}=\mathbb{C}\oplus\mathcal{H}_1\), that is, if we would
like the matrix of the projector \(P\) to be of the form %
\begin{math}
P=%
\begin{bmatrix}
0&0\\
0&I_{n-1}
\end{bmatrix}\,,
\end{math}
then we have to replace the original matrix \(U\) with the
matrix %
\begin{equation}
\label{6} %
U^0=\begin{bmatrix}
 1 &0_{n\times{}n}\ \\[0.5ex]
0_{n\times{}n}    & V^\ast\end{bmatrix}\,%
U\,%
\begin{bmatrix}
 1 &0_{n\times{}n}\ \\[0.5ex]
0_{n\times{}n}   & V
\end{bmatrix}\,,\quad%
U^0=\begin{bmatrix}
 s_0 &B_0 \\[0.5ex]
C_0    & D_0\end{bmatrix}
\end{equation}
where \(V\in\mathfrak{M}_{n\times{}n}\) is a unitary matrix such
that %
\begin{equation}
\label{7} %
V^\ast{}PV=\begin{bmatrix}
0&0\\
0&I_{n-1}
\end{bmatrix}.%
\end{equation}
 The condition \(BP=0\) implies the condition
\(B_0\begin{bmatrix}
 1 &0\\[0.5ex]
0   & I_{n-1}
\end{bmatrix}=0\).
The last equality means that \(B_0\) is of the form %
\(B_0=\begin{bmatrix} b&0_{1\times(n-1 )}
\end{bmatrix}\). %
 Since the matrix \(U^0\) is unitary, we have
\(|s_0|^2+B_0B_0^\ast=1\). Therefore, \(B_0\) must be of the form
\begin{equation*}
B_0=
\begin{bmatrix} \delta(1-|s_0|^2)^{1/2}
&0_{1\times(n-1 )}
\end{bmatrix},
\end{equation*}
where \(\delta\) is a unimodular complex number. The unitary
matrix \(V\) from \eqref{6} is not unique: In this case, the
degrees of freedom are clear, when we consider the replacement
\begin{math}
V\to{}V\cdot\begin{bmatrix}
 \varepsilon &0\\[0.5ex]
0   & v
\end{bmatrix},
\end{math}
where \(\varepsilon\) is an arbitrary unimodular complex number
and \(v,\ v\in\mathfrak{M}_{(n-1)\times{}(n-1)}\) are
unitary matrices. Choosing the number \(\varepsilon\) appropriately,
we can ensure that \(B_0\) is of the form
\begin{equation}
\label{8}%
 B_0=
\begin{bmatrix} (1-|s_0|^2)^{1/2}
&0_{1\times(n-1 )}
\end{bmatrix}.
\end{equation}
 Let us
decompose the matrices \(C_0,\,D_0\), which appear as the entries
of the matrix \(U_0\) from \eqref{5}:
\begin{equation}%
\label{9}%
 C_0=
\begin{bmatrix}
c_1\\[0.5ex]
c_2
\end{bmatrix}\,,
\qquad D_0=
\begin{bmatrix}
d_{11}&d_{12}\\[0.5ex]
d_{21}&d_{22}
\end{bmatrix}\,,
\end{equation}
\[c_1\in\mathfrak{M}_{1\times{}1},\ c_2\in\mathfrak{M}_{(n-1)\times{}1},\]
\[ D_{11}\in\mathfrak{M}_{1\times{}1},\ D_{12}\in\mathfrak{M}_{1\times{}(n-1)},\
D_{21}\in\mathfrak{M}_{(n-1)\times{}1},\
D_{22}\in\mathfrak{M}_{(n-1)\times{}(n-1)}\,.\] The equalities
\eqref{5} (where \(B,\,C,\,D\) are replaced by
\(B_0,\,C_0,\,D_0\)) now take the form
\begin{alignat}{2}
\alpha&=\frac{1}{\sqrt{1-|s_0|^2}}c_1,\,\quad&%
\beta&=d_{12},\notag\\
 \gamma&=\frac{1}{\sqrt{1-|s_0|^2}}\,c_2,\quad&\delta&=d_{22}\,,
 \label{10}
\end{alignat}
Thus, the matrix
\begin{equation}
\label{11}%
 U^1=\begin{bmatrix}
\alpha&\beta\\[0.5ex]
\gamma&\delta
\end{bmatrix},
\end{equation}
from \cite{AADL}, whose entries appear in the representation \eqref{4}
of the function \(s_1(z)\) is the same matrix which appears in our
Theorem 7.2 as the matrix (7.38). (The matrix \(V_1\) from page 11
of \cite{AADL} can be considered as a \textit{coordinate-free}
expression for the colligation matrix representing the function
\(s_1(z)\).)
 The difference between our work
and the work \cite{AADL} is not in the results but in the methods. The
reason for choosing the
expression for the colligation matrix \(U_1\) given in \cite{AADL} is
not fully explained. The facts that the matrix \(U_1\) is unitary and
that the matrix \(U_1\) represents the Schur transform \(s_1\) of
the function \(s\) are obtained as the result of a long chain of
\textit{formal} calculations. These calculations come across as
somewhat contrived and do not serve to further our understanding
of the subject at hand.

The state system approach is much more transparent. The fact that
the matrix \(U_1\) is unitary is an immediate consequence of our
formula (7.6). The fact that the matrix \(U_1\) represents the
function \(s_1\) is a consequence of the interpretation of the
linear fractional transform (6.6)-(6.7) in terms of the Redheffer
coupling of the appropriate colligation.

The paper \cite{KiNe} can also be considered as relevant to our paper.
In \cite{KiNe} the system representation of Schur functions is not considered
at all. Nevertheless, in this work the Householder algorithm is used to
calculate the sequence of numbers, which can be identified with the Schur
parameters of the rational inner function naturally related to the appropriate
unitary matrix. Namely, given a unitary matrix
\(U \in \mathfrak{M}_{(n+1)\times(n+1)}\), the measure \(\mu\)
on the unit circle is related to \(U\) in the following way:
\(\mu(dt) = \left( E(dt)e_1, \, e_1 \right)\), where \(E(dt)\) is the
sprectral measure of the matrix \(U\) and \(e_1 = (1, 0,  \ldots , 0)^T\),
\(e_1 \in \mathfrak{M}_{(n+1) \times 1}\). It is assumed that \(e_1\) is a
cyclic vector of \(U\). The following equality holds:
\begin{equation}
\label{12}%
 e_1^* \frac{I + zU}{I - zU} e_1
 = \displaystyle \int_\mathbb{T}{\frac{1 + zt}{1 - zt}}\mu(dt)
\end{equation}
The measure \(\mu\) generates the (finite) sequence of polynomials orthogonal
on the unit circle. These orthogonal polynomials (\(\Phi_k\) is monic of degree
\(k\)) satisfy the recurrence relations
\begin{equation}
\label{13}%
  \Phi_{k + 1}(z) = z \Phi_{k}(z) - \overline{s}_k \Phi^*_{k}(z)
\end{equation}
\begin{equation}
\label{14}%
  \, \Phi^*_{k + 1}(z) = z \Phi^*_{k}(z) - s_k z \Phi_{k}(z)
\end{equation}
where \(s_k\), \, \(k = 0, 1, \ldots, n\) are some recurrence coefficients.
There are many different names for these coefficients. Recently dubbed
`Verblunsky parameters' by Barry Simon in \cite{Sim}.
On the other hand, the function in \eqref{12}, which we denote by \(p(z)\)
is holomorphic in the unit disc \(\mathbb{D}\) and has the following
properties.
\begin{equation*}
  p(0)  =   1,  \qquad  p(z) + \overline{p(z)} \geq 0 \quad (z \in \mathbb{D}).
\end{equation*}
Therefore \(p(z)\) is representable in the form
\begin{equation}
\label{15}%
  p(z)  =   \frac{1 + z s(z)}{1 - z s(z)},
\end{equation}
where \(s(z)\) is a function holomorphic and contractive in \(\mathbb{D}\).
Ya. L. Geronimus established that the Verblunsky coefficients \(s(z)\) in
the recurrence relations \eqref{13} - \eqref{14} are also the Schur parameters
of the functions \(s(z)\), which appear in \eqref{15}. From \eqref{12} and
\eqref{15} it follows that
\begin{equation}
\label{16}%
   e_1^* \frac{I + zU}{I - zU} e_1  =   \frac{1 + z s(z)}{1 - z s(z)}.
\end{equation}
In Lemma 3.2 of \cite{KiNe}, the following method for finding Schur (=Verblunsky)
parameters was proposed: First, the given unitary matrix \(U\) should be converted
to Hessenberg form:
\begin{equation}
\label{17}%
   U^0
   =    \begin{bmatrix}
            1                               &   0_{1 \times n}      \\
            0_{n \times 1}  &   V^*
        \end{bmatrix}
        U
        \begin{bmatrix}
            1                               &   0_{1 \times n}      \\
            0_{n \times 1}  &   V
        \end{bmatrix},
\end{equation}
In \cite{KiNe} it is claimed that the entries of the (lower
Hessenberg) matrix \(U^0\) are of the form
\eqref{ESchPa}-\eqref{Def}, from which the Schur-Verblunsky
parameters \(s_k\) can be found. However, it follows from
\eqref{16} that
\begin{equation}
\label{18}%
   s(z) = A + s B(I - zD)^{-1}C ,
\end{equation}
where
\begin{equation}
\label{19}%
   U
   =    \begin{bmatrix}
            A                               &   B       \\
            C                               &   D
        \end{bmatrix}
\end{equation}
\begin{equation*}
A    \in \mathfrak{M}_{1 \times 1}, \quad%
B\in\mathfrak{M}_{1\times n},\quad%
C\in\mathfrak{M}_{n\times 1},\quad%
D\in\mathfrak{M}_{n\times n\,.}
\end{equation*}
Thus the formula \eqref{18} can be interpreted as the system
representation of the function \(s(z)\). The formula \eqref{18},
where \(s(z)\) is defined by \eqref{16} from \(U\) was unfamiliar
to us, but we do not think that this formula is new.

 In his forthcoming
paper [Arl] Yu. M. Arlinskii studied a related question for
operator-valued Schur functions \(\Theta\) acting between
separable Hilbert spaces. These investigations correspond to the
operator generalization of the classical Schur algorithm which is
due to Constantinescu (see Section 1.3 in \cite{BC}.) Yu. M.
Arlinskii presents a construction of conservative and simple
realizations of the Schur algorithm iterates \(\Theta_n\) of
\(\Theta\) by means of the conservative and simple realization of
\(\Theta\).

\vspace{20pt}

\begin{center}
\textbf{\large Appendix:\\[1.0ex]
 System Realizations of Inner Rational Functions\label{App}.}
\end{center}
\setcounter{equation}{0}%
\renewcommand{\theequation}{\textup{A}.\arabic{equation}}

We prove that every complex-valued (i.e. scalar) inner
rational function of degree \(n\) can be represented as the
characteristic function of the minimal unitary colligation
associated with some  unitary matrix \(U\in\mathfrak{M}_{(n+1)\times{}(n+1)}\).
 Let us denote a given rational inner function by \(S\).
 The operator colligation whose characteristic
 function is \(S\) will be constructed as the `left shift'
 operator in the appropriate space of analytic functions
 constructed from \(S\). A similar construction appears in a paper
 by B.SzNagy-C.Foias. See \cite[Chapter VI]{SzNFo}. The
 construction of B.SzNagy-C.Foias was adapted to unitary
 colligations in \cite{BrSv2}.

\vspace{12pt}

\textbf{1.} \textbf{The space \(\boldsymbol{K_S}\)}. The most important
part of our construction is the Hilbert space \(K_S\) of rational functions.
 We consider \(S\) as a function defined
on the unit circle \(\mathbb{T}\), i.e. \(S:\,\mathbb{T}\to\mathbb{T}\).
As usual,
\begin{center}
\(L^2=\{x:\mathbb{T}\to\mathbb{C}, \|x\|<\infty\}\), \qquad where
\(\|x\|^2=\langle{}x\,,\,{}x\rangle\) \\ and
\[\langle{}x\,,\,{}y\rangle=\int\limits_{\mathbb{T}}x(t)\,\overline{y(t)}\,m(dt)\,,\]
\end{center}
\vspace{-2.0ex}
 \(m(dt)\) is the normalized Lebesgue measure on
\(\mathbb{T}\). Let \(H^2_+\) and \(H^2_-\) be the Hardy subspaces
of the space \(L^2\):

\[H^2_+=\{x\in L^2: \langle{}x(t)\,,\,{}t^k\rangle=0,\ \ k=-1,-2,\,\ldots\,\}\,.\]
\[H^2_-=\{x\in L^2: \langle{}x(t)\,,\,{}t^k\rangle=0,\ \ k=0,\,1,\,2,\,\ldots\,\ \}\,.\]
Clearly,
\[L^2=H^2_+\oplus H^2_-\,.\]
It is also convenient to consider the functions from \(H^2_+\) and from \(H^2_-\)
as functions holomorphic in \(\mathbb{D}\) and in \(\mathbb{D}^{-}\), respectively.
In particular, the evaluation \(f\to f(0)\) is defined for every \(f\) in \(H^2_+\)
and \(f(\infty)=0\) for every \(f\) in \(H^2_-\).

The space \(K_S\) is defined as
\begin{equation}
\label{DMSp}%
 K_S=H^2_+\ominus\, S\,H^2_+\,,
\end{equation}
where \(S\,H^2_+= \{S(t)\,h(t): \,h \in H^2_+\}\).
Another description of the space \(K_S\) is:
\begin{equation}
\label{ODMSp}%
K_S=\{x\in{}L^2:\,x\in H^2_+,\, xS^{-1}\in H^2_-\,\}.
\end{equation}
It can be shown that the space \(K_S\) consists of rational
functions whose poles are contained in the set of poles of the function
\(S\) and that \(\dim K_S=\deg S\). If all zeros \(z_k\) of
\(S\) are simple (see \eqref{BP}), then the space \(K_S\) is
generated by the functions \(\{(1-t\overline{z_k})^{-1}\}_{1\leq k
\leq n}\). If \(S\) has non-simple zeros, the modification of this
statement is clear. The space \(K_S\) is a reproducing kernel Hilbert space.
If \(f\in K_S\), then
\begin{equation}
\label{RKI}%
f(z)=\langle f(t), K(t,z)\rangle\,,
\end{equation}
where the reproducing kernel \(K(t,z)\) is:
\begin{equation}
\label{RKE}
K(t,z)=\frac{1-S(t)\overline{S(z)}}{1-t\overline{z}}\,.
\end{equation}

\vspace{12pt}

\textbf{2.} \textbf{The left shift operator.} The left shift operator
\(T\) is defined as
\begin{equation}
\label{DLSO}%
 T(f)(t)=(f(t)-f(0)e(t) )\cdot t^{-1}  \quad
\text{for} \ f\in H^2_+\,.
\end{equation}
where
\begin{equation}
\label{Uni}%
e(t) = 1\ \  \ \ \forall t \in \mathbb{T}\,.
\end{equation}
 This operator is contractive:
\begin{equation}
\label{CoSO}
\|Tf\|^2=\|f\|^2-|f(0)|^2 \quad \forall f\in H^2_+\,.
\end{equation}
The space \(K_S\), considered as a subspace of \(H^2_+\), is an invariant subspace
of the left shift operator \(T\). This is evident from the description \eqref{ODMSp}
of the space \(K_S\).

\vspace{12pt}

\textbf{3.} \textbf{The construction of the unitary colligation
\(U\).} The unitary colligation \((\mathcal{E},\
\mathcal{H},\,U)\) (see Definition \ref{DOC}) is defined as
follows: Let the state space \(\mathcal{H}\) be the space \(K_S\) and
let the principal operator \(D\) be the left shift operator \(T\),
\eqref{DLSO}, restricted to \(K_S\):
\begin{equation}
\label{IOMC} \mathcal{H}=K_S, \quad Df(t)=(f(t)-f(0)e(t))\cdot
t^{-1} \ \ \forall f\in \mathcal{H}\,.
\end{equation}
The equality \(\|Df\|^2+|f(0)|^2=\|f\|^2\), together with the
requirement that the colligation operator \(U\),
\eqref{UnCol}-\eqref{CODB}, be unitary, prompts us to define the
exterior space \(\mathcal{E}\) and the channel operator
\(B:\,\mathcal{H}\to\mathcal{E}\) as follows:

Let \(\mathcal{E}\) be a one-dimensional Hilbert space which is
identified with the \emph{vector space} \(\mathbb{C}\) over the
\emph{field \(\mathbb{C}\) of scalars}.
 We choose the number \(\beta=1\) as a basis \emph{vector} in
\(\mathbb{C}\) and will denote this basis vector by \Bbo.
Every number \(\varepsilon\in\mathbb{C}\), considered as an
element \(\Bbe\) of the \emph{vector space} \(\mathbb{C}\), can be
presented as \(\Bbe=\varepsilon\Bbo\), where the factor
 in front of \(\Bbo\) is the same \emph{number} \(\varepsilon\),
 but considered as an element of the \emph{field of scalars} \(\mathbb{C}\).

The channel operator \(B\) is:
\begin{equation}
\label{DChO1} (Bf)(t)=f(0)\Bbo,\quad \forall f\in\mathcal{H}\,.
\end{equation}
Equation \eqref{CoSO} ensures that
\begin{equation*}
\|Bf\|_{\mathcal{E}}^2+\|Df\|_{\mathcal{H}}^2=\|f\|^2_{\mathcal{H}}\,,
 \quad \forall f\in\mathcal{H}\,.
\end{equation*}
\(f(0)\), which appears in \eqref{IOMC} and \eqref{DChO1},
can be represented using the reproducing kernel \eqref{RKI}-\eqref{RKE}.
Let
\begin{equation}
\label{ChVe1}%
 k(t)=1-s(t)\,\overline{s(0)},\quad (\,=K(t,0)\,)\,.
\end{equation}
Then
\begin{equation}
\label{ChOp1}%
 Bf=\langle\,f,k\,\rangle\,\Bbo\,.
\end{equation}
The operator \(A: \mathcal{E}\to\mathcal{E}\)
(as is the case for every operator in \(\mathcal{E}:\dim \mathcal{E}=1\))
is of the form
\[A\Bbe=\alpha\langle\,\Bbe,\Bbo\,\rangle\,\Bbo,\ \ \Bbe\in\mathcal{E},\]
where \(\alpha\in\mathbb{C}\). Since the vector \(\Bbo\), which
generates \(\mathcal{E}\), is orthogonal to \(\mathcal{H}\) in the
orthogonal sum \(\mathcal{E}\oplus\mathcal{H}\), the unitary property of
\(U\) implies that
\begin{equation}
\label{OREU}%
\langle{}Bf\,,\,A\Bbo{}\rangle+\langle{}Df\,,C\Bbo\,\rangle=0\quad
\forall \,f\in\mathcal{H}\,.
\end{equation}
Therefore
\begin{equation*}
\alpha\langle{}Bf\,,\,\Bbo{}\rangle+\langle{}Df\,,C\Bbo\,\rangle=0\quad
\forall \,f\in\mathcal{H}\,.
\end{equation*}
Let us denote
\begin{equation*}%
C\Bbo=l,\,\quad l\in\mathcal{H}\,.
\end{equation*}%
Equation \eqref{OREU} means that
\begin{equation*}%
\overline{\alpha}\langle\,f,k\,\rangle+\langle\,Df,l\,\rangle=0\,,\quad {}\forall f\in\mathcal{H}\,.
\end{equation*}%
Thus, one should take
\begin{equation}%
\label{ECV}%
 l=-\alpha{}(D^{\ast})^{-1}k,
 \end{equation}
 where \(D^{\ast}\) is the adjoint to the operator \(D\), with respect
 to the scalar product \(\langle\,,\,\rangle\).

 We now look to determine the operator \(D^{\ast}\). The equality
 \begin{equation*}\langle\,Df,g\,\rangle=
 \langle\,f,D^{\ast}g\,\rangle\ \ \forall{}f,g\in\mathcal{H}
 \end{equation*}
 means that
 \begin{equation*}%
 \langle\,(f(t)-f(0)e(t))\,t^{-1},g(t)\,\rangle=
\langle\,f(t),(D^{\ast}g)(t)\,\rangle\,.
\end{equation*}
The last equality implies that
\begin{equation}
(D^{\ast}g)(t)=P(tg(t)),\quad \forall g\in{}K_S
\end{equation}
where \(P\) is the orthogonal projector from \(L^2\) onto \(K_S\).
Clearly, \[\langle\,h(t),S(t)\,\rangle=0 \ \ \forall h\in{}K_S\,, \]
and
\[tg(t)-\langle\,tg(t),S(t)\,\rangle S(t)  \in K_S \ \forall g\in K_S\,.\]
Therefore,
\[P(tg(t))=tg(t)-\langle\,tg(t),S(t)\,\rangle S(t)\  \ \forall g\in K_S\,,\]
that is
\begin{equation}%
\label{EAPO}%
(D^{\ast}g)(t)=tg(t)-\langle\,tg(t),S(t)\,
\rangle S(t)\  \ \forall g\in\mathcal{H}\,.
\end{equation}%
From \eqref{ECV}, we obtain
\begin{equation*}%
l(t)=\frac{\alpha}{S(0)}\frac{S(t)-S(0)e(t)}{t}\,.
\end{equation*}%
\(\|Ae\|^2+\|Ce\|^2=1\) gives us
\(|\alpha|=|S(0)|\). We choose \[\alpha=S(0)\] (Later we see that
this is the only possible choice for \(\alpha\).) We set
\begin{equation}%
\label{ECVe}%
l(t)=\frac{S(t)-S(0)e(t)}{t}\,.
\end{equation}%
(In intermediate steps we assumed that \(S(0)\not=0\), but this does not appear
in the final expression \eqref{ECVe} for \(l(t)\).)
Thus,
\begin{multline}%
\label{FEC}
A\Bbe=S(0)\langle\,\Bbe,\Bbo\,\rangle\,\Bbo,\quad
 Bf=\langle\,f,k\,\rangle\,\Bbo,\quad
 C\Bbe=\langle\,\Bbe,\Bbo\,\rangle\,l,\\[1.0ex]
(Df)(t)=(\,f-\langle\,f,k\,\rangle e(t)\,)\,t^{-1}\quad \forall
\Bbe\in\mathcal{E},\,f\in\mathcal{H}\,\,.
\end{multline}
or
\begin{multline}%
\label{FEO} A\Bbe=\varepsilon{}S(0)\Bbo\,,\quad Bf=f(0)\Bbo,
\quad (C\Bbe))=\varepsilon{}l(t),\\[1.0ex]
(Df)(t)=\big(\,f(t)-f(0)e(t)\,\big)\,t^{-1}\,,\quad \forall
\Bbe=\varepsilon\Bbo\in\mathcal{E},\,f\in\mathcal{H}\,.
\end{multline}

From \eqref{FEO} it follows that the block-operator
 \begin{equation}%
 \label{BOFC}%
 U=\begin{bmatrix}
 A&B\\C&D
\end{bmatrix}
 \end{equation}%
  is unitary (After the block \(D\) was chosen, the other blocks
 \(A,\,B,\,C\) were chosen to ensure that \(U\) be a unitary operator.)
 The characteristic function \(S_U(z)\),
 \begin{equation}%
 \label{CFMC}
 S_U(z)=A+zB(I-zD)^{-1}C\,,
 \end{equation}
 of the colligation \(U\) coincides with the original rational inner function
 \(S(z)\). This can be checked by direct calculation of \(S_U(z)\) using the expression
 \eqref{FEO} for blocks of the colligation operator \(U\). The expression for the
 operator \((I-zD)^{-1}\), which is needed for this calculation, is
 \begin{equation}
 \label{ReMO}
 \big((I-zD)^{-1}f\big)(t)=\frac{tf(t)-zf(z)}{t-z},\quad \forall f\in\mathcal{H}\,.
 \end{equation}
 In what follows we also need the expression for the operator \((I-zD^{\ast})^{-1}\):
 \begin{equation}
 \label{ReMOa}
 \big((I-zD^{\ast})^{-1}f\big)(t)=
 \frac{f(t)-(fS^{-1})(z^{-1})S(t)}{1-tz},\quad \forall f\in\mathcal{H}\,.
\end{equation}
(Since the function \(fS^{-1}\) belongs to \(H^2_-\)\,, the evaluation
\(fS^{-1}\to{}(fS^{-1})(z^{-1})\) is defined for \(z\in\mathbb{D}\).)

Choosing an orthogonal basis in the (\(n\)-dimensional) Hilbert space \(K_S\),
we realize that the unitary operator \(U\), \eqref{FEO}-\eqref{BOFC},
originally constructed as an operator
acting in \textit{a functional} space, is {a matrix} operator acting in
\(\mathbb{C}^{n+1}=\mathbb{C}\oplus\mathbb{C}^{n}\).
\begin{definition}
\label{DeMoCo}%
 The colligation \eqref{FEO}-\eqref{BOFC} is called the
 \textsf{model unitary colligation constructed from the rational inner function
 \(S\).}
\end{definition}

\textbf{4.} \textbf{Minimality of the model unitary colligation
\(\boldsymbol{(\mathcal{E},\ \mathcal{H},\,U)}\).} We look to prove
that the model colligation \(U\), \eqref{FEO}-\eqref{BOFC}, is
controllable and observable. In view of the expression for the
channel operator \(C\) (one-dimensional), controllability of \(U\)
can be formulated as
follows:
\begin{equation}
\label{Contr}%
 \textit{The set of vectors
\(\{(I-zD)^{-1}l\}_{z\in\mathbb{D}}\) \ %
  generates the space \(K_S\).}
\end{equation}
From \eqref{ECVe} and \eqref{ReMO} it follows that
\begin{equation*}
\big((I-zD)^{-1}l\big)(t)=\frac{S(t)-S(z)}{t-z}\,.
\end{equation*}
Let \(f\in{}L^2\) be such that
  \begin{equation}
 \label{OrtCont}%
\int\limits_{\mathbb{T}}\frac{S(t)-S(z)}{t-z}\,\overline{f(t)}\,m(dt)=0\
\ \forall\,z\in\mathbb{D}
 \end{equation}
 If \(f\in{}H^2_+\), then
  \(\int\limits_{\mathbb{T}}\frac{\overline{f(t)}}{t-z}\,m(dt)=0\
  \forall{}z\in\mathbb{D}\),
  hence,
 \(\int\limits_{\mathbb{T}}\frac{S(t)\overline{f(t)}}{t-z}\,\,m(dt)=0\ %
 \forall\,z\in\mathbb{D}\). The last equality implies that
\(\overline{f(t)S^{-1}(t)}\in{}H^2_-\,.\)
 (Here we use that \(S(t)=\overline{S^{-1}(t)}\) for \(t\in\mathbb{T}\).)
 If also \(f(t)S^{-1}(t) \in{}H^2_-\,,\) then
  \(f(t)S^{-1}(t)\equiv{}0 \) and \(f\equiv{}0\).
Therefore, if the condition \eqref{OrtCont}
 holds for some \(f\in{}K_S\), then \(f\equiv{}0\).
   Controllability of the colligation
  \((\mathcal{E},\ \mathcal{H},\,U)\) is thus proved.

      Observability of this colligation can be proved analogously.
 According to \eqref{ChOp1},
 \(B^{\ast}f=\langle\,f,e\,\rangle\,k\). Therefore the observability
 criterion is reduced to the statement:
\begin{equation}
\label{Observ}%
 \textit{The set of vectors
\(\{(I-zD^{\ast})^{-1}k\}_{z\in\mathbb{D}}\) \ %
  generates the space \(K_S\).}
\end{equation}
Using expressions \eqref{ReMOa} and \eqref{ChVe1}, we obtain:
\begin{equation*}
\big((I-zD^{\ast})^{-1}k\big)(t)=\frac{1-S(t)\,S^{-1}(z^{-1})}{1-tz}\,.
\end{equation*}
Let \(f\in{}L^2\) is such that
 \begin{equation}
 \label{OrtObs}%
\int\limits_{\mathbb{T}}\frac{1-S(t)S^{-1}(z^{-1})}{1-tz}\,\overline{f(t)}\,m(dt)=0\
\ \forall\,z\in\mathbb{D}\,.
 \end{equation}
 If \(f(t)S^{-1}(t)\in{}H^2_-\),
 then
 \(\int\limits_{\mathbb{T}}\frac{S(t)\overline{f(t)}}{1-tz}=0\),
 hence
 \(\int\limits_{\mathbb{T}}\frac{\overline{f(t)}}{1-tz}m(dt)=0\
 \forall{}z\in\mathbb{D}\) and \(f(t)\in{}H^2_-\). If also
 \(f\in{}H^2_+\), then \(f\equiv{}0\). Therefore if the condition \eqref{OrtObs}
 holds for some \(f\in{}K_S\), then \(f\equiv{}0\).

\vspace{12pt}

 \textbf{5.} \textbf{Uniqueness of simple realization}.
 The uniqueness of the minimal realization is, in fact, a version
 of a result by M.S.\,Livshitz, which, in the language of
 M.S.Livshitz, claims that the characteristic function uniquely determines
 (up to unitary equivalence) the operator colligation without complementary
 component.

 Let
 \begin{math}
  U_1=\begin{bmatrix}
 A_1&B_1\\C_1&D_1
\end{bmatrix}
 \end{math}
 and
\begin{math}
  U_2=\begin{bmatrix}
 A_2&B_2\\C_2&D_2
 \end{bmatrix}
 \end{math}
 be two unitary matrices divided into blocks,
\begin{multline}
\label{BlDiv}%
    A_j\in\mathfrak{M}_{1\times{}1}, \ \
    B_j\in\mathfrak{M}_{1\times{}n_j}, \
\  C_j\in\mathfrak{M}_{n_j\times{}1}, \ \
D_j\in\mathfrak{M}_{n_j\times{}n_j},\\[0.5ex]
\text{where} \ n_j,\,j=1,\,2, \ \text{are natural numbers.}
\end{multline}
We \textsf{do not} assume that \(n_1=n_2\).

Let
\[ S_i(z)=A_i+B_i(I-zD_i)^{-1}C_i,\quad i=1,2,\]
be the characteristic functions of the unitary colligations
associated with the matrices \(U_1\) and \(U_2\) respectively.
Suppose that
\begin{itemize}
\item[\textup{1.}] The characteristic functions are equal.
\begin{equation}
\label{EqChF}%
S_1(z)\equiv{}S_2(z).\end{equation}%
\item[\textup{2.}] Each of the matrices \(U_1\) and \(U_2\) is
simple in the sense of Definition \ref{DMM}.
\end{itemize}

We prove that under these assumptions the matrices \(U_1\) and
\(U_2\) are equivalent in the sense of Definition \ref{DEUM}, in
particular, \(n_1=n_2\).

To prove this, we have to first of all construct a unitary
mapping \(V\) of the space \(\mathbb{C}^{n_2}\) onto the space
\(\mathbb{C}^{n_1}\). Assume that the matrices \(U_1\) and \(U_2\)
are simple. Let us consider the vectors
\(f_j^k,\,g_j^l\in\mathbb{C}^{n_j}\,(=\mathfrak{M}_{n_j\times{}1}),\,\,j=1,2,\,
0\leq{}k,\,l\,\):
\begin{equation}%
\label{DeGeVe}%
f_j^k=D_j^kC_j,\ \ g_j^l=(D_j^{\ast})^lB_j^{\ast},\,j=1,\,2, \
0\leq{}k,l\,.
\end{equation}%
By the assumption, for each \(j=1,\,2\), the vectors
\(f_j^k,\,g_j^k,\,0\leq{}k\leq\max(n_1,\,n_2)\) generate the space
\(\mathbb{C}^{n_j}\). The equality \eqref{EqChF} implies
\begin{equation}
\label{ECTe} A_1=A_2\,,
\end{equation}
and the equalities
 \begin{equation}
 \label{FEvC}%
B_1D_1^{p}C_1=B_2D_2^pC_2,\,\,0\leq{}p ,
\end{equation}%
 or
\begin{equation*}
 B_1{D_1^{l}D_1^{k}}C_1=B_2D_2^lD_2^{k}C_2,
  \ \ 0\leq k,l\,.
\end{equation*}
The latter equalities can be interpreted as
\begin{subequations}
\label{ECP}
\begin{equation}%
\label{ECPa}%
\langle\,f_1^k,\,g_1^l\rangle_{\mathbb{C}^{n_1}}=
\langle\,f_2^k,\,g_2^l\rangle_{\mathbb{C}^{n_2}}\,,\quad \forall
k,l: 0\leq{}k,\,0\leq{}l\,.
\end{equation}%
Moreover, the equalities \eqref{EqChF} imply that
\begin{equation*}
1-S_1^{\ast}(\zeta)S_1(z)\equiv{}1-S_2^{\ast}(\zeta)S_2(z),\quad
1-S_1(z)S_1^{\ast}(\zeta)\equiv{}1-S_2(z)S_2^{\ast}(\zeta)\,.
\end{equation*}
In view of \eqref{Id1}, \eqref{Id2}, the latter equalities imply
that
\begin{multline*}
 C_1^{\ast}(D_1^{\ast})^qD_1^pC_1=C_2^{\ast}(D_2^{\ast})^qD_2^pC_2,\
 \ \text{and} \ \  B_1D_1^q(D_1^{\ast})^pC_1^{\ast}=
 B_2D_2^q(D_2^{\ast})^pC_2^{\ast},\
\\[0.5ex] 0\leq p,\,q\,.
\end{multline*}
This can, in turn, be interpreted as
\begin{multline}%
\label{ECPb}%
\langle\,f_1^p,\,f_1^q\rangle_{\mathbb{C}^{n_1}}=
\langle\,f_2^p,\,f_2^q\rangle_{\mathbb{C}^{n_2}}\,,\quad
\text{and}\quad \langle\,g_1^p,\,g_1^q\rangle_{\mathbb{C}^{n_1}}=
\langle\,g_2^p,\,g_2^q\rangle_{\mathbb{C}^{n_2}}\,,\\[0.5ex]
\forall p,q: 0\leq{}p,\,0\leq{}q\,.
\end{multline}%
\end{subequations}
From \eqref{ECP} it follows that for arbitrary
\(\alpha_k,\beta_l\) (such that only finitely many of them differ
from zero),
\begin{equation}
\label{NoEq}
\big\|\sum\alpha_kf^k_1+\sum\beta_lg^l_1\big\|_{\mathbb{C}^{n_1}}=
\big\|\sum\alpha_kf^k_2+\sum\beta_lg^l_2\big\|_{\mathbb{C}^{n_2}}\,.
\end{equation}
Let us define the operator
\(V:\,\mathbb{C}^{n_2}\to\mathbb{C}^{n_1}\) first as
\begin{subequations}
\label{DeV}
\begin{equation}
\label{DeV1}%
 Vf^k_2=f^k_1, \ \ Vg^l_2=g^l_1, \ \ \forall\ \
k\geq{}0,\,l\geq{}0,
\end{equation}
and then extend this operator by linearity to all vector columns
\(h\in\mathbb{C}^{n_2}\) representable as a finite linear
combination of the form \(h=\sum\alpha_kf^k_2+\sum\beta_lg^l_2\).
Thus,
\begin{equation}%
\label{DeV2}%
V\big(\sum\alpha_kf^k_2+\sum\beta_lg^l_2\big)=\sum\alpha_kf^k_1+\sum\beta_lg^l_1\,.
\end{equation}
\end{subequations}
If some \(h\in\mathbb{C}^{n_2}\) admits two different
representations, say
\begin{equation*}%
h=\sum\alpha_k^{\prime}f^k_2+\sum\beta_l^{\prime}g^l_2, \ \
\text{and} \ \
h=\sum\alpha_k^{\prime\prime}f^k_2+\sum\beta_l^{\prime\prime}g^l_2\,,
\end{equation*}
then \(Vh\) also admits two different representations:
\begin{equation*}%
Vh=\sum\alpha_k^{\prime}f^k_1+\sum\beta_l^{\prime}g^l_1, \ \
\text{and} \ \
Vh=\sum\alpha_k^{\prime\prime}f^k_1+\sum\beta_l^{\prime\prime}g^l_1\,.
\end{equation*}
However, since \(\sum\alpha_kf^k_2+\sum\beta_lg^l_2=0\), where
\(\alpha_k=\alpha_k^{\prime\prime}-\alpha_k^{\prime}\),
\(\beta_l=\beta_l^{\prime\prime}-\beta_l^{\prime}\), the equality
\eqref{NoEq} implies that
\(\sum\alpha_kf^k_1+\sum\beta_lg^l_1=0\), i.e.
\[\sum\alpha_k^{\prime}f^k_1+\sum\beta_l^{\prime}g^l_1=
\sum\alpha_k^{\prime\prime}f^k_1+\sum\beta_l^{\prime\prime}g^l_1\,.\]
The definition \eqref{DeV} of \(V\) is thus non-contradictory.

The operator \(V\) is defined on the linear hull of all vectors
\(\lbrace{}f^k_2,g^l_2{}\rbrace_{k,l}\) and isometrically maps its
definition domain onto the linear hull of all vectors
\(\lbrace{}f^k_1,g^l_1{}\rbrace_{k,l}\)\,. If both the matrices
\(U^2,\,U^1\) are simple, then these linear hulls are the whole
spaces \(\mathbb{C}^{n_2}\) and \(\mathbb{C}^{n_1}\), respectively.
In this case \(n_1=n_2\,(\stackrel{\text{\tiny def}}{=}n)\) and
\begin{equation}
\label{UnV}%
 V^{\ast}V=I_n,\quad VV^{\ast}=I_n
\end{equation}

 We now prove the intertwining relation
\begin{equation*}
  \begin{bmatrix}
 A_1&B_1\\C_1&D_1
\end{bmatrix}
\begin{bmatrix}
 1&0\\0&V
\end{bmatrix}
=\begin{bmatrix}
 1&0\\0&V
\end{bmatrix}\begin{bmatrix}
 A_2&B_2\\C_2&D_2
 \end{bmatrix}\,,
 \end{equation*}
which can be rewritten as follows:
\begin{equation}
\label{IntW1}%
 C_1=VC_2 ,\quad B_2=B_1V ,
\end{equation}
 and
\begin{equation}
\label{IntW2}%
VD_2=D_1V.
\end{equation}
The first of the equalities \eqref{IntW1} corresponds to the first of the
equalities \eqref{DeGeVe} for \(k=0\) (See \eqref{DeV1} for
\(k=0\).) The second of the equalities \eqref{IntW1} relates to the second
of equality in \eqref{DeGeVe} for \(l=0\) (See \eqref{UnV}.)

To check the splitting relation \eqref{IntW2}, it is enough to
check that
\begin{subequations}
\begin{equation}
\label{ChSpRe1} VD_2f^k_2=D_1Vf^k_2\quad \textup{for} \ \ \forall
k\geq{}0,
\end{equation}
and
\begin{equation}
\label{ChSpRe2} VD_2g^l_2=D_1Vg^l_2\quad \textup{for} \ \ \forall
l\geq{}0,
\end{equation}
\end{subequations}
The  equality \eqref{ChSpRe1} is an obvious consequence of the
definitions of the operator \(V\) and vectors \(f_j^k\). Indeed,
\(D_2f^k_2=f^{k+1}_{\,2},\,Vf^{k+1}_{\,2}=f^{k+1}_{\,1}\).
On the other hand,
\(Vf^{k}_{\,2}=f^{k}_{\,1},\,D_1f^{k}_{\,1}=f^{k+1}_{\,1}\).
Therefore, \eqref{ChSpRe1} holds.

 Our approach to checking the condition \eqref{ChSpRe2} will
be different in the two cases \(l=0\) and \(l>0\). For \(l=0\) the
equality \eqref{ChSpRe2} takes the form
\(VD_2B_2^{\ast}=D_1VB_2^{\ast}\). \eqref{DeV} for \(l=0\) means
that \(VB_2^{\ast}=B_1^{\ast}\), so we should check that
\(VD_2B_2^{\ast}=D_1B_1^{\ast}\,.\) Since
\(D_jB_j^{\ast}=-C_jA_j^{\ast},\,j=1,2\), the last equality is
equivalent to \(VC_2A_2^{\ast}=C_1A_1^{\ast}\). The
latter equation is a consequence of the first of the equalities
\eqref{IntW1}. Thus, \eqref{ChSpRe2} holds for \(l=0\).

We check condition \eqref{ChSpRe2} for \(l>0\). Since the
matrices \(U_j\) are unitary, we have that \(D_jD_j^{\ast}=I-C_jC_j^{\ast}\).
Thus, \eqref{ChSpRe2} is equivalent to
\begin{equation*}V(I-C_2C_2^{\ast})(D_2^{\ast})^{l-1}B^{\ast}_2
=(I-C_1C_1^{\ast})(D_1^{\ast})^{l-1}B^{\ast}_1\,.
\end{equation*}
This equation is a consequence of the following three equalities:
\begin{subequations}
\label{CoThE}
\begin{equation}
\label{CoThE1}
V(D_2^{\ast})^{l-1}B^{\ast}_2=(D_1^{\ast})^{l-1}B^{\ast}_1\,,
\end{equation}
\begin{equation}
\label{CoThE2} VC_2=C_1\,,
\end{equation}
and
\begin{equation}
\label{CoThE3}
C_2^{\ast}(D_2^{\ast})^{l-1}B^{\ast}_2=C_1^{\ast}(D_1^{\ast})^{l-1}B^{\ast}_1\,.
\end{equation}
\end{subequations}
\eqref{CoThE1} holds, because it can be written as
\(Vg^{l-1}_2=g^{l-1}_1\), which is part of the definition
\eqref{DeV} of the operator \(V\). \eqref{CoThE2} has already been
checked: This is the first of the relations \eqref{IntW1}.
\eqref{CoThE3} is the same as \eqref{FEvC} for \(p=l-1\).
The condition \eqref{ChSpRe2} has also been checked for \(l>0\).

\vspace{12pt}

 \textbf{6.} \textbf{Simple realization is minimal.}
 Let \(U\in\mathfrak{M}_{(1+n)\times(1+n)}\) be a simple unitary
 matrix. The matrix \(U\) is then minimal. Indeed, let \(S(z)\) be
 the characteristic function of the unitary colligation associated
 with \(U\). \(S(z)\) is a rational inner function, \(\deg{}S\leq{}n\).
 Let \((\mathbb{C},K_S,T)\) be the model colligation
 constructed from this \(S\). We established that the model unitary
 colligation is minimal (in particular, simple) and that its characteristic function is
 the function \(S\), from which it was constructed.
 Both colligations (the original colligation and the model colligation)
 have the same characteristic function and both are simple.
 Hence, these colligations are equivalent. Since the model
 colligation is minimal, the original colligation is also minimal.

\begin{center}
\textbf{\large Acknowledgements.}
\end{center}
\hspace{20pt} We thank Professors I. Gohberg and M. A. Kaashoek for useful remarks
and suggestions on the history and scope of the state space method
and its applications.

Moreover, we thank Professor D. Alpay for drawing our attention to
the paper \cite{AADL}.

We thank Armin Rahn for his careful reading of the manuscript and his
help in improving the English in this paper.

\end{document}